# THE MONODROMY MATRIX FOR A FAMILY OF ALMOST PERIODIC SCHRÖDINGER EQUATIONS IN THE ADIABATIC CASE

ALEXANDER FEDOTOV AND FRÉDÉRIC KLOPP

ABSTRACT. This work is devoted to the study of a family of almost periodic one-dimensional Schrödinger equations. We define a monodromy matrix for this family. We study the asymptotic behavior of this matrix in the adiabatic case. Therefore, we develop a complex WKB method for adiabatic perturbations of periodic Schrödinger equations. At last, the study of the monodromy matrix enables us to get some spectral results for the initial family of almost periodic equations.

RÉSUMÉ. Ce travail est consacré à l'étude d'une famille d'équations de Schrödinger quasi-périodiques en dimension 1. Nous définissons une matrice de monodromie pour cette famille. Nous étudions le comportement asymptotique de cette matrice dans le cas adiabatique. A cette fin, nous développons une version de la méthode WKB complexe pour l'équation de Schrödinger périodique avec une perturbation adiabatique. Enfin, l'étude de la matrice de monodromie nous permet d'obtenir des résultats spectraux pour la famille d'équations quasi-périodiques initiale.

## 1. INTRODUCTION AND STATEMENT OF THE MAIN RESULTS

In this paper, we study the family of Schrödinger equations

$$(1.1) \quad -\frac{d^2}{dx^2}\psi(x) + (V(x-\phi) + W(\varepsilon x))\psi(x) = E\psi(x), \quad x \in \mathbb{R}.$$

Here
- $\phi \in \mathbb{R}$ is a parameter numbering the equations,
- $V$ and $W$ are two periodic functions satisfying

$$(1.2) \quad V(x+1) = V(x), \quad W(\xi + 2\pi) = W(\xi), \quad x, \xi \in \mathbb{R},$$

- $\varepsilon$ is a fixed positive number.

The ratio of the periods of the potential in (1.1) is equal to $2\pi/\varepsilon$. If $2\pi/\varepsilon \notin \mathbb{Q}$, the potential $V(x-\phi) + W(\varepsilon x)$ is almost periodic.

Our ultimate goal is to understand the spectral properties of the family (1.1) assuming that $V$ is a real-valued measurable bounded function of $x \in \mathbb{R}$ and that $W$ is analytic in $\xi \in \mathbb{C}$.

1991 *Mathematics Subject Classification.* 34E05, 34E20, 34L05.

This work was begun while A.F. was staying at the Université Paris 12-Val de Marne. A.F. thanks Prof. Guillopé for the kind invitation and hospitality during his stay in Paris. A.F. and F.K. acknowledges the support of the Fields Institute (Toronto) where this work was completed. This work was also partially supported by the INTASS grant RFBR 96-0414.

It is well known that for one-dimensional periodic Schrödinger equation, the key to the understanding of its spectral properties is the monodromy matrix. It is not possible to define this object for a single almost periodic differential equation. In this paper, we will show that, for a family of almost periodic differential equations of type (1.1), a monodromy matrix can be defined and used for its spectral analysis.

Of course, a thorough understanding of the structure of monodromy matrix is crucial. Therefore, in this paper, we will need to take the adiabatic limit $\varepsilon \to 0$. As we shall see, the monodromy matrix always enables us to reduce the study of (1.1) to the study of a difference equation. In the adiabatic case, this equation takes a standard form and can be effectively analyzed.

## 1.1. Reduction to a difference equation.
We begin by introducing several natural notions related to the family of equations (1.1).

1.1.1. For any fixed $\phi$, let $\psi_{1,2}(x, \phi)$ be two linearly independent solutions of equation (1.1). We say that they form a *consistent basis* if

$$\psi_{1,2}(x, \phi + 1) = \psi_{1,2}(x, \phi), \quad \forall x, \phi, \tag{1.3}$$

and if their Wronskian defined by

$$w(\psi_1, \psi_2) = \psi_1' \psi_2 - \psi_2' \psi_1. \tag{1.4}$$

is independent of $\phi$. For convenience, we shall assume that

$$w(\psi_1, \psi_2) = 1, \quad \forall \phi. \tag{1.5}$$

Remark that the solutions verifying the Cauchy conditions

$$\psi_1(0, \phi) = 0, \quad \psi_1'(0, \phi) = 1,$$
$$\psi_2(0, \phi) = 1, \quad \psi_2'(0, \phi) = 0$$

form a consistent basis. However, as we shall see in the next sections, this choice of a consistent basis is not the most natural, and, for a while, we do not choose any particular construction of the solutions $\psi_1$ and $\psi_2$.

1.1.2. In view of (1.2), the functions $\psi_{1,2}(x + 2\pi/\varepsilon, \phi + 2\pi/\varepsilon)$ are also two solutions of equation (1.1). Therefore, one can write

$$\Psi(x + 2\pi/\varepsilon, \phi + 2\pi/\varepsilon) = M(\phi) \Psi(x, \phi), \tag{1.6}$$

where $\Psi$ is the vector

$$\Psi^T(x, \phi) = (\psi_1(x, \phi), \psi_2(x, \phi)), \tag{1.7}$$

where $^T$ is the symbol of transposition, and $M(\phi)$ is a matrix with coefficients independent of $x$. We call this matrix *the monodromy matrix* corresponding to the consistent basis $\psi_{1,2}$.

Note that this is a natural generalization of the canonical definition of the monodromy matrix: if, in the Schrödinger equation (1.1), $V \equiv 0$, then the potential is periodic, the basis solutions $\psi_{1,2}$ are independent of $\phi$, and formula (1.6) defines the standard monodromy matrix corresponding to two linearly independent solutions of the periodic Schrödinger equation, see [21].

It follows from the definition that, for any consistent basis,
$$\det M(\phi) \equiv 1, \quad M(\phi+1) = M(\phi), \quad \forall \phi.$$

1.1.3. Set

(1.8) $$h = \frac{2\pi}{\varepsilon} \bmod 1.$$

Let $M$ be a monodromy matrix corresponding to a consistent basis $\psi_{1,2}$. In Section 6, we shall see that the spectral analysis of (1.1) can be reduced to the investigation of the solutions of the equation

(1.9) $$\chi(\phi + h) = M(\phi)\chi(\phi), \quad \forall x \in \mathbb{R}.$$

Assume that $\chi$ is a *fundamental* solution i.e. that

(1.10) $$\det \chi(\phi) = 1, \quad \forall \phi \in \mathbb{R}.$$

One can easily obtain that

(1.11) $$\Psi(x + 2\pi n/\varepsilon, \phi) = \chi(\phi)\chi^{-1}(\phi - nh)\Psi(x, \phi - nh).$$

It is well known that the spectral properties of one dimensional Schrödinger equations can be completely described in terms of the behavior of its solutions as $x \to \pm\infty$ (see [12]). Remembering that $\Psi$ is 1-periodic in $\phi$, we see that relation (1.11) reduces the analysis of the behavior for $x \to \pm\infty$ of the solutions of the input differential equation (1.1) to the analysis of solutions of the difference equation (1.9).

The passage from equation (1.1) to (1.9) is done in the same spirit as the monodromization transformation [8]. For example, this latter transformation reduces the spectral analysis of a difference Schrödinger equation to the analysis of a difference equation of the form (1.9); in this case too, the matrix defining the resulting equation is the monodromy matrix corresponding to one of the fundamental solutions of the initial difference equation.

Here, as in the case of the monodromization, one can define Bloch solutions for the almost periodic family (1.1), and show that these Bloch solutions can be constructed in terms of Bloch solutions of the difference equation (1.9) (for the definition of the latter, see [8]). These ideas proved to be effective for the spectral analysis of the Harper equation [6]. We will discuss and use these ideas in section 6.

1.2. **Complex WKB method.** To make the above ideas work, one has to find an effective way to control the monodromy matrices for family of equations (1.1). For this, it suffices to have a good enough control of some consistent basis solutions. In this paper, we construct and study certain consistent basis for small values of $\varepsilon$ using the ideas of the complex WKB method. The complex constructions appear to be inevitable since some of the coefficients of the monodromy matrices are exponentially small when $\varepsilon \to 0$. The results and the constructions developed first for (1.1) are very general, and hopefully, they may be applied to study a wide class of problems for adiabatically perturbed one-dimensional periodic Schrödinger equations. So, we formulate the results in a more or less general form without assuming that the function $W$ is periodic. We also note that our assumptions on the regularity of $V$ and $W$ can be essentially weakened.

1.2.1. We begin with some preliminary remarks concerning the standard version of the complex WKB method developed for ordinary differential equations of the form (see, for example [11])

$$(1.12) \qquad -\varepsilon^2 \frac{d^2}{dz^2}\psi(z) + U(z)\psi(z) = E\psi(z), \quad z \in \mathbb{C}.$$

Here, $U$ is an analytic function of $z$, and $\varepsilon$ is the standard small semi-classical parameter. For (1.12), the complex WKB method consists in the description of certain *canonical* domains in the complex plane $z \in \mathbb{C}$, and in constructing analytic solutions of (1.12) having standard asymptotic behavior for $\varepsilon \to 0$ on a given canonical domain. For a given canonical domain $K$, one constructs two linearly independent solutions $\psi_\pm$ with the asymptotics

$$\psi_\pm(z) = \frac{1}{\sqrt{p(z)}} e^{\pm \frac{i}{\varepsilon} \int_{z_0}^z p(z)\,dz + o(1)}, \quad z \in K, \quad \varepsilon \to 0.$$

Here $z_0$ is a fixed point in $K$ and $p$ is the complex momentum defined by the relation

$$p^2(z) + U(z) = E$$

(i.e. by the symbol of the equation (1.12)). The complex momentum is the main analytic object of the WKB constructions. The canonical domain can be characterized as a domain such that, through any point of it, one can draw a smooth curve along which the function $\operatorname{Im} \int_{z_0}^z p(z)\,dz$ is monotonically increasing. For precise definitions, we again refer to [11].

1.2.2. For the family of equations (1.1) the role of the complex momentum is played by the function $\kappa(\varphi)$ satisfying the relation

$$(1.13) \qquad \mathcal{E}(\kappa) + W(\varphi) = E,$$

where $\mathcal{E}(k)$ is the dispersion law of the periodic Schrödinger equation

$$(1.14) \qquad -\frac{d^2}{du^2}\psi(x) + V(x) = E\psi(x), \quad x \in \mathbb{R}.$$

In other words,

$$(1.15) \qquad \kappa(\varphi) = k(E - W(\varphi)),$$

where $k(E)$ is the Bloch quasi-momentum of (1.14), see subsection 2.2. As $p(z)$, the complex momentum $\kappa$ is a multi-valued analytic function. Its branching points are related to the branching points of the quasi-momentum.

1.2.3. The simplest asymptotic behavior of analytic solutions of (1.1) is predicted by the already classical ansatz suggested by V. S. Buslaev. In his papers (see, for example, [1], [2], [3], [4]), V. S. Buslaev has studied the asymptotic behavior along the real line $x \in \mathbb{R}$ of solutions of the equation

$$(1.16) \qquad -\frac{d^2}{dx^2}\psi(x) + U(x, \varepsilon x)\psi(x) = E\psi(x), \quad x \in \mathbb{R}.$$

In this equation, $U(x, \xi)$ is a sufficiently smooth function periodic in $x$. Between the branching points of the complex momentum, the leading term of the asymptotics of two linearly independent solutions has the form

$$(1.17) \qquad \psi_\pm \sim A_\pm(\varepsilon x) e^{\pm \frac{i}{\varepsilon} \int_{\varepsilon x_0}^{\varepsilon x} \kappa(\xi)\,d\xi} p_\pm(x, \varepsilon x), \quad \varepsilon \to 0.$$

Here, $p_{\pm}(x, \xi)$ are the periodic component of the Bloch solutions (see subsection 3.1 and 3.3) of the equation

$$-\frac{d^2}{dx^2}\psi(x) + U(x, \xi)\psi(x) = E\psi(x), \quad x \in \mathbb{R}, \tag{1.18}$$

where the parameter $\xi \in \mathbb{R}$ is frozen and $A_{\pm}(\xi)$ are two "normalization constants" described in terms of the periodic components and the complex momentum.

Buslaev's ansatz played a very stimulating role for our constructions; moreover the asymptotic formulas of the complex WKB method developed here look similar to the ones of Buslaev. However, Buslaev's constructions can not be used as a base for the constructions of the complex WKB method needed here.

1.2.4. We will now formulate the main theorem of the complex WKB method for the family of equations (1.1). For this rewrite (1.1) in terms of the variables

$$u = x - \phi, \quad \varphi = \varepsilon\phi.$$

The equation takes the form:

$$-\frac{d^2}{du^2}\psi(u) + (V(u) + W(\varepsilon u + \varphi))\psi(u) = E\psi(u), \quad u \in \mathbb{R}, \tag{1.19}$$

and the consistency condition (1.3) comes to

$$\psi_{1,2}(u - 1, \varphi + \varepsilon) = \psi_{1,2}(u, \varphi), \quad \varphi, u \in \mathbb{R}. \tag{1.20}$$

Recall that the function $W$ is analytic. So, we may choose the parameter $\varphi$ to be complex. The goal of our version of the complex WKB method is twofold:

1. describe certain canonical domains on the complex plane $\varphi \in \mathbb{C}$,
2. on any given canonical domain, construct a consistent basis of solutions of (1.19) having the simplest asymptotic behavior when $\varepsilon \to 0$ ($u$ being in a fixed interval $[-U, U] \subset \mathbb{R}$).

When describing the constructions of the complex WKB method we shall refer to the equation

$$-\frac{d^2}{du^2}\psi + V(u)\psi = \mathcal{E}\psi, \quad u \in \mathbb{R}, \tag{1.21}$$

in which

$$\mathcal{E} = E - W(\varphi), \tag{1.22}$$

and $\varphi \in \mathbb{C}$ is fixed. The main theorem of the complex WKB method in the adiabatic case is

**Theorem 1.1.** *Fix $U \in \mathbb{R}$. Let $K$ be a bounded admissible canonical domain for the family of equations (1.1). For sufficiently small $\varepsilon$, there exist a consistent basis $f_{\pm}$ defined for $u \in \mathbb{R}$ and $\varphi \in K$ and having the following properties:*

- *for any fixed $u \in \mathbb{R}$, the functions $f_{\pm}(u, \varphi)$ are analytic in $\varphi \in K$.*
- *for $-U \leq u \leq U$, the functions $f_{\pm}(u, \varphi)$ have the asymptotic representations*

$$f_{\pm}(u, \varphi) = e^{\pm\frac{i}{\varepsilon}\int_0^{\varphi} \kappa \, d\varphi}(\Psi_{\pm}(u, \varphi) + o(1)), \quad \varepsilon \to 0. \tag{1.23}$$

*Here $\Psi_{\pm}$ are the canonical Bloch solutions of equation (1.21) corresponding to the domain $K$.*

- *the error estimates in (1.23) are uniform in $\varphi$ and $u$ and may differentiated once in $u$.*

Roughly, a canonical domain is a domain such that through any point of it, one can draw a smooth curve along which the function $\operatorname{Im} \int_{\varphi_0}^{\varphi} \kappa \, d\varphi$ is monotonically increasing and the function $\operatorname{Im} \int_{\varphi_0}^{\varphi} (\kappa - \pi) \, d\varphi$ is monotonically decreasing. An admissible sub-domain of a canonical domain is this canonical domain without a fixed $\delta$-vicinity of its boundary. The Bloch solutions mentioned in the theorem have the form
$$\Psi_\pm(u, \varphi) = A_\pm(\varphi) \, \psi_\pm(u, \varphi).$$
Here $\psi_\pm$ are two linearly independent Bloch solutions of equation (1.21) defined for $\varphi \in K$. The factors $A_\pm$ are constructed in terms of the periodic components of these Bloch solutions by the same formulae as in the case of (1.17). The reader will find precise definitions in section 4.

We prove Theorem 1.1 in section 3. Let us emphasize here the crucial role played by the consistency condition; indeed, it is this condition that allows us to control the dependence of the solutions on $\varphi$ and leads to the asymptotics (1.23) and to the definition of the canonical domains.

In section 4, we consider the equation
$$(1.24) \qquad -\frac{d^2}{dx^2}\psi(x) + V(x - \phi) + \cos(\varepsilon x)\psi(x) = E\psi(x), \quad x \in \mathbb{R},$$
analyze the analytic properties of the complex momentum $\kappa$ and describe all the related geometric constructions and some of the canonical domains. These canonical domains are shown on Fig. 7 and Fig. 9. They remind us of the canonical domains arising in the semi-classical analysis of Harper equation
$$(1.25) \qquad \frac{\psi(x+h) + \psi(x-h)}{2} + \lambda \cos x \, \psi(x) = E\psi(x),$$
in which $h$ is the semi-classical parameter (the case of $\lambda = 1$ was investigated in [7]). We will discuss this similarity later.

1.3. **Asymptotics of the monodromy matrices.** For any fixed $u \in \mathbb{R}$, the WKB solutions described in the Theorem 1.1 appear to be analytic in a whole strip
$$Y_1 < \operatorname{Im} \varphi < Y_2, \quad Y_1, Y_2 \in \mathbb{R},$$
of the complex plane. Indeed, being analytic in a given canonical domain, they can continued analytically outside just by means of formulae (1.20). As a result, one can see that the corresponding monodromy matrices are analytic in the same strip. In sections 5 – 6, we apply the above constructions to equation (1.24). The structure of the canonical domains for this equation and, thus, the asymptotics of the monodromy matrix strongly depends on the value of the spectral parameter $E$. Denote by $E_1$ and $E_2$ the ends of the first spectral interval of the periodic Schrödinger equation (1.14). For the sake of the definiteness, we assume that
$$(1.26) \qquad E - 1 \leq E_1 - \delta, \quad E_1 + \delta < E + 1 < E_2 - \delta, \quad E_2 < E_3$$
where $\delta$ is a fixed positive number. Starting with a consistent basis corresponding to one of the canonical domains, we prove in section 6 the

**Theorem 1.2.** *Under the condition (1.26), for sufficiently small positive $\varepsilon$, there exists a consistent basis and a positive number $Y$ such that the corresponding monodromy matrix is analytic in the strip $|\operatorname{Im} \varphi| \leq Y$ and can be represented in the form*

$$(1.27) \qquad M = T \cdot \begin{pmatrix} a_0 & b_0 + b_1 e^{iz} \\ c_0 + c_1 e^{-iz} & d_0 + d_1 e^{iz} + d_{-1} e^{-iz} \end{pmatrix} \cdot T^{-1}.$$

*Here*

- *$T$ is a constant $2 \times 2$ diagonal matrix.*
- *$z = 2\pi\phi + C$ ($C$ is a real constant),*
- *the coefficients $a$, $b$, $c$ and $d$ admit the asymptotics:*

$$d_0 = \tfrac{2}{t} \cos(\phi_1 + o(1))\,(1+o(1)), \quad d_1 = -\tfrac{t_1}{t}\,(1+o(1)), \quad d_{-1} = -\tfrac{t_1}{t}\,(1+o(1)),$$

$$b_0 = ie^{i\phi_1}\,(1+o(1)), \quad b_1 = -it_1\,(1+o(1)),$$

$$c_0 = -ie^{i\phi_1}\,(1+o(1)), \quad c_1 = it_1\,(1+o(1)),$$

$$a_0 = te^{i\phi_1}\,(1+o(1)),$$

*where $\phi_1$, $t$ and $t_1$ are nonzero real coefficients depending only on $E$. The asymptotics are uniform in $E$ and in $\varphi$.*

Let us discuss the above result. First, we describe the coefficients $\phi_1$, $t$ and $t_1$. The discussion given below easily follows from formulae (5.1). Remember that the complex momentum $\kappa(\varphi)$ is a multi-valued function. Under the condition (1.26), it has branching points at the real line (see Fig. 2–4). They are described by the relation

$$E_1 + \cos \varphi = E.$$

There is only one branching point on the interval $(0, \pi)$. Denote it by $\varphi_1$. All the other real branching points are the points $\varphi = \pm\varphi_1 + 2\pi l$, $l \in \mathbb{Z}$ (see Fig. 2–4). There is a branch of the complex momentum that belongs to $i\mathbb{R}_+$ on the intervals $(-\varphi_1 + 2\pi l, \varphi_1 + 2\pi l)$ and to $\mathbb{R}_+$ on the intervals $(\varphi_1 + 2\pi l, 2\pi - \varphi_1 + 2\pi l)$, $l \in \mathbb{Z}$. Denote this branch by $\kappa_*$.

Consider the real isoenergy curve (1.13) (i.e. the intersection of the curve (1.13) and the plane $(\varphi, \kappa) \in \mathbb{R}^2$). Under the condition (1.26), it corresponds to the figure given bellow. This curve is $2\pi$-periodic in $\varphi$ and $\kappa$. The branch of this curve connecting the points $\varphi_1$ and $2\pi - \varphi_1$ is described by the formula

$$\kappa = \kappa_*(\varphi), \quad \varphi_1 \leq \varphi \leq 2\pi - \varphi_1.$$

The coefficient $t$ is the tunneling coefficient describing the tunneling between the branches of the isoenergy curve separated by the points $-\varphi_1$ and $\varphi_1$. The coefficient $t_1$ is the tunneling coefficient describing the tunneling of between the branches of the isoenergy curve separated by the points $\kappa_1$ and $2\pi - \kappa_1$. The $\phi_1$ is the phase integral associated to the closed branch connecting the points $\varphi_1$ and

$2\pi - \varphi_1$. More precisely,

(1.28)
$$\phi_1 = \frac{1}{\varepsilon} \int_{\varphi_1}^{2\pi-\varphi_1} \kappa_* d\varphi, \ t = \exp(\frac{i}{\varepsilon} \int_{-\varphi_1}^{\varphi_1} \kappa_* d\varphi), \ t_1 = \exp(\frac{i}{\varepsilon} \int_{\kappa_1}^{2\pi-\kappa_1} (\varphi^* - \pi) \, d\kappa).$$

Here, $\varphi^*(\kappa)$ is a complex branch of the isoenergy curve connecting the points $\varphi = \pi$, $\kappa = \kappa_1$ and $\varphi = \pi$, $\kappa = 2\pi - \kappa_1$. Along this branch $\varphi^*(\kappa) - \pi \in i\mathbb{R}_+$. Note that the tunneling coefficients $t$ and $t_1$ are exponentially small in $1/\varepsilon$, and that the phase integral $\phi_1$ is real.

The integrals defining $\phi_1$, $t_1$ and $t$ can also be written as contour integrals in $(\kappa, \varphi)$-complex plane. For $\phi_1$ (resp. $t$, $t_1$), one integrates the form $\kappa d\varphi$ along the contour $\gamma_1$ (resp. $\gamma_2$, $\gamma_3$). The contour $\gamma_2$ (resp. $\gamma_3$) is imaginary in $\kappa$ (resp. in $\varphi - \pi$) and real in $\varphi$ (resp. in $\kappa$) (see figure below).

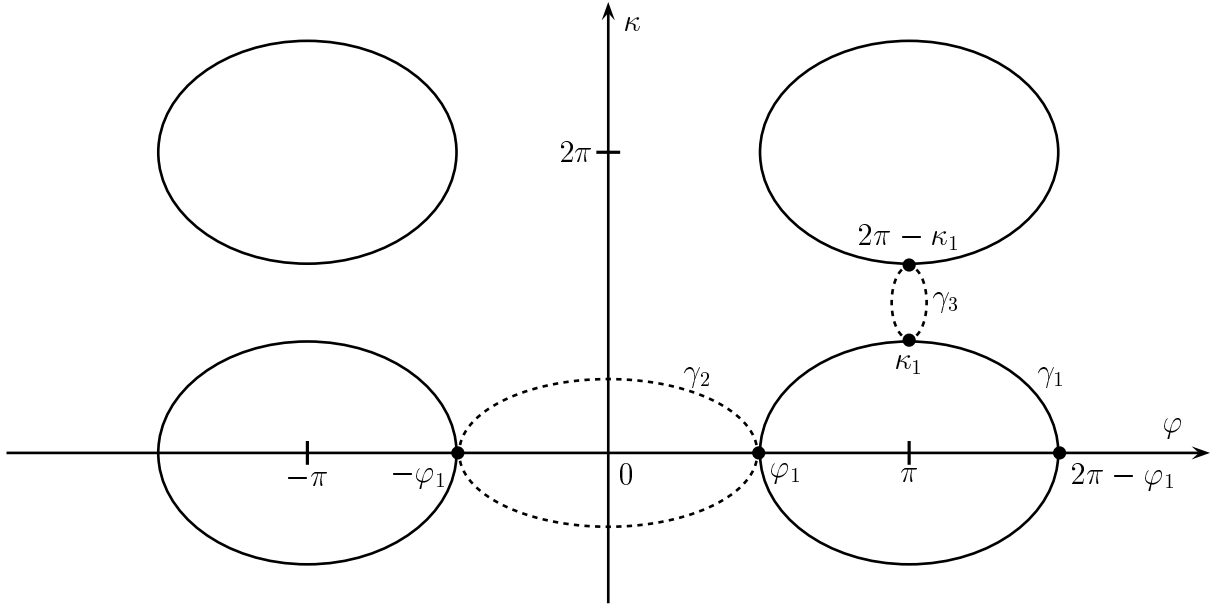

The topological structure of the real isoenergy curve (1.13) is the same as the one of the curve $\cos p + \cos x = E$ with $-2 < E < 2$. The latter one is the real isoenergy curve defined by the symbol of Harper equation (1.25) ($\lambda = 1$). This leads to a similarity of the above mentioned canonical domains for equations (1.25) and (1.1) in a neighborhood of the real line.

We prove Theorem 1.2 in section 5. The calculations leading to Theorem 1.2 use standard ideas of the classical complex WKB method. Due to the similarity of the canonical domains, they are parallel to the calculations of made in [5] for the semi-classical analysis of monodromy matrices for Harper equation.

Let us note that it is Theorem 1.1 which reduces the asymptotic analysis of the almost periodic family of equations (1.1) to the investigation of objects, most of which well known in classical versions of the complex WKB method.

We finish this subsection by discussing the asymptotic analytic structure of the monodromy matrix (1.27). In the Theorem 1.2, we have described the asymptotics of the first Fourier coefficients of the monodromy matrix. We can show that if the periodic Schrödinger operator with the potential $V$ has infinite number of gaps in its spectrum, then the Fourier series are also infinite. However, in the strip $|\operatorname{Im} \varphi| \leq Y$, the leading term of this matrix is a trigonometric polynomial of order one.

Consider a family of monodromy matrices generated by Harper equation in course of the monodromization renormalization procedure, see [9]. The matrices of this family have the form

$$
(1.29) \qquad \begin{pmatrix} a - 2\lambda \cos z & s + te^{-iz} \\ -s - te^{iz} & st/\lambda \end{pmatrix},
$$

where $a, s, t$ are constant parameters. Since the determinant of this matrix equals to one, only two of them are independent. Comparing (1.27) with (1.29), we see that one can hope to carry out the spectral analysis of the family of equations (1.1) by means of the monodromization method developed for Harper equation.

### 1.4. First spectral results and main conjectures.

1.4.1. Let $J_\delta$ be the interval of energies $E$ defined by (1.26). In Section 6 we show that Theorem 1.2 implies

**Theorem 1.3.** *Consider the points $E^{(l)}$, $l \in \mathbb{Z}$, defined by*

$$(1.30) \qquad \phi_1(E^{(l)}) = \pi/2 + \pi l,$$

*and situated in $J_\delta$. Let $L$ be the set of integers $l$ such that $E^{(l)}$ exists. Then, for sufficiently small $\varepsilon$,*

- *the distances between these points satisfy the inequalities*

$$c_1 \varepsilon \leq |E^{(l)} - E^{(l-1)}| \leq c_2 \varepsilon,$$

  *where $c_1$ and $c_2$ are two positive constant independent of $\varepsilon$.*
- *there exists a collection of intervals $(I_l)_{l \in L}$ such that $I_l \subset J_\delta$ and such that the spectrum of equation (1.24) in $J_\delta$ is contained in $\cup_{l \in L} I_l$.*

*Moreover for $l \in L$, we have*

- *the interval $I_l$ is in an $o(\varepsilon)$-vicinity of $E^{(l)}$,*
- *the length of $I_l$ is equals to $(t(E^{(l)}) + t_1(E^{(l)})$ (1+o(1))*.

This theorem is similar to a result on the spectrum of Harper equation first proved in [14]. We prove it by using the idea on the relation between the Bloch solutions of (1.1) and (1.9) discussed above.

1.4.2. Consider the monodromy matrix described in Theorem 1.2 for $E$ being in one of the subintervals $I_l$ described in Theorem 1.3. Clearly, the leading terms of its asymptotics have the form

$$(1.31) \qquad M \sim \begin{pmatrix} 0 & \pm 1 \\ \mp 1 & F_l - 2\lambda_l \cos z \end{pmatrix},$$

where

$$\lambda_l = t_1(E^{(l)})/t(E^{(l)}), \quad F_l = \cos \phi_1(E)/t(E^{(l)}).$$

If we forget about the error terms and let $x = 2\pi\varphi/\varepsilon + C$, then equation (1.9) becomes equivalent to the Harper equation

(1.32) $$f(x+h) + f(x-h) - 2\lambda_l \cos x \, f(x) = F_l(E) f(x).$$

This leads to a list of conjectures on the spectral properties of family of equations (1.1). Before discussing these conjectures, we have to make a remark about the individual equations of this family.

1.4.3. For a general metrically transitive family of equations, the spectral properties have to be formulated in probabilistic terms (see, for example, [20]). In particular, certain basic spectral results appear to be valid with probability one (i.e. for almost all the individual equations). Here, it is the parameter $\phi$ that plays the role of the "random variable". Trying to study the spectral properties of one of the equations (1.1), i.e. the equation corresponding to a given value of $\phi \in \mathbb{R}$, one has to restrict the solutions of (1.9) on the lattice

$$x_n = nh - \phi, \quad n \in \mathbb{Z}.$$

The restrictions of the solutions

$$\chi_n = \chi(nh - \phi), \quad n \in \mathbb{Z},$$

satisfy the finite difference equation

(1.33) $$\chi_{n+1} = M(nh - \phi) \chi_n, \quad n \in \mathbb{Z}.$$

So, trying to analyze typical properties of individual equations of the family (1.1), one has to study solutions of equation (1.33).

1.4.4. Consider the finite difference equation

(1.34) $$\frac{\psi_{n+1} + \psi_{n-1}}{2} + \lambda \cos(nh + \varphi) \psi_n = E\psi_n, \quad \psi \in l_2.$$

It is the well known and intensively studied Almost Mathieu equation (for an almost complete list of references see [16], [15]).
Clearly, (1.34) coincides with the restriction of (1.32) to a lattice. It is believed (and, actually, proved in many situations) that

- if $\lambda > 1$ and $h$ is Diophantine, then the spectrum of (1.34) is pure point for almost all $\varphi \in \mathbb{R}$;
- if $0 < \lambda < 1$ and $h \notin \mathbb{Q}$, then the spectrum is absolutely continuous for almost all $\varphi \in \mathbb{R}$.
- if $h \notin \mathbb{Q}$ and $\lambda > 0$, then the spectrum is a Cantor set.

This immediately leads to three conjectures on the spectral properties of (1.1). However, in [13], it is proved that, for $\lambda = 1$ and a full measure set of $h$, the spectrum of (1.34) is singularly continuous. But since the spectrum of (1.34) crucially depends on the value of $\lambda$, the analysis of the spectrum of (1.1) on the intervals where $\lambda_m = t_1(E_l)/t(E_l)$ is exponentially close to 1, or where the difference $1 - t_1(E)/t(E)$ even changes its sign, has to be very delicate and, so, the results known for the Almost Mathieu equation lead here not to a conjecture, but to an open question.

## Contents

## 2. Periodic Schrödinger operators

Here, we collect some well known facts (see, for example, [17], [18]) on one-dimensional periodic Schrödinger operators (1.14). We assume that $V$ is sufficiently regular real valued function,

$$V(x+1) = V(x), \quad x \in \mathbb{R}.$$

The results described below are valid if $V \in L^2_{\text{loc}}(\mathbb{R})$.

2.1. **Bloch solutions.** Let $\psi$ be a solution of (1.14) satisfying the relation

(2.1) $$\psi(x+1) = \lambda\psi(x), \quad \forall x \in \mathbb{R},$$

with $\lambda$ independent of $x$. Such solution is called a *Bloch* solution, and the number $\lambda$ is called the *Floquet multiplier*.

Consider the spectrum of equation (1.14) in $L^2(\mathbb{R})$. It consists of intervals $[E_1, E_2]$, $[E_3, E_4]$ ... $[E_{2n+1}, E_{2n+2}]$, ..., of the real axis such that

(2.2) $$E_1 < E_2 \leq E_3 < E_4 \ldots E_{2n} \leq E_{2n+1} < E_{2n+2} \leq \ldots$$

(2.3) $$E_n \to +\infty, \quad n \to +\infty.$$

The intervals $[E_{2n-1}, E_{2n}]$ $(n \geq 1)$ are called the *spectral zones*, and the intervals $(E_{2n}, E_{2n+1})$, $n \geq 1$, are called the *spectral gaps*. Usually, the number of the gaps is infinite.

Consider two copies of the complex plane $E \in \mathbb{C}$ cut along the spectral zones. Paste them together to get a Riemann surface with square root branching points. We call his Riemann surface $\Gamma$. Denote one of the copies of this surface by $\Gamma_+$ and the other one by $\Gamma_-$. $\Gamma_+$ (resp. $\Gamma_-$) will be refered to as the upper (resp. lower) sheet.

One can construct a Bloch solution $\psi(x, E)$ of the equation (1.14) meromorphic on this Riemann surface. The poles of this solution are located in the spectral gaps. More precisely, each spectral gap contains precisely one of these poles. All the poles are simple. They are located either on $\Gamma_+$ or on $\Gamma_-$. The position of a pole is independent of $x$. The Bloch solution $\psi(x, E)$ has the asymptotics

(2.4) $$\psi(x, E) = e^{i\sqrt{E}\,x}(1 + o(1)), \quad E \to \infty,$$

(2.5) $$E \in \Gamma, \quad 0 < \arg E < 2\pi.$$

We define the branch of the square root so that, for $E < 0$, $i\sqrt{E} > 0$ if $E \in \Gamma_+$.

Denote by $\psi_\pm$ the restrictions of the above defined Bloch solution on $\Gamma_\pm$. Outside the ends of the spectral zones, these functions are linearly independent solutions of (1.14). Along the cuts, we have

(2.6) $$\psi_\pm(x, E - i0) = \psi_\mp(x, E + i0), \quad E_{2l+1} < E < E_{2l+2}, \quad l \in \mathbb{N}.$$

2.2. **Bloch quasi-momentum.**

2.2.1. Consider the Bloch solution $\psi(x, E)$ introduced in the previous subsection. The corresponding Floquet multiplier $\lambda(E)$ is analytic on $\Gamma$. Represent it in the form

(2.7) $$\lambda(E) = \exp(ik(E)).$$

The function $k(E)$ is called *Bloch quasi-momentum*. It has the same branching points as $\psi(x, E)$, but the corresponding Riemann surface is more complicated. To describe the main properties of $k$, consider the complex plane cut along the real line from $E_1$ to $+\infty$. Denote the cut plane by $\mathbb{C}_0$. One can fix there a single valued branch of the quasi-momentum by the condition

(2.8) $$ik_0(E) < 0, \quad E < E_1.$$

The image of $\mathbb{C}_0$ by $k_0$ is located in the upper half of the complex plane,

(2.9) $$\operatorname{Im} k_0(E) > 0, \quad E \in \mathbb{C}_0.$$

In Fig. 1, we drew several curves in $\mathbb{C}_0$ and their images under the transformation $E \mapsto k_0(E)$.

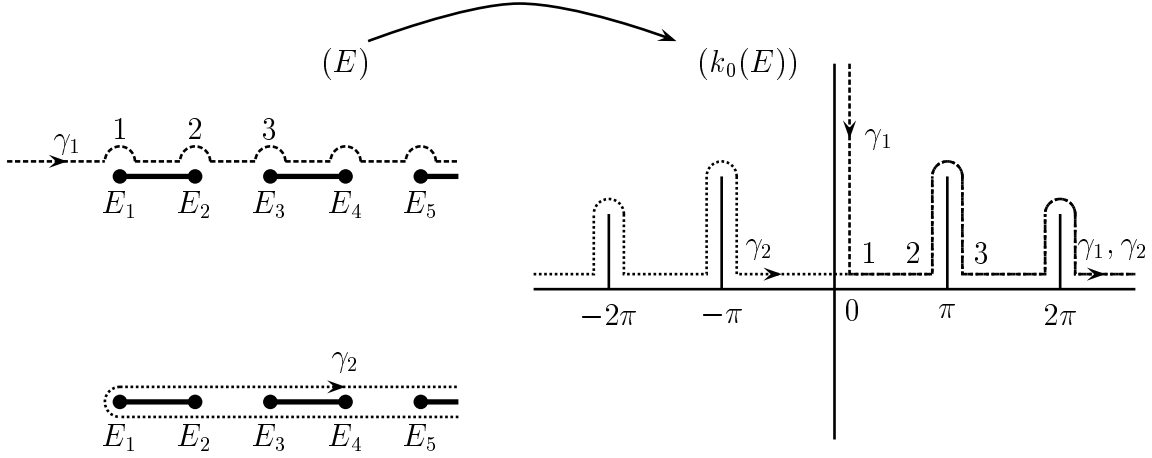

Figure 1.

The quasi-momentum $k_0(E)$ is real along the spectral zones, and, along the spectral gaps, its real part is constant; in particular, we have

(2.10) $$k_0(E_1) = 0, \quad k_0(E_{2l} \pm i0) = k_0(E_{2l+1} \pm i0) = \pm \pi l, \quad l \in \mathbb{N}.$$

Let $E_m$ be one of the branching points of $k$. Then both of the functions

(2.11) $$f_m^\pm(E) = (k_0(E \pm i0) - k_0(E_m \pm i0))/\sqrt{E - E_m}, \quad E \in \mathbb{R},$$

can be analytically continued in a small vicinity of the branching point $E_m$.

Finally note that

(2.12) $$k_0(E) = \sqrt{E} + O(1/\sqrt{E}), \quad |E| \to \infty,$$

where $E \in \Gamma$ and $0 < \arg E < 2\pi$. This representation is uniform in $\arg E$ if $0 < \delta < \arg E < 2\pi - \delta$ for any fixed $\delta$.

2.2.2. The values of the quasi-momentum $k_0$ on the different sides of the cut $E_1 < E < +\infty$ are related to each other by the formula

(2.13) $$k_0(E + i0) = -\overline{k_0(E - i0)}, \quad E_1 \leq E.$$

Consider the spectral gap $E_{2l} < E < E_{2l+1}$, $l \in \mathbb{N}$. Let $\mathbb{C}_l$ be the complex plane cut from $-\infty$ to $E_{2l}$ and from $E_{2l+1}$ to $+\infty$. Denote by $k_l$ the branch of the quasimomentum analytic on $\mathbb{C}_l$ and coinciding with $k_0$ for $\operatorname{Im} E > 0$. One has

(2.14) $$k_l(E + i0) + \overline{k_l(E - i0)} = 2\pi l, \quad \text{if } E < E_{2l} \text{ or } E > E_{2l+1}.$$

This relation follows from (2.10) and from the analyticity of the functions given by (2.11).

2.3. **Periodic components of the Bloch solution.** Fix $l \in \mathbb{N} \cup \{0\}$. Consider the solutions $\psi_\pm$ introduced in subsection 3.1. They can be represented in the form

$$\psi_\pm(x, E) = e^{\pm i k_l(E) x} p_l^\pm(x, E), \quad E \in \mathbb{C}_l, \tag{2.15}$$

where $p_l^\pm(x, E)$ are functions periodic in $x$,

$$p_l^\pm(x + 1, E) = p_l^\pm(x, E), \quad \forall x \in \mathbb{R}. \tag{2.16}$$

2.4. **Useful formulae.** We shall use the following two integral relations.

**Lemma 2.1.** *Let $E \neq E_l$, $l \in \mathbb{N}$, and let $x \in \mathbb{R}$. Then*

$$\int_{x-1}^{x} \psi_+(x, E) \, \psi_-(x, E) \, dx = -i k'(E) \, w(E), \tag{2.17}$$

*where $w = w(\psi_+, \psi_-)$ and*

$$\begin{aligned}
2 \int_{x-1}^{x} \frac{\partial \psi_+}{\partial E}(x, E) \, \psi_-(x, E) \, dx = \\
= -(i k''(E) + (k'(E))^2) \, w(E) - 2 i k'(E) \, w\left(\frac{\partial \psi_+}{\partial E}(x, E), \psi_-(x, E)\right).
\end{aligned} \tag{2.18}$$

*Here $w(f(x), g(x)) = f'(x) g(x) - g'(x) f(x)$.*

*Proof.* Though these two formulae are known, we describe their simple proofs here. Begin with the formula (2.17). Derivating (1.14) in $E$, we get:

$$-\frac{\partial^3 \psi_+}{\partial x^2 \, \partial E}(x, E) + V(x) \frac{\partial \psi_+}{\partial E}(x, E) = E \frac{\partial \psi_+}{\partial E}(x, E) + \psi_+(x, E). \tag{2.19}$$

Multiplying this equality by $\psi_-(x, E)$, integrating from $x - 1$ to $x$, then integrating twice by parts, we obtain

$$-w\left(\frac{\partial \psi_+}{\partial E}(x', E), \psi_-(x', E)\right)\bigg|_{x'=x-1}^{x'=x} = \int_{x-1}^{x} \psi_+(x', E) \, \psi_-(x', E) \, dx'. \tag{2.20}$$

Now, by means of the formulae

$$\psi_\pm(x - 1, E) = e^{\mp i k(E)} \psi_\pm(x, E),$$

we can express all the functions in the left hand side of (2.20) in terms of functions of the argument $x' = x$. This leads to

$$w\left(\frac{\partial \psi_+}{\partial E}(x', E), \psi_-(x', E)\right)\bigg|_{x'=x-1}^{x'=x} = i k'(E) \, w(E),$$

which implies (2.17). The formula (2.18) can be obtained similarly. But now, instead of (2.19), we start with the relation

$$-\frac{\partial^4 \psi_+}{\partial x^2 \, \partial E^2}(x, E) + V(x) \frac{\partial^2 \psi_+}{\partial E^2}(x, E) = E \frac{\partial^2 \psi_+}{\partial E^2}(x, E) + 2 \frac{\partial \psi_+}{\partial E}(x, E).$$

We omit the elementary details. □

In section 4, we shall also need the

**Corollary 2.1.** *Under the assumptions of Lemma 2.1,*

(2.21)
$$-w\left(\psi_+(x,E),\ \frac{\partial \psi_-}{\partial E}(x,E)\right) + \int_{x-1}^{x} \psi_+(x',E)\,\psi_-(x',E)\,x'\,dx' =$$
$$= -w'(E) + w(E)\left(-\frac{1}{2}\frac{d}{dE}\ln k'(E) + \frac{i}{2}k'(E) - \omega^+(E)\right),$$

*where*

$$\omega^+(E) = -\int_0^1 \frac{\partial p_+}{\partial E}(x,E)\,p_-(x,E)\,dx\ \left(\int_0^1 p_+(x,E)\,p_-(x,E)\,dx\right)^{-1}.$$

This corollary easily follows from formulae (2.17) and (2.18). Note that, though the Bloch quasi-momentum $k$ is defined up to an integral multiple of $2\pi$ (and, therefore, the periodic components $p_\pm$ of $\psi_\pm$ are defined up to some factor $e^{\mp 2\pi i l x}$), all the objects in (2.21) are uniquely defined.

## 3. Main theorem of the complex WKB method

In this section, we consider the family of equations (1.19), define the corresponding canonical domains and, then, prove Theorem 1.1.

### 3.1. Canonical domains.
The canonical domain is the main geometric notion of the complex WKB method. Its origin is tightly connected with the construction of the analytic solutions having standard asymptotic behavior in the parameter $\varphi \in \mathbb{C}$.

3.1.1. The canonical domains can be described in terms of the complex momentum $\kappa(\varphi)$. Remind that it is defined by formula (1.15) in terms of the Bloch quasi-momentum $k(E)$ of the operator $H_0 = -\dfrac{d^2}{du^2} + V(u)$ with the same $V$ as in (1.19). The complex momentum $\kappa$ is a multi-valued analytic function. Its branching points are related to the branching points of the quasi-momentum by the relations

(3.1) $$E_l = E - W(\varphi),\quad l \in \mathbb{N},$$

where $E_l$ are the ends of the spectral zones of the operator $H_0$. The conformal properties of $\kappa$ and the structure of its Riemann surface depend on the analytic properties of the function $W$ and on the value of the spectral parameter $E$. For $W(\varphi) = \cos\varphi$, the analytic properties of the complex momentum and of its Riemann surface are described in subsection 4.1.

3.1.2. In what follows, $\zeta_1$ and $\zeta_2$ are two points in $\overline{\mathbb{C}}$ such that

$$\operatorname{Im}\zeta_1 < \operatorname{Im}\zeta_2.$$

These points can either be situated in the complex plane or be the points $\pm i\infty$. We shall denote by $\gamma$ a smooth curve going from $\zeta_1$ to $\zeta_2$. When we will integrate along such a curve, the curve will always be oriented from $\zeta_1$ to $\zeta_2$. We now state the definitions.

3.1.3. We call a smooth curve $\gamma \subset \mathbb{C}$ *vertical* if it intersects the lines $\operatorname{Im} z = \operatorname{Const}$ at non-zero angles $\theta$.

If there is a positive number $\delta_1$ such that, at any point of $\gamma$, the intersection angle $\theta$ satisfies the inequality

(3.2) $$\delta_1 < \theta < \pi - \delta_1$$

then, we call $\gamma$ *strictly vertical*.

3.1.4. Let $\gamma$ be a vertical curve. On $\gamma$, fix a continuous branch of the momentum $\kappa$. We call $\gamma$ *canonical* if, along $\gamma$, $\operatorname{Im} \int^{\varphi} \kappa \, d\varphi$ is monotonically increasing, and $Im \int^{\varphi} (\kappa - \pi) \, d\varphi$ is monotonically decreasing.

Assume that $\gamma$ is strictly vertical. If there exist two positive numbers $\delta_2$ and $\delta_3$ such that

(3.3) $$\operatorname{Im} \int_{\varphi}^{\varphi'} \kappa \, d\varphi \geq \delta_2 \operatorname{Im} (\varphi' - \varphi)$$

(3.4) $$\operatorname{Im} \int_{\varphi}^{\varphi'} (\kappa - \pi) \, d\varphi \leq -\delta_3 \operatorname{Im} (\varphi' - \varphi),$$

then, we call $\gamma$ *strictly canonical*.

3.1.5. Let $K$ a simply connected domain containing no branching points of the complex momentum $\kappa$. On $K$, fix a continuous branch $\kappa$. We call $K$ a *canonical domain* if it is the union of canonical curves connecting some two points $\zeta_1$ and $\zeta_2$ fixed as above. To be precise, we have to exclude these two points from $K$.

If there are three positive numbers $\delta_1$, $\delta_2$ and $\delta_3$ such that $K$ is a union of canonical curves satisfying (3.2) – (3.4), then, we call $K$ *strictly canonical*.

3.1.6. Assume that $K \subset \mathbb{C}$ is a canonical domain. Denote by $\partial K$ its boundary. Fix a positive number $\delta_4$. We call the domain

(3.5) $$\mathcal{C} = \{z \in K \,|\, \operatorname{dist}(z, \partial K) > \delta_4\}$$

an *admissible sub-domain* of $K$.

Note that the branching points of the complex momentum are outside of $\mathcal{C}$, at a distance greater than $\delta_4$.

**Remark 3.1.** We will describe the other classical geometric objects and constructions usually related to the complex WKB method (for example, Stokes lines) only in the case when $W(\varphi) = \cos(\varphi)$. These constructions will not be needed neither to formulate nor to prove our main theorem. However they are essential to understand the geometry of canonical domains for a given equation.

3.2. **Canonical Bloch solutions.** To describe the asymptotic formulae of of the complex WKB method for equation (1.1), one needs specially normalized Bloch solutions of the equation (1.21).

3.2.1. Assume that $D$ is a simply connected domain in the plane $\varphi \in \mathbb{C}$ containing no branching points of $\kappa$.

Consider (1.21) as an equation with spectral parameter $\mathcal{E}$ and consider the associated Riemann surface $\Gamma$. Let $\psi(u, \mathcal{E})$ be the Bloch solution of equation (1.21) meromorphic on $\Gamma$.

The mapping $\varphi \to E - W(\varphi)$ maps $D$ onto a domain $\hat{D} \subset \mathbb{C}$. Clearly, $\hat{D}$ is simply connected and does not contain any branching point of $\psi(u, \mathcal{E})$. So, since $\Gamma$ is a two sheeted Riemann surface, we can also consider two different copies $\hat{D}_+$ and $\hat{D}_-$ of the domain $\hat{D}$.

Denote the restrictions of $\psi(u, \mathcal{E})$ on $\hat{D}_\pm$ by $\psi_\pm$. These functions are two Bloch solutions of (1.21). Fix on $\hat{D}_+$ a continuous branch of the Bloch quasi-momentum $k(\mathcal{E})$. Clearly,
$$\psi_\pm(u, \mathcal{E}) = e^{\pm ik(\mathcal{E})u} p_\pm(u, \mathcal{E}),$$
where $p_\pm$ are 1-periodic in $u$,
$$p_\pm(u+1, \mathcal{E}) = p_\pm(u, \mathcal{E}).$$
So, $\pm k(\mathcal{E})$ are the quasi-momenta of $\psi_\pm(u, \mathcal{E})$. Note also, that $p_\pm(u, \mathcal{E})$ are meromorphic in $\mathcal{E} \in \hat{D}$.

3.2.2. Let
$$(3.6) \qquad g_\pm(\mathcal{E}) = -\frac{\int_0^1 p_\mp(x, \mathcal{E}) \frac{\partial p_\pm}{\partial \mathcal{E}}(x, \mathcal{E}) dx}{\int_0^1 p_+(x, \mathcal{E}) p_-(x, \mathcal{E}) dx}, \quad \mathcal{E} \in \hat{D}.$$

In view of (2.17), $g_\pm$ are meromorphic in $\hat{D}$.

Define
$$(3.7) \qquad f_\pm(u, \mathcal{E}) = \sqrt{k'_E(\mathcal{E})} \, e^{\int^\mathcal{E} g_\pm(e)de} \, \psi_\pm^0(u, \mathcal{E}), \quad \mathcal{E} \in \hat{D}.$$

We call the functions
$$(3.8) \qquad \Psi_\pm(u, \varphi) = f_\pm(u, \mathcal{E}), \quad \mathcal{E} = E - W(\varphi),$$
*canonical Bloch solutions associated to the domain $D$*. Note that they are independent of the choice of the branch $k$. One has

**Lemma 3.1.** *The Bloch solutions $\Psi_\pm(u, \varphi)$ are analytic in $\varphi \in D$.*

*Proof.* It suffices to prove that the functions $f_\pm(u, \mathcal{E})$ are analytic in $\hat{D}$. Let us prove it for $f_+$. The case of $f_-$ can be dealt with in the same way. In what follows, we denote by $\sim$ the equalities up to some terms analytic in $\mathcal{E} \in \hat{D}$.

Obviously, the formulae from Lemma 2.1 remain valid for $f_\pm(u, \mathcal{E})$ when $\mathcal{E} \in \hat{D}$. They imply that
$$(3.9) \qquad g_\pm(\mathcal{E}) \sim -\frac{1}{2}\frac{d}{d\mathcal{E}} \ln k'(\mathcal{E}) - \frac{1}{w(\psi_+, \psi_-)} w\left(\frac{\partial \psi_+}{\partial \mathcal{E}}(u, \mathcal{E}), \psi_-(u, \mathcal{E})\right)\bigg|_{u=1}.$$

We have omitted the term $\frac{i}{2} k'(\mathcal{E})$ since it is analytic outside the branching points.

Consider the second term in the (3.9). Remember that the poles of $\psi_\pm(u, \mathcal{E})$ are independent of $u$. Let $\mathcal{E}_0 \in \hat{D}$ be a pole of $\psi_\pm(u, \mathcal{E})$. In a neighborhood of $\mathcal{E}_0$, we have

$$\psi_\pm(u, \mathcal{E}) = \frac{1}{(\mathcal{E} - \mathcal{E}_0)^{\nu_\pm}} \tilde{\psi}_\pm(u, \mathcal{E}),$$

where $\nu_\pm$ is the multiplicity of $\mathcal{E}_0$ as a pole of $\psi_\pm(u, \mathcal{E})$, and $\tilde{\psi}_\pm(u, \mathcal{E})$ is a solution of (1.21) analytic in a neighborhood of $\mathcal{E}_0$.

$\tilde{\psi}_\pm(u, \mathcal{E})$ is a Bloch solution with the same Bloch quasi-momentum $\pm k(\mathcal{E})$ as $\psi_\pm(u, \mathcal{E})$. Since $k(\mathcal{E}) \in \pi \mathbb{Z}$ only at the branching points, the solutions $\tilde{\psi}_\pm(u, \mathcal{E})$ are linearly independent. These arguments imply that, in a vicinity of $\mathcal{E}_0$,

$$\frac{1}{w(\psi_+, \psi_-)} w \left( \frac{\partial \psi_+}{\partial \mathcal{E}}(u, \mathcal{E}), \psi_-(u, \mathcal{E}) \right) \bigg|_{u=1} \sim -\frac{\nu_\pm}{\mathcal{E} - \mathcal{E}_0},$$

and, thus, (3.9) takes the form

$$g_\pm(\mathcal{E}) \sim -\frac{1}{2} \frac{d}{d\mathcal{E}} \ln k'(\mathcal{E}) + \frac{\nu_\pm}{\mathcal{E} - \mathcal{E}_0}.$$

This implies the analyticity of $f_\pm$ in a vicinity of any point $\mathcal{E}_0 \in \hat{D}$. □

In the sequel, we shall use the notation

(3.10) $$\omega_\pm(\varphi) = -W'(\varphi) \cdot g_\pm(E - W(\varphi)).$$

3.3. **Formulation of the main theorem of the complex WKB method.** We prove

**Theorem 3.1.** *Fix a positive $U$. Let $K$ be a bounded strictly canonical domain for (1.19).*

*For sufficiently small positive $\varepsilon$, there exists a consistent basis $f_\pm(u, \varphi)$ for the family of equations (1.19) with the following properties:*

- *for any fixed $u \in \mathbb{R}$, the functions $f_\pm(u, \varphi)$ are analytic in $\varphi \in \mathcal{C}$;*
- *the functions $f_\pm(u, \varphi)$ have the asymptotic representations*

(3.11) $$f_\pm(u, \varphi) = e^{\pm \frac{i}{\varepsilon} \int_0^\varphi \kappa \, d\varphi} (\Psi_\pm(u, \varphi) + o(1)), \quad \varepsilon \to 0,$$

$$-U \leq u \leq U, \quad \varphi \in \mathcal{C},$$

  *where $\mathcal{C}$ is a fixed admissible sub-domain for $K$ and $\Psi_\pm$ are the canonical Bloch solutions associated to the domain $K$.*
- *the error estimates in (3.11) are uniform in $\varphi$ and $u$.*

*Moreover the asymptotic representations (3.11) are once differentiable in $u$.*

The rest of this section is devoted to the proof of this theorem. For the sake of brevity, we do not check explicitly the uniformity of errors nor the differentiability of the asymptotics: this does not require any additional estimates or discussions, and the reader will easily fill this gap.

3.4. **Sketch of the proof.** Before beginning the actual proof of Theorem 3.1, let us outline the main ideas of this proof.

3.4.1. Let us assume that $\varphi$ belongs to a bounded domain $D \subset \mathbb{C}$ and that $u$ is in a fixed interval $[-U, U] \subset \mathbb{R}$. Clearly,
$$W(\varepsilon u + \varphi) = W(\varphi) + O(\varepsilon), \quad \varepsilon \to 0.$$
So, we may "replace" (1.19) by equation (1.21) with $\mathcal{E}$ given by (1.22). Assuming that the domain $D$ does not contain any branching points of the complex momentum $\kappa(\varphi) = k(E - W(\varphi))$, we introduce two linearly independent Bloch solutions $\psi_\pm^0(u, \varphi)$ of (1.21)–(1.22). Assuming that these two solutions have no singularities in $D$, we prove that the initial equation (1.19) has two linearly independent solutions $\psi_\pm^1$ close to $\psi_\pm^0$ when $\varepsilon \to 0$.

3.4.2. Though the constructed solutions $\psi_\pm^1(u, \varphi)$ are linearly independent and analytic in $\varphi \in D$, they need not be consistent. If there exists any solution $f(u, \varphi)$ of (1.19) satisfying the consistency condition (1.20), it certainly can be represented in the form
$$f(u, \varphi) = A(\varphi)\,\psi_+^1(u, \varphi) + B(\varphi)\,\psi_-^1(u, \varphi)$$
where the coefficients $A$ and $B$ are independent of $u$. The consistency condition is equivalent to a difference equation of the form

(3.12) $$\begin{pmatrix} A(\varphi + \varepsilon) \\ B(\varphi + \varepsilon) \end{pmatrix} = T(\varphi) \begin{pmatrix} A(\varphi) \\ B(\varphi) \end{pmatrix}, \quad \varphi, \varphi + \varepsilon \in D.$$

The coefficient of the matrix $T$ can be calculated in terms of the Wronskians of the solutions $\psi_\pm^1(u - 1, \varphi + \varepsilon)$ and $\psi_\pm^1(u, \varphi)$.

3.4.3. Then, we compute the asymptotics of the matrix $T$ to find out that
$$T = \begin{pmatrix} e^{i\kappa(\varphi) + O(\varepsilon)} & O(\varepsilon) \\ O(\varepsilon) & e^{-i\kappa(\varphi) + O(\varepsilon)} \end{pmatrix}, \quad \varepsilon \to 0,$$
Using a transformation of the form
$$A_1(\varphi) = e^{\frac{i}{\varepsilon}\int_{\varphi_0}^\varphi \kappa(\varphi)\,d\varphi + O(\varepsilon)}\,A(\varphi), \quad B_1(\varphi) = e^{-\frac{i}{\varepsilon}\int_{\varphi_0}^\varphi \kappa(\varphi)\,d\varphi + O(\varepsilon)}\,B(\varphi),$$
we come to the equation

(3.13) $$\begin{pmatrix} A_1(\varphi + \varepsilon) \\ B_1(\varphi + \varepsilon) \end{pmatrix} - \begin{pmatrix} A(\varphi) \\ B(\varphi) \end{pmatrix} = T_1(\varphi) \begin{pmatrix} A(\varphi) \\ B(\varphi) \end{pmatrix}, \quad \varphi, \varphi + \varepsilon \in D,$$

where the norm of the matrix $T_1$ is small.

3.4.4. The operator in the left-hand side of (3.13) can be easily inverted (see, for example [7]); so that, finally, we come to the equation

(3.14) $$\begin{pmatrix} A_1 \\ B_2 \end{pmatrix} = \begin{pmatrix} 1 \\ 0 \end{pmatrix} + K \begin{pmatrix} A_1 \\ B_2 \end{pmatrix},$$

where $K$ is a singular integral operator acting in a space of Hölder functions defined on a vertical curve $\gamma \subset D$. If the curve $\gamma$ is strictly canonical then the norm of the operator is small. So, we can construct a solution of (3.14). This solution can be analytically continued in the domain $D$, and, as a result, one gets an analytic solution satisfying the condition (1.20).

3.5. **Solutions of equation (1.19) close to Bloch solutions of equation (1.21).** Let $D$ be a bounded domain in the complex plane $\varphi$ not containing any branching points of the complex momentum. Define, as in section 3.2.1, the functions $\psi_\pm(u, \mathcal{E})$, and let

(3.15) $$\psi_\pm^0(x, \varphi) = \psi_\pm(x, E - W(\varphi)), \quad \varphi \in D.$$

In the sequel, we shall assume that $\overline{D}$ does not contain any poles of $\psi_\pm^0(x, \varphi)$[1].

The functions $\psi_\pm^0(x, \varphi)$ are two linearly independent Bloch solutions of the equation (1.21) (with $\mathcal{E}$ given by (1.22)) with the Bloch quasi-momenta $\pm\kappa(\varphi)$ related to the quasi-momenta of $\psi_\pm(u, \mathcal{E})$ by the formula

$$\kappa(\varphi) = k(E - W(\varphi)).$$

We denote the Wronskian of the solutions $\psi_\pm^0(x, \varphi)$ by $w_0(\varphi)$ i.e.

$$w_0(\varphi) = w(\psi_+^0(., \varphi), \psi_-^0(., \varphi)).$$

Clearly, $w_0(\varphi)$ is analytic and has no zeros in $D$.

Fix $U > 0$. Let us construct two solutions $\psi_\pm^1$ of equation (1.19) close to $\psi_\pm^0$ for $\varepsilon \to 0$ uniformly in $\varphi \in D$ and $u \in [-U, U]$.

3.5.1. Fix a positive $U_1$ so that $0 < U_1 < U$. Define an integral operator acting on $C([-U, U])$ by the formula

(3.16)
$$R f(u, \varphi) = \psi_+^0(u, \varphi) \int_{-U_1}^u \psi_-^0(x, \varphi) f(x) \, dx - \psi_-^0(u, \varphi) \int_{U_1}^u \psi_+^0(x, \varphi) f(x) \, dx.$$

This operator coincides with a Green's function of equation (1.21) – (1.22). Indeed, letting

$$H_0 = -\frac{d^2}{du^2} + V(u) + W(\varphi) - E,$$

we get

$$H_0 R f = -w_0 f.$$

Clearly, under our assumptions on $V$ and $W$ and the domain $D$, the operator $R$ is bounded,

(3.17) $$\sup_{\varphi \in D} \|R\| \leq C.$$

3.5.2. For $\psi_\pm^1$, take the ansatz

(3.18) $$\psi_\pm^1 = \psi_\pm^0 + \varepsilon R g_\pm$$

Substituting this ansatz into (1.19) and assuming that $\psi_\pm^1$ are solutions of this equation, we get

(3.19) $$-\Delta \cdot \psi_\pm^0 = -w_0 g_\pm + \varepsilon\Delta \cdot R g_\pm,$$

---

[1] From the very beginning, we could have considered the canonical Bloch solutions associated to the domain $D$ (since they are analytic in $D$). However, we did not want to introduce these objects as an "ansatz" but show how they arise naturally in the constructions of the complex WKB method. This will neither complicate nor enlarge our analysis.

where

(3.20) $$\Delta \equiv \Delta\left(\varepsilon u,\,\varphi\right) = \frac{W\left(\varepsilon u + \varphi\right) - W\left(\varphi\right)}{\varepsilon}.$$

3.5.3. Since the Wronskian $w_0\left(\varphi\right)$ has no zeros in $\overline{D}$, then

$$\left|\frac{\Delta}{w_0}\right| \leq C, \quad -U \leq u \leq U, \quad \varphi \in D,$$

where $C$ is clearly independent of $\varepsilon$, $u$ and $\varphi$. Therefore, for sufficiently small $\varepsilon$, equation (3.19) has a unique solution in $C([-U, U])$; this solution satisfies

(3.21) $$\sup_{\varphi \in D} \|g_\pm\|_\infty \leq C.$$

The formula (3.19) implies that

(3.22) $$g_\pm = \frac{\Delta}{w_0} \psi_\pm^0 + O\left(\varepsilon\right), \quad \varphi \in D, \quad -U \leq u \leq U$$

in the sense of the supremum norm.

3.5.4. By formula (3.18), we reconstruct the function $\psi_\pm^1$ in terms of $g_\pm$. This leads to the following results:
1. for any fixed $\varphi \in D$, $\psi_\pm^1$ satisfies (1.19) in $C([-U, U])$;
2. for any fixed $u \in [-U, U]$, $\psi_\pm^1$ is analytic in $\varphi \in D$;
3. if $-U \leq u \leq U$ and $\varphi \in D$, then

(3.23) $$\psi_\pm^1(u, \varphi) = \psi_\pm^0(u, \varphi) + O\left(\varepsilon\right),$$

the error estimate being uniform in $\varphi$ and $u$.

3.5.5. We shall need the asymptotics for the Wronskian $w_1\left(\varphi\right)$ of the solutions $\psi_+^1(x, \varphi)$ and $\psi_-^1(x, \varphi)$. Therefore we prove the

**Lemma 3.2.** *For $\varepsilon \to 0$,*

(3.24) $$w_1\left(\varphi\right) = w_0\left(\varphi\right) + \varepsilon\, W'(\varphi) \int_{-U_1}^{U_1} \psi_+^0(x, \varphi)\, \psi_-^0(x, \varphi)\, x\, dx + O\left(\varepsilon^2\right)$$

*uniformly in $\varphi \in D$.*

*Proof.* By (3.18),

(3.25) $$\begin{aligned} w_1 &= w\left(\psi_+^0 + \varepsilon\, R\, g_+,\, \psi_-^0 + \varepsilon\, R\, g_-\right) \\ &= w_0 + \varepsilon\left(w\left(\psi_+^0,\, R\, g_-\right) + w\left(R\, g_+,\, \psi_-^0\right)\right) + O\left(\varepsilon^2\right). \end{aligned}$$

Now, using the definition (3.16) of $R$, we obtain

$$w\left(\psi_+^0,\, R\, g_-\right) + w\left(R\, g_+,\, \psi_-^0\right) = -w_0 \int_{U_1}^{u} \psi_+^0\, g_-\, dx + w_0 \int_{-U_1}^{u} \psi_-^0\, g_+\, dx.$$

Substituting the expressions given by (3.22) for $g_\pm$ in this formula, we get

$$w\left(\psi_+^0,\, R\, g_-\right) + w\left(R\, g_+,\, \psi_-^0\right) = \int_{-U_1}^{U_1} \Delta\, \psi_+^0\, \psi_-^0\, dx + O\left(\varepsilon\right).$$

This formula and formula (3.25) imply the desired result as

$$\Delta\left(\varepsilon x,\, \varphi\right) = W'(\varphi)\, x + O\left(\varepsilon\right).$$

Note that, since $w_0$ has no zeros in $\overline{D}$, Lemma 3.2 implies in particular that, for sufficiently small $\varepsilon$, the newly constructed solutions are linearly independent for any $\varphi \in D$.

3.6. **Consistency condition.** The functions $\psi_\pm^1(u, \varphi)$ satisfy (1.19), but they need not form a consistent basis of solutions for the family of equations (1.19). However, in terms of these two solutions, we will construct two consistent ones.

Let $A(\varphi)$ and $B(\varphi)$ be analytic functions of $\varphi$. Define

$$(3.26) \qquad \psi(u, \varphi) = A(\varphi) \psi_+^1(u, \varphi) + B(\varphi) \psi_-^1(u, \varphi), \quad \varphi \in D,$$

Clearly, $\psi(u, \varphi)$ verifies (1.19). This solution satisfies the consistency condition (1.20) if and only if, for $\varphi, \varphi + \varepsilon \in D$, $u \in \mathbb{R}$, $A$ and $B$ satisfy

$$(3.27) \quad A(\varphi + \varepsilon) \psi_+^1(u - 1, \varphi + \varepsilon) + B(\varphi + \varepsilon) \psi_-^1(u - 1, \varphi + \varepsilon)$$
$$= A(\varphi) \psi_+^1(u, \varphi) + B(\varphi) \psi_-^1(u, \varphi).$$

This relation implies that

$$(3.28) \qquad \begin{pmatrix} A(\varphi + \varepsilon) \\ B(\varphi + \varepsilon) \end{pmatrix} = T(\varphi) \begin{pmatrix} A(\varphi) \\ B(\varphi) \end{pmatrix}, \quad \varphi, \varphi + \varepsilon \in D,$$

where $T$ is a $2 \times 2$-matrix with coefficients

$$(3.29)$$
$$T_{11} = \frac{w(\psi_+^1(\,.\,, \varphi), \psi_-^1(\,.\,-1, \varphi + \varepsilon))}{w_1(\varphi + \varepsilon)}, \quad T_{12} = \frac{w(\psi_-^1(\,.\,, \varphi), \psi_-^1(\,.\,-1, \varphi + \varepsilon))}{w_1(\varphi + \varepsilon)},$$
$$T_{21} = \frac{w(\psi_+^1(\,.\,-1, \varphi + \varepsilon), \psi_+^1(\,.\,, \varphi))}{w_1(\varphi + \varepsilon)}, \quad T_{22} = \frac{w(\psi_+^1(\,.\,-1, \varphi + \varepsilon), \psi_-^1(\,.\,, \varphi))}{w_1(\varphi + \varepsilon)}.$$

Now, our aim is to construct two analytic solutions of the difference equation (3.28) to get a consistent basis for the family of equations (1.19). Therefore, we shall first study the matrix $T$.

3.7. **Asymptotics for the coefficients of the matrix $T$.** When constructing the solutions of the difference equation (3.28), we shall need to know the diagonal elements of the matrix $T$ up to an error of size $O(\varepsilon^2)$, and the anti-diagonal elements up to an error of size $O(\varepsilon)$.

One has

**Proposition 3.1.** *For sufficiently small $\varepsilon$,*

$$T_{11} = \exp\left(i\kappa(\varphi + \varepsilon/2) + \frac{\varepsilon}{2} \frac{d}{d\varphi} \ln k_E + \varepsilon \omega_+^0(\varphi)\right) + O(\varepsilon^2),$$

$$T_{22} = \exp\left(-i\kappa(\varphi + \varepsilon/2) + \frac{\varepsilon}{2} \frac{d}{d\varphi} \ln k_E + \varepsilon \omega_-^0(\varphi)\right) + O(\varepsilon^2),$$

$$T_{12} = O(\varepsilon), \qquad T_{21} = O(\varepsilon),$$

where the error estimates are uniform in $\varphi \in D$. Here, the functions $\omega_\pm^0$ are defined by formulae (3.10), and $k_E = k'(E - W(\varphi)) = \left(\dfrac{d}{dE}k\right)(E - W(\varphi))$.

*Proof.* The estimates of $T_{12}$ and $T_{21}$ are elementary. For example, as $\psi_\pm^1$ are analytic in $\varphi \in D$,

(3.30)
$$T_{12} = \frac{w\left(\psi_-^1(\,.\,,\varphi),\,\psi_-^1(\,.\,-1,\varphi+\varepsilon)\right)}{w_1(\varphi+\varepsilon)} = \frac{w\left(\psi_-^1(\,.\,,\varphi),\,\psi_-^1(\,.\,-1,\varphi)\right)}{w_1(\varphi)} + O(\varepsilon),$$

and since
$$\psi_-^1(x-1,\varphi)) = e^{i\kappa(\varphi)}\psi_-^1(x,\varphi), \quad x \in \mathbb{R}, \quad \varphi \in D,$$
the leading term in (3.30) is zero.

To derive the asymptotics of $T_{11}$ and $T_{22}$, we shall use Lemma 2.1. We shall check only the formula for $T_{11}$. The second one is obtained in a similar way.

In view of (3.18),

$$\begin{aligned}
\psi_-^1(u-1,\varphi+\varepsilon) &= \psi_-^0(u-1,\varphi+\varepsilon) + \varepsilon\,(R\,g_-)(u-1,\varphi+\varepsilon) \\
&= e^{i\kappa(\varphi+\varepsilon)}\psi_-^0(u,\varphi+\varepsilon) + \varepsilon\,(R\,g_-)(u-1,\varphi+\varepsilon) \\
&= e^{i\kappa(\varphi)}\psi_-^0(u,\varphi) + \varepsilon\left(i\kappa'(\varphi)\,e^{i\kappa(\varphi)}\psi_-^0(u,\varphi) + \right. \\
&\qquad \left. + e^{i\kappa(\varphi)}\frac{\partial \psi_-^0}{\partial \varphi}(u,\varphi) + (R\,g_-)(u-1,\varphi)\right) + O(\varepsilon^2).
\end{aligned}$$

Now, again using (3.18) for $\psi_+^1$ and the result of the above calculation, we get

(3.31)
$$\begin{aligned}
w\left(\psi_+^1(\,.\,,\varphi),\psi_-^1(\,.\,-1,\varphi+\varepsilon)\right) &= e^{i\kappa}w_0 + \varepsilon\left(i\kappa'\,e^{i\kappa}w_0 + \right. \\
&\quad + e^{i\kappa}w\left(\psi_+^0(u),\frac{\partial \psi_-^0}{\partial \varphi}(u)\right) + w\left(Rg_+(u),\,e^{i\kappa}\psi_-^0(u)\right) + \\
&\quad \left. + w\left(\psi_+^0(u),\,Rg_-(u-1)\right)\right) + O(\varepsilon^2).
\end{aligned}$$

Here, for the sake of brevity, we do not write the argument $\varphi$ anymore. Elementary calculations show that

$$w\left(Rg_+(u),\,\psi_-^0(u)\right) = w_0 \int_{-U_1}^{u} \psi_-^0(x)\,g_+(x)\,dx,$$

$$w\left(\psi_+^0(u),\,Rg_-(u-1)\right) = -w_0 e^{i\kappa}\int_{U_1}^{u-1}\psi_+^0(x)\,g_-(x)\,dx,$$

and, using formula (3.22) as when proving Lemma 3.2, and the analyticity of $W$, we get

$$\begin{aligned}
w\left(Rg_+(u),\,e^{i\kappa}\psi_-^0(u)\right) + w\left(\psi_+^0(u),\,Rg_-(u-1)\right) &= \\
= e^{i\kappa}\,W'\left(\int_{-U_1}^{U_1}\psi_+^0(x)\,\psi_-^0(x)\,x\,dx + \int_{u-1}^{u}\psi_+^0(x)\,\psi_-^0(x)\,x\,dx\right) + O(\varepsilon).
\end{aligned}$$

Substituting this expression into (3.31) and taking into account Lemma 3.2, we obtain

(3.32)
$$w\left(\psi_+^1(\,.\,,\varphi),\psi_-^1(\,.\,-1,\varphi+\varepsilon)\right) = e^{i\kappa}\,w_1 + \varepsilon\,e^{i\kappa}\left(i\kappa'\,w_0 + \right.$$
$$\left. +w\left(\psi_+^0(u),\frac{\partial\psi_-^0}{\partial\varphi}(u)\right) + W'\int_{u-1}^u \psi_+^0(x)\,\psi_-^0(x)\,x\,dx\right) + O\left(\varepsilon^2\right).$$

But, in view of Corollary 2.1, we have

$$w\left(\psi_+^0(u),\frac{\partial\psi_-^0}{\partial\varphi}(u)\right) + W'\int_{u-1}^u \psi_+^0(x)\,\psi_-^0(x)\,x\,dx =$$
$$= w_0' + w_0\left(\frac{1}{2}\frac{d}{d\varphi}\ln k_E - \frac{i}{2}\kappa'(\varphi) + \omega_+^0\right).$$

Substituting this result in (3.32), one easily gets

$$w\left(\psi_+^1(\,.\,,\varphi),\,\psi_-^1(\,.\,-1,\varphi+\varepsilon)\right)$$
$$= e^{i\kappa}\left(w_1 + \varepsilon\left(i\kappa'w_0/2 + w_0' + w_0\left(\frac{1}{2}\frac{d}{d\varphi}\ln k_E + \omega_+^0\right)\right)\right) + O\left(\varepsilon^2\right)$$
$$= e^{i\kappa}\,w_1(\varphi+\varepsilon)\left(1 + \varepsilon\left(i\kappa'/2 + \frac{1}{2}\frac{d}{d\varphi}\ln k_E + \omega_+^0\right)\right) + O\left(\varepsilon^2\right)$$
$$= w_1(\varphi+\varepsilon)\,\exp\left(i\kappa\,(\varphi+\varepsilon/2) + \frac{\varepsilon}{2}\frac{d}{d\varphi}\ln k_E + \varepsilon\omega_+^0\right) + O\left(\varepsilon^2\right).$$

We have indicated explicitly the cases where the argument of the functions differs from $\varphi$. This, together with the definition for $T_{11}$ given by (3.29), implies the desired result. □

**3.8. Transformation of the difference equation into an integral one.** To construct analytic solutions of equation (3.28), we transform it into an integral equation. This transformation consists in three steps.

3.8.1. Note that, in view of Proposition 3.1, the matrix $T$ is close to a diagonal one. First, we transform the functions $A$ and $B$ so that the diagonal elements of the new matrix be close to 1. Therefore, we define

(3.33) $\qquad A(\varphi) = e^{i\theta_A(\varphi)}\,A_1(\varphi), \quad B(\varphi) = e^{i\theta_B(\varphi)}\,B_1(\varphi).$

So that, equation (3.28) takes the form

(3.34) $\qquad \begin{pmatrix} A_1(\varphi+\varepsilon) \\ B_1(\varphi+\varepsilon) \end{pmatrix} = \hat{T}(\varphi)\begin{pmatrix} A_1(\varphi) \\ B_1(\varphi) \end{pmatrix}, \quad \varphi,\ \varphi+\varepsilon \in D,$

with

$$\hat{T} = \begin{pmatrix} e^{-\theta_A(\varphi+\varepsilon)+\theta_A(\varphi)}\,T_{11} & e^{-\theta_A(\varphi+\varepsilon)+\theta_B(\varphi)}\,T_{12} \\ e^{-\theta_B(\varphi+\varepsilon)+\theta_A(\varphi)}\,T_{21} & e^{-\theta_B(\varphi+\varepsilon)+\theta_B(\varphi)}\,T_{22} \end{pmatrix}.$$

Fix a point $\varphi_0$ in $D$ and set

$$\theta_A(\varphi) = \frac{i}{\varepsilon} \int_{\varphi_0}^{\varphi} \kappa \, d\varphi + \frac{1}{2} \ln k_E + \int_{\varphi_0}^{\varphi} \omega_+^0 \, d\varphi \tag{3.35}$$

$$\theta_B(\varphi) = -\frac{i}{\varepsilon} \int_{\varphi_0}^{\varphi} \kappa \, d\varphi + \frac{1}{2} \ln k_E + \int_{\varphi_0}^{\varphi} \omega_-^0 \, d\varphi. \tag{3.36}$$

In this case, in view of Proposition 3.1,

$$\hat{T} = I + \hat{T}_1, \quad \hat{T}_1 = \begin{pmatrix} O(\varepsilon^2) & e^{-2\theta_0} O(\varepsilon) \\ e^{2\theta_0} O(\varepsilon) & O(\varepsilon^2) \end{pmatrix}, \tag{3.37}$$

where

$$\theta_0 = \frac{i}{\varepsilon} \int_{\varphi_0}^{\varphi} \kappa \, d\varphi, \tag{3.38}$$

and the symbol $O(\varepsilon^l)$ denotes functions satisfying the estimate

$$|O(\varepsilon^l)| \leq C \varepsilon^l, \quad \varphi \in D,$$

uniformly in $\varphi$.

3.8.2. A simple example of how to pass from a difference equation to an integral equation is given by the following simple

**Lemma 3.3.** *Let $f(z)$ be a function of $z \in \mathbb{C}$ analytic in the strip $Y_1 \leq \operatorname{Im} z \leq Y_2$, $Y_1, Y_2 \in \mathbb{R}$. Fix $z_1$ and $z_2$ so that $\operatorname{Im} z_1 = Y_1$, $\operatorname{Im} z_2 = Y_2$. Let $\gamma_z$ be a smooth vertical curve going from $z_1$ to $z_2$ via $z$. The function*

$$g(z) = \frac{1}{2i} \int_{\gamma_z} \cot \pi(\zeta - z - 0) f(\zeta) \, d\zeta$$

*is a solution of the equation*

$$g(z+1) - g(z) = f(z)$$

*analytic in the strip $Y_1 < \operatorname{Im} z < Y_2$.*

The proof of this lemma immediately follows from the residue theorem. For analogous statements see, for example, [8].

We will follow the idea of this lemma to transform equation (3.34) into an integral equation.

3.8.3. We are ready to make the last two steps of the desired transformation. However, to construct the solutions $f_\pm$ described in Theorem 3.1, one has to make these two steps in two different ways depending on whether one is dealing with $f_+$ or $f_-$. We shall concentrate on $f_+$. The solution $f_-$ can be constructed similarly and we shall discuss the details of this construction in the end of section 3.

Fix two points $\zeta_1$ and $\zeta_2$ in $D$ so that $\operatorname{Im} \zeta_1 < \operatorname{Im} \zeta_2$ and so that they can be connected by a vertical curve $\gamma \subset D$. For a function $f(\varphi)$ sufficiently regular on $\gamma$, define

$$L_+ f(\varphi) = \frac{1}{2i\varepsilon} \int_\gamma \left( \cot \left[ \frac{\pi(\varphi' - \varphi - 0)}{\varepsilon} \right] - i \right) f(\varphi') \, d\varphi'. \tag{3.39}$$

To construct a solution of (3.34) analytic in the domain
$$\{\varphi \in D;\ \operatorname{Im} \zeta_1 < \operatorname{Im} \varphi < \operatorname{Im} \zeta_2\},$$
we consider the equation

(3.40) $$\begin{pmatrix} A_1 \\ B_1 \end{pmatrix} = \begin{pmatrix} 1 \\ 0 \end{pmatrix} + L_+ \hat{T}_1 \begin{pmatrix} A_1 \\ B_1 \end{pmatrix}, \quad \varphi \in \gamma \subset D.$$

3.8.4. We need to make one more transformation. We fix $\alpha$ such that
$$0 < \alpha < 1,$$
and let

(3.41) $$B_1 = \varepsilon^\alpha e^{2\theta_0} B_2,$$

so that equation (3.40) be transformed into

(3.42) $$\begin{pmatrix} A_1 \\ B_2 \end{pmatrix} = \begin{pmatrix} 1 \\ 0 \end{pmatrix} + \begin{pmatrix} L_+ O(\varepsilon^2) & L_+ O(\varepsilon^{1+\alpha}) \\ K_+ O(\varepsilon^{1-\alpha}) & K_+ O(\varepsilon^2) \end{pmatrix} \begin{pmatrix} A_1 \\ B_2 \end{pmatrix},$$

where $\varphi \in \gamma \subset D$ and

(3.43) $$K_+ f(\varphi) = e^{-2\theta_0(\varphi)} L_+ \left( e^{2\theta_0(.)} f(.) \right)(\varphi).$$

In the next subsection, we shall show that, if $\gamma$ is a *strictly canonical curve*, then, after a proper choice of the functional space, the operator in the right hand side of (3.42) is small, and thus, that the equation (3.42) has a unique solution. Then we show that this solution can be analytically continued outside the contour $\gamma$. This will allow us to construct an analytic solution of the system (3.28).

We note that an operator similar to $K_+$ was studied in [7]. It was the integral operator with the kernel
$$\mathcal{L}(z, z') = \frac{1}{2ih} \left( \tan \frac{\pi(z'-z)}{h} - i \right) e^{\frac{i}{h} \int_z^{z'} p(\zeta) d\zeta},$$
acting in the space of functions analytic in a vicinity of a strictly vertical curve going from $-i\infty$ to $+i\infty$. One can immediately see that, if, along this curve, $\operatorname{Im} \int^z p\, d\zeta$ is strictly monotonically increasing and $\operatorname{Im} \int^z (p-\pi)\, d\zeta$ is strictly monotonically decreasing, then $|\mathcal{L}(z, z')| \leq \operatorname{Const} h^{-1} e^{-\delta|\operatorname{Im} z - \operatorname{Im} z'|/h}$, $\delta > 0$, and, thus, for small $h$, the corresponding operator is bounded. In our case, the estimates are more complicated since $K_+$ is a singular integral operator.

3.9. **Estimates of the norms of the integral operators.** In this subsection we estimate the norms of $L_+$ and $K_+$ as the operators acting in a proper functional space. We assume that the curve $\gamma$ defined above is sufficiently smooth. In the sequel, the letter $C$ will denote a positive constants independent of $\varepsilon$.

3.9.1. *Functional space.* Let $\gamma'$ be the curve $\gamma$ without the end points. Fix two positive numbers $a$ and $b$ so that

(3.44) $$0 < b < a < 1.$$

For a function $f$ defined on $\gamma'$, we let

(3.45) $$\|f\|_{a,b} = \sup_{\varphi \in \gamma'} |f_a(\varphi)| + \sup_{\varphi, \varphi' \in \gamma'} \frac{|f_a(\varphi) - f_a(\varphi')|}{|\varphi - \varphi'|^b},$$

where
$$f_a(\varphi) = |\varphi - \zeta_1|^a |\varphi - \zeta_2|^a f(\varphi).$$

We say that $f$ belongs to $C_{a,b}(\gamma)$ if $\|f\|_{a,b} < \infty$. One can easily see that

1. The set $C_{a,b}(\gamma)$ equipped with the norm (3.45) is a Banach space.
2. For any $f \in C^1(\gamma)$,

(3.46) $$\|f\|_{a,b} \leq C \|f\|_{C^1}$$

and

(3.47) $$\|f\|_{a,b} \leq C \|f\|_C + C \|f\|_C^{1-b} \|f'\|_C^b.$$

with some $C > 0$ independent of $f$.

We have chosen the space $C_{a,b}(\gamma)$ as, in particular, it is well known (see, for example, [19]) that the singular integral operator defined by the formula

(3.48) $$Sf(\varphi) = \frac{1}{2\pi i} \int_\gamma \frac{f(\varphi')}{\varphi' - \varphi - 0} d\varphi', \quad \varphi \in \gamma',$$

is a bounded operator in this space.

3.9.2. *Estimates of the norms.* The operators $L_+$ and $K_+$ are bounded operators acting in $C_{a,b}(\gamma)$. The estimate of the norm of $L_+$ does not require any additional condition on $\gamma$. The estimate of the norm of $K_+$ is the principal point in the proof of Theorem 3.1. It is in this estimate where one has to use the definition of the canonical curves. First, we shall describe the estimates of the norms of these two operators, and then, we shall prove them. The estimates given below are rather rough, but we do not need any sharper ones.

Begin with the operator $L_+$.

**Lemma 3.4.** *Let $\gamma$ be sufficiently smooth strictly vertical curve. The norm of the operator $L_+$ in $C_{a,b}(\gamma)$ satisfies the estimate*

(3.49) $$\|L_+\| \leq C \varepsilon^{-1-a}.$$

One can see that the norm of $L_+$ can be estimated from below be $C/\varepsilon$. At the same time, if $\gamma$ is a strictly canonical curve, then the norm of $K_+$ satisfies a much better estimate. We recall that $a$ and $b$ are two arbitrary positive numbers fixed so that $0 < b < a < 1$. One has

**Proposition 3.2.** *Let $\gamma$ be sufficiently smooth strictly canonical curve. The norm of the operator $K_+$ in $C_{a,b}(\gamma)$ satisfies the estimate*

(3.50) $$\|K_+\| \leq C \varepsilon^{-a-b}.$$

**Proof of Lemma 3.4.** For $f \in C_{a,b}(\gamma)$, represent $L_+ f$ in the form
$$L_+ f = R_L f + S f,$$
where $S$ is the singular integral operator (3.48) and
$$R_L f(\varphi) = \frac{1}{2i\varepsilon} \int_{\gamma'} r_L\left(\frac{\varphi' - \varphi}{\varepsilon}\right) f(\varphi') \, d\varphi',$$
where
$$r_L(\varphi) = \frac{1}{\varphi} \left( \varphi \left( \cot(\pi\varphi) - i \right) - \frac{1}{\pi} \right).$$
The norm of $S$ is bounded by a constant independent of $\varepsilon$. Let us estimate the norm of $R_L$.

Note that, along any given strictly vertical curve passing via the point $\varphi_0 = 0$,

(3.51) $$|r_L(\varphi)| \leq C,$$

and

(3.52) $$|r'_L(\varphi)| \leq C \begin{cases} 1 & |\eta| \leq 1, \\ e^{-2\pi|\eta|} + \frac{1}{|\eta|^2}, & |\eta| \geq 1, \end{cases} \quad \eta = \mathrm{Im}\, \varphi.$$

The constants here depend only on the parameter $\delta_1$ from the inequality (3.2) for a strictly vertical curve. Below, we shall denote $\mathrm{Im}\, \varphi$ by $\eta$.

First, we estimate $\|R_L f\|_C$. Fix $\varphi' \in \gamma$. Since the transformation $\varphi \mapsto \frac{\varphi - \varphi'}{\varepsilon}$ does not change the angles on the complex plane $\mathbb{C}$, the inequality (3.51) takes the form
$$\left| r_L\left(\frac{\varphi - \varphi'}{\varepsilon}\right) \right| \leq C, \quad \varphi, \varphi' \in \gamma,$$
with $C$ is independent of $\varepsilon$, $\varphi$ and $\varphi'$. This immediately implies that

(3.53) $$\|R_L f\|_C \leq C\,\varepsilon^{-1} \int_{\eta_1}^{\eta_2} \frac{d\eta'}{(\eta' - \eta_1)^a (\eta_2 - \eta')^a} \|f\|_{a,b} = C\,\varepsilon^{-1} \|f\|_{a,b}.$$

Here $\eta_{1,2} = \mathrm{Im}\, \zeta_{1,2}$.

Similarly, using (3.51), we get
$$|(R_L f)'(\varphi)| \leq C\,\varepsilon^{-2} \left( I_1(\eta) + I_2(\eta) + I_3(\eta) \right) \|f\|_{a,b}, \quad \eta = \mathrm{Im}\, \varphi,$$
where
$$I_1 = \int_{\eta_1,\, |\eta - \eta'| \leq \varepsilon}^{\eta_2} \frac{d\eta'}{(\eta' - \eta_1)^a (\eta_2 - \eta')^a},$$
$$I_2 = \varepsilon^2 \int_{\eta_1,\, |\eta - \eta'| \geq \varepsilon}^{\eta_2} \frac{d\eta'}{|\eta' - \eta|^2 (\eta' - \eta_1)^a (\eta_2 - \eta')^a},$$
$$I_3 = \int_{\eta_1,\, |\eta - \eta'| \geq \varepsilon}^{\eta_2} \frac{e^{-2\pi|\eta - \eta'|/\varepsilon} \, d\eta'}{(\eta' - \eta_1)^a (\eta_2 - \eta')^a},$$
Now, elementary estimates show that
$$\sup_{\eta_1 < \eta < \eta_2} |I_j| \leq C\,\varepsilon^{1-a}, \quad j = 1, 2, 3.$$

Thus,

(3.54) $$\|(R_L f)'\|_C \leq C\varepsilon^{-1-a}\|f\|_{a,b}.$$

In view of (3.46), (3.53) and (3.54) imply the desired result. □

**Proof of Proposition 3.2.** As in the proof of Lemma 3.4, we represent $K_+ f$ in the form
$$K_+ f = R_K f + S f,$$
where
$$R_K f(\varphi) = \frac{1}{2i\varepsilon} \int_{\gamma'} r_K(\varphi, \varphi') f(\varphi') d\varphi',$$
$$r_K(\varphi, \varphi') = \frac{\varepsilon}{\varphi' - \varphi} \left( \frac{(\varphi' - \varphi)}{\varepsilon} e(\varphi, \varphi') - \frac{1}{\pi} \right),$$
$$e(\varphi, \varphi') = \left( \cot \frac{\pi(\varphi' - \varphi)}{\varepsilon} - i \right) e^{\frac{2i}{\varepsilon} \int_\varphi^{\varphi'} \kappa(s) ds}.$$

We begin by estimating $\|R_K f\|_C$. As before, let $\eta = \operatorname{Im} \varphi$. Note that
$$\left| \int_\varphi^{\varphi'} \kappa(s) ds \right| \leq C |\varphi - \varphi'|, \quad \varphi, \varphi' \in \gamma,$$

Therefore,
$$e(\varphi, \varphi') = \frac{\varepsilon}{\varphi - \varphi'} \left( \frac{1}{\pi} + \frac{1}{\varepsilon} O(\varphi - \varphi') \right), \quad \varphi, \varphi' \in \gamma, \quad |\eta - \eta'| \leq \varepsilon,$$

where the error estimate is uniform in $\varepsilon$, $\varphi$ and $\varphi'$. So,

(3.55) $$|r_K(\varphi, \varphi')| \leq C, \quad |\eta - \eta'| \leq \varepsilon, \quad \varphi, \varphi' \in \gamma.$$

Now, we assume that $\gamma$ is strictly canonical. Due to the definition of a strictly canonical curve, we show that

(3.56) $$|e(\varphi, \varphi')| \leq C e^{-\frac{\delta}{\varepsilon} |\eta - \eta'|}, \quad |\eta - \eta'| \geq \varepsilon, \quad \varphi, \varphi' \in \gamma,$$

where $\delta$ is a positive constant independent of $\varepsilon$.

First, we consider the case where $\eta - \eta' \geq \varepsilon$. In this case, since $\gamma$ is strictly vertical,

(3.57) $$\left| \cot \frac{\pi(\varphi' - \varphi)}{\varepsilon} - i \right| \leq C e^{-\frac{2\pi}{\varepsilon}(\eta - \eta')},$$

and, therefore,
$$|e(\varphi, \varphi')| \leq C e^{\frac{2}{\varepsilon} \operatorname{Im} \int_{\varphi'}^\varphi (\kappa(s) - \pi) ds}, \quad \varphi, \varphi' \in \gamma.$$

In view of the definition of a strictly canonical curve,
$$\operatorname{Im} \int_{\varphi'}^\varphi (\kappa(s) - \pi) ds \leq -\delta_2 (\eta - \eta'),$$

with some positive constant $\delta_2$ (see (3.3)). This implies (3.56).

Now, we consider the case where $\eta - \eta' \leq -\varepsilon$. In this case, instead of (3.57), we get
$$\left| \cot \frac{\pi(\varphi' - \varphi)}{\varepsilon} - i \right| \leq C,$$
and, therefore,
$$|e(\varphi, \varphi')| \leq C\, e^{-\frac{2}{\varepsilon} \operatorname{Im} \int_\varphi^{\varphi'} \kappa(s)\,ds}, \qquad \varphi, \varphi' \in \gamma.$$
But, since $\gamma$ is a strictly canonical curve,
$$\operatorname{Im} \int_\varphi^{\varphi'} \kappa(s)\,ds \geq \delta(\eta' - \eta),$$
with some positive $\delta$, which again leads to (3.56).

The estimate (3.56) and the definition of $r_K$ imply the inequality
$$(3.58) \quad |r_K(\varphi, \varphi')| \leq C\, e^{-\frac{\delta}{\varepsilon}|\eta - \eta'|} + \frac{C\varepsilon}{|\eta - \eta'|}, \qquad |\eta - \eta'| \geq \varepsilon, \qquad \varphi, \varphi' \in \gamma.$$

The inequalities (3.55), (3.58) and the definition of $R_K$ imply that
$$|R_K f(\varphi)| \leq C\, (I_4(\eta) + I_5(\eta) + I_6(\eta))\, \|f\|_{a,b}, \quad \varphi \in \gamma,$$
where
$$I_4 = \frac{1}{\varepsilon} \int_{\eta_1,\, |\eta - \eta'| \leq \varepsilon}^{\eta_2} \frac{d\eta'}{(\eta' - \eta_1)^a (\eta_2 - \eta')^a},$$
$$I_5 = \int_{\eta_1,\, |\eta - \eta'| \geq \varepsilon}^{\eta_2} \frac{d\eta'}{|\eta' - \eta| (\eta' - \eta_1)^a (\eta_2 - \eta')^a},$$
$$I_6 = \frac{1}{\varepsilon} \int_{\eta_1,\, |\eta - \eta'| \geq \varepsilon}^{\eta_2} \frac{e^{-\delta|\eta - \eta'|/\varepsilon}\, d\eta'}{(\eta' - \eta_1)^a (\eta_2 - \eta')^a},$$

Elementary estimates show that
$$\sup_{\eta_1 < \eta < \eta_2} |I_j| \leq C\, \varepsilon^{-a}, \quad j = 4, 5, 6.$$

And, as a result, we get
$$(3.59) \qquad \|R_K f\|_C \leq C\, \varepsilon^{-a} \|f\|_{a,b}.$$

Having estimated $\|R_K f\|_C$, we estimate $\|(R_K f)'\|_C$. We first show that, for $\varphi, \varphi' \in \gamma$,

$$(3.60) \quad \left| \frac{\partial r_K}{\partial \varphi}(\varphi, \varphi') \right| \leq C\, \varepsilon^{-1} \begin{cases} 1 & \text{when } |\eta - \eta'| \leq \varepsilon, \\ e^{-2\pi|\eta - \eta'|/\varepsilon} + \frac{\varepsilon^2}{|\eta - \eta'|^2}, & \text{when } |\eta - \eta'| \geq \varepsilon, \end{cases}$$

Note that
(3.61)
$$\frac{\partial r_K}{\partial \varphi}(\varphi, \varphi') = \left( \frac{\pi}{\varepsilon} \left( \sin^2 \zeta\, (\cot \zeta - i) \right)^{-1} - \frac{2i}{\varepsilon} \kappa(\varphi) \right) e(\varphi, \varphi') - \frac{\varepsilon}{\pi(\varphi - \varphi')^2},$$
where
$$\zeta = \frac{\pi(\varphi' - \varphi)}{\varepsilon}.$$

Now, first, we assume that $|\eta-\eta'|\leq\varepsilon$. We recall that $\gamma\subset D$, and, thus, the curve $\gamma$ is situated at some distance from the branching points of the quasi-momentum $\kappa$. Therefore,

$$\int_{\varphi}^{\varphi'}\kappa(s)\,ds=(\varphi'-\varphi)\kappa(\varphi)+O((\varphi-\varphi')^2),\quad \varphi,\varphi'\in\gamma,$$

with the error estimate uniform in $\varphi'$ and $\varepsilon$. This allows us to get

$$e(\varphi,\varphi')=\frac{\varepsilon}{\pi(\varphi'-\varphi)}\left(1-i\frac{\pi(\varphi'-\varphi)}{\varepsilon}+2i\kappa(\varphi)\frac{(\varphi'-\varphi)}{\varepsilon}+O\left(\frac{(\varphi-\varphi')^2}{\varepsilon^2}\right)\right),$$
$$\varphi,\varphi'\in\gamma,\quad |\eta-\eta'|\leq\varepsilon.$$

And, substituting this representation in (3.61), one easily comes to the first inequality of (3.60).

To treat the case where $|\eta-\eta'|\geq\varepsilon$, one notes that, for $|\eta-\eta'|\geq\varepsilon$ and $\varphi,\varphi'\in\gamma$, one has

$$\left|\sin^2\zeta\,(\cot\zeta-i)\right|^{-1}\leq C.$$

This estimate and the estimate (3.56) lead to the second inequality from (3.60).

Using the estimates (3.60), one obtains

$$|(R_K f)'(\varphi)|\leq C\,\varepsilon^{-2}\,(\,I_1(\eta)+I_2(\eta)+I_3(\eta)\,)\,\|f\|_{a,b},\quad \varphi\in\gamma,$$

with the same $I_1$, $I_2$ and $I_3$ as in the proof of Lemma 3.4. This gives the inequality

(3.62) $$\|(R_K f)'\|_C\leq C\,\varepsilon^{-1-a}.$$

similar to (3.54). Now, combining (3.59) and (3.62), by means of (3.47), we get

$$\|R_K f\|_{a,b}\leq C\left(\varepsilon^{-a}+\varepsilon^{-a(1-b)}\,\varepsilon^{-(1+a)b}\right)\|f\|_{a,b}\leq C\,\varepsilon^{-a-b}\|f\|_{a,b}.$$

This implies the announced estimate of $\|K_+\|$. □

3.10. **Analytic solution of equation (3.28).** We are now ready to solve equation (3.28). In this subsection, we prove

**Proposition 3.3.** *Let $K$ be a strictly canonical domain related to some points $\zeta_1,\zeta_2\in\mathbb{C}$, $\operatorname{Im}\zeta_1<\operatorname{Im}\zeta_2$. Let $\mathcal{C}$ be an admissible sub-domain of $K$. Assume that $\mathcal{C}$ does not contain any pole of $\psi_\pm^0$. For sufficiently small $\varepsilon>0$, there exists a solution of (3.28) analytic in the strip*

(3.63) $$\operatorname{Im}\zeta_1<\operatorname{Im}\varphi<\operatorname{Im}\zeta_2.$$

*This solution has the asymptotics:*

(3.64) $$\begin{pmatrix}A(\varphi)\\ B(\varphi)\end{pmatrix}=e^{\theta_A(\varphi)}\begin{pmatrix}(1+O(\varepsilon^\beta)),\\ O(\varepsilon^\beta)\end{pmatrix},\quad \varphi\in\mathcal{C},\quad \varepsilon\to 0.$$

*Here, $\beta$ is any fixed number,*

$$0<\beta<1,$$

*and $\theta_A$ is the function defined by (3.35) The error estimate is uniform in $\varphi\in\mathcal{C}$.*

The rest of this subsection will be devoted to the proof of Proposition 3.3.

3.10.1. Let the domain $D$ in all the previous constructions coincide with $K$, and let $\gamma \in K$ be one of the strictly canonical curves connecting $\zeta_1$ and $\zeta_2$. First, construct a solution of (3.42) with the components $A_1$ and $B_2$ from the class $C_{a,b}(\gamma)$. In view of Lemma 3.4 and Proposition 3.2, for sufficiently small $\varepsilon$, the components of the matrix integral operator from the right hand side of equation (3.42) satisfy the estimates

$$
(3.65) \qquad \|L_+ \left(O\left(\varepsilon^2\right) \cdot \right)\| \leq C\,\varepsilon^{1-a}, \qquad \|L_+ \left(O\left(\varepsilon^{1+\alpha}\right) \cdot \right)\| \leq C\,\varepsilon^{\alpha-a},
$$
$$
\|K_+ \left(O\left(\varepsilon^{1-\alpha}\right) \cdot \right)\| \leq C\,\varepsilon^{1-\alpha-a-b}, \quad \|K_+ \left(O\left(\varepsilon^2\right) \cdot \right)\| \leq C\,\varepsilon^{2-a-b},
$$

where $\|\cdot\|$ denotes the operator norm induced by $\|\cdot\|_{a,b}$. Fix the numbers $\alpha$, $a$ and $b$ so that

$$(3.66) \qquad 0 < b < a < \alpha < 1, \quad \alpha + a + b < 1.$$

Clearly, in this case, for sufficiently small $\varepsilon$, equation (3.42) has a unique solution in $C_{a,b}$. This solution satisfies

$$A_1 = 1 + O\left(\varepsilon^{\alpha-a}\right), \quad B_2 = O\left(\varepsilon^{1-\alpha-a-b}\right).$$

Here the error terms are estimated in the $\|\cdot\|_{a,b}$-norm. Letting

$$(3.67) \qquad \beta_1 = \min\{\alpha - a,\ 1 - a - b\}$$

we get

$$(3.68) \qquad A_1 = 1 + O\left(\varepsilon^{\beta_1}\right), \quad B_2 = O\left(\varepsilon^{\beta_1-\alpha}\right).$$

Note that $0 < \beta_1 < 1$.

3.10.2. Now, let us check that $A_1$ and $B_2$ can be analytically continued in the strip defined by (3.63). Clearly, the poles in $\varphi'$ of the kernel of the operator $L_+$ are situated at the points $\varphi + 0 + \varepsilon l$, $l \in \mathbb{Z}$. Therefore, the right hand side in (3.42) can be analytically continued in the strip

$$S_\varepsilon = \{\varphi \in \mathbb{C} : \exists \zeta \in \gamma',\ 0 < \operatorname{Re}\varphi - \operatorname{Re}\zeta < \varepsilon\},$$

where, as before, $\gamma'$ is the curve $\gamma$ without the endpoints. Deformating the integration contour in (3.39) inside $S_\varepsilon$, we see that, in fact, $A_1$ and $B_2$ can be analytically continued in the strip $S_{2\varepsilon}$. Repeating these arguments, we prove that these functions can be continued in the part of the strip defined by (3.63) situated at the right of $\gamma$.

We check that $A_1$ and $B_2$ are analytic at the left of $\gamma$. Denote the second term in (3.42) by $\mathcal{L} \begin{pmatrix} A_1 \\ B_2 \end{pmatrix}(\varphi)$. We assume first that $\varphi$ is located at the right of $\gamma$ and deform the integration contour so that it passes between the points $\varphi$ and $\varphi + h$. Since the kernels of $L_+$ and $K_+$ have simple poles at $\varphi' = \varphi$, instead of (3.42), one obtains

$$(3.69) \qquad [I + O\left(\varepsilon^{1-\alpha}\right)] \begin{pmatrix} A_1(\varphi) \\ B_2(\varphi) \end{pmatrix} = \begin{pmatrix} 1 \\ 0 \end{pmatrix} + \tilde{\mathcal{L}} \begin{pmatrix} A_1 \\ B_2 \end{pmatrix}(\varphi).$$

Here, the second term in the brackets in the left hand side is a $2 \times 2$-matrix with analytic coefficients which can be estimated by $C\varepsilon^{1-\alpha}$ for $\varphi \in K$. This term is the residue that appeared in the course of the contour deformation. We have continued $A_1$ and $B_2$ to the right of $\gamma$ by using (3.42). Applying almost the same arguments to (3.69), one continues these functions to the left of $\gamma$.

3.10.3. Now, we discuss analytic properties of the functions $A_1$ and $B_2$ in some small vicinities of the end points of the curve $\gamma$. Clearly, for $\varphi \in \gamma$ being in a sufficiently small vicinity $V_1$ of $\zeta_1$, equation (3.42) can be rewritten in the form

$$\begin{pmatrix} A_1 \\ B_2 \end{pmatrix} = \begin{pmatrix} 1 \\ 0 \end{pmatrix} + S \begin{pmatrix} A_1 \\ B_2 \end{pmatrix} + \mathcal{L}_1 \begin{pmatrix} A_1 \\ B_2 \end{pmatrix},$$

where $S$ is the singular integral operator (3.48), and $\mathcal{L}_1$ is a matrix integral operator with some kernel analytic in $\varphi \in V_1$ and in $\varphi' \in \gamma$. Therefore,

$$\begin{pmatrix} A_1 \\ B_2 \end{pmatrix}(\varphi) = S \begin{pmatrix} A_1 \\ B_2 \end{pmatrix}(\varphi) + f(\varphi),$$

where $f$ is analytic in $V_1$. We recall that $A_1$ and $B_2$ belong to $C_{a,b}(\gamma)$. The analytic behavior of the singular integral in a vicinity of the ends of the integration contour is well known (see, e.g, [19]). Using classical results, we see that $\zeta_1$ is a branching point of $A_1$ and $B_2$, and that these functions admit the representations

(3.70)
$$(\varphi - \zeta_1)^a A_1(\varphi) = c_1 + O(|\varphi - \zeta_1|^{b_1}), \quad (\varphi - \zeta_1)^a B_2(\varphi) = c_2 + O(|\varphi - \zeta_1|^{b_1}),$$

where $c_{1,2}$ are some complex constants, and $b_1$ is a positive number. Similarly, in a sufficiently small vicinity $V_2$ of $\zeta_2$,

(3.71)
$$(\varphi - \zeta_2)^a A_1(\varphi) = c_3 + O(|\varphi - \zeta_2|^{b_2}), \quad (\varphi - \zeta_2)^a B_2(\varphi) = c_4 + O(|\varphi - \zeta_2|^{b_2}).$$

We can assume that $0 < b_1 = b_2 \leq b$. We set

(3.72) $$\beta = \min\{\alpha - a, \, 1 - a - b_1\}.$$

Clearly, $0 < \beta < 1$.

3.10.4. Consider the functions $A_1$ and $B_2$ along some other strictly canonical curve $\gamma_1 \in K$. We show that, for sufficiently small $\varepsilon$, they have the asymptotic representation

(3.73) $$A_1 = 1 + O(\varepsilon^\beta), \quad B_2 = O(\varepsilon^{\beta-\alpha}) \quad \text{in} \quad C_{a,b_1}(\gamma_1)$$

similar to (3.68)

Therefore, instead of (3.42), we consider the analogous equation along the curve $\gamma_1$. For sufficiently small $\varepsilon$, this new equation has a unique solution in $C_{a,b_1}(\gamma_1)$, and the new functions $A_1$ and $B_2$, e.g. the components of this solution, do have the desired asymptotics. But, the representations (3.70) – (3.71) for the "old" functions $A_1$ and $B_2$ imply that they belong to this functional space. And just deforming the curve $\gamma$ to $\gamma_1$, one can see these two "old" functions satisfy the "new" equation. This implies the desired result.

3.10.5. Consider an admissible domain $\mathcal{C}$ of $K$. Now, we can easily show that

(3.74) $\qquad A_1(\varphi) = 1 + O(\varepsilon^\beta), \quad B_2(\varphi) = O(\varepsilon^{\beta-\alpha}), \quad \varphi \in \mathcal{C},$

with the uniform error estimates.

Let $\varphi_0 \in \overline{\mathcal{C}}$. One can consider $\varphi_0$ as an internal point of a domain $K_0 \subset K$ with the boundary $\Gamma_0$ consisting of two segments of two strictly canonical curves $\gamma_1$ and $\gamma_2$ situated in $K$. We represent the functions $A_1$ and $B_2$ by the Cauchy integrals along $\Gamma_0$. Along $\Gamma_0 \cap \gamma_1$ (resp. $\Gamma_0 \cap \gamma_2$), the functions $A_1$ and $B_2$ admit the representation (3.74) in $C_{a,\,b_1}(\gamma_1)$ (resp. $C_{a,\,b_1}(\gamma_2)$). This leads to the desired asymptotic formulae for any fixed $\varphi \in K_0$. Moreover, the error estimates are uniform in $\varphi$ in a sufficiently small fixed vicinity of $\varphi_0$.

Thus, we see that any point $\varphi_0$ of the compact $\overline{\mathcal{C}}$ has an open vicinity $V_0$ such that the desired asymptotic formulae are uniform in $\varphi \in V_0$. This implies the needed result.

3.10.6. To finish the proof, we reconstruct in terms of $B_2$ the function $B_1$ by formula (3.41). The vector $\begin{pmatrix} A_1 \\ B_1 \end{pmatrix}(\varphi)$ satisfies equation (3.40). By means of the residue theorem, one checks that

$$L_+\hat{T}_1 \begin{pmatrix} A_1 \\ B_1 \end{pmatrix}(\varphi + \varepsilon) - L_+\hat{T}_1 \begin{pmatrix} A_1 \\ B_1 \end{pmatrix}(\varphi) = \hat{T}_1 \begin{pmatrix} A_1 \\ B_1 \end{pmatrix}(\varphi),$$

$$\operatorname{Im}\zeta_1 < \operatorname{Im}\varphi < \operatorname{Im}\zeta_2.$$

This implies that this vector satisfies (3.34), and, by means of (3.33) we reconstruct a vector solution $\begin{pmatrix} A \\ B \end{pmatrix}$ of the difference equation (3.28). Clearly, this solution is analytic in $K$; moreover the asymptotic formulae (3.74) and the relations (3.33), (3.41) imply (3.64). $\square$

3.11. **Consistent basis solutions.** In this subsection, we complete the proof of Theorem 3.1.

3.11.1. First, we construct the solution $f_+$ described in Theorem 3.1. Denote the components of the vector solution of (3.28) described in Proposition 3.3 by $A^+$ and $B^+$. Let $\psi_\pm^1$ be the solutions of (1.19) constructed in subsection 3.5. We consider the function

(3.75) $\qquad f_+(u, \varphi) = A^+(\varphi)\psi_+^1(u, \varphi) + B^+(\varphi)\psi_-^1(u, \varphi),$

where $\operatorname{Im}\zeta_1 < \operatorname{Im}\varphi < \operatorname{Im}\zeta_2$, $u \in \mathbb{R}$. This function is a solution of (1.19) analytic in $\varphi \in \mathcal{C}$. By construction, it satisfies the condition (1.20)
Finally, we note that the asymptotics (3.64) and (3.23) imply

$$f_+(u, \varphi) = e^{\theta_A(\varphi)}\left(\psi_+^0(u, \varphi) + O(\varepsilon^\beta)\right), \quad \varphi \in \mathcal{C}, \quad -U \le u \le U.$$

Remembering the definitions of $\theta_A$, $\psi_\pm^0$ and $\Psi_\pm^0$, and equation (3.35), we see that this formula coincides with (3.11).

3.11.2. Up to this moment, we were assuming that $\mathcal{C}$ does not contain any pole of $\psi_\pm^0$. Now, we get formula (3.11), and see that the leading term in this formula is the canonical Bloch solution. By Lemma 3.1, it has no poles. Having found the canonical Bloch solutions, we could have used them everywhere instead of $\psi_\pm^0$. This remark enables us to remove the assumption about the absence of the poles of $\psi_\pm^0$ in $\mathcal{C}$.

3.11.3. Let us discuss the construction of the second solution of (1.19) described in Theorem 3.1. It is similar to the one leading to $f_+$. So, we shall only single out its new elements.

Now, the central objects of the analysis are the operators

$$L_- f(\varphi) = \frac{1}{2i\varepsilon} \int_\gamma \left( \cot \frac{\pi(\varphi' - \varphi - 0)}{\varepsilon} + i \right) f(\varphi') \, d\varphi',$$

and

$$K_- f(\varphi) = e^{2\theta_0(\varphi)} L_- \left( e^{-2\theta_0(.)} f(.) \right)(\varphi).$$

acting on the "old" space $C_{a,b}(\gamma)$. As in subsection 3.9.2, one checks that, for $\gamma$ being a strictly canonical curve,

$$\|L_-\| \leq C \, \varepsilon^{-1-a}, \quad \|K_-\| \leq C \, \varepsilon^{-a-b}.$$

To reduce the analysis of (3.34) to the analysis of these two operators, we consider instead of (3.40) the equation

$$\begin{pmatrix} A_1 \\ B_1 \end{pmatrix} = \begin{pmatrix} 0 \\ 1 \end{pmatrix} + L_- \hat{T}_1 \begin{pmatrix} A_1 \\ B_1 \end{pmatrix}, \quad \varphi \in \gamma \subset K.$$

Note that it differs from (3.40) not only by replacing $L_+$ by $L_-$, but also by replacing the term $\begin{pmatrix} 1 \\ 0 \end{pmatrix}$ by $\begin{pmatrix} 0 \\ 1 \end{pmatrix}$. The latter allows to obtain, by means of the transformation

$$A_1 = \varepsilon^\alpha e^{-2\theta_0} A_2,$$

the integral equation

$$\begin{pmatrix} A_2 \\ B_1 \end{pmatrix} = \begin{pmatrix} 0 \\ 1 \end{pmatrix} + \begin{pmatrix} K_- O(\varepsilon^2) & K_- O(\varepsilon^{1-\alpha}) \\ L_- O(\varepsilon^{1+\alpha}) & L_- O(\varepsilon^2) \end{pmatrix} \begin{pmatrix} A_2 \\ B_1 \end{pmatrix},$$

for $\varphi \in \gamma \subset K$.

The analysis of this equation is similar to the analysis of (3.42), and, as a result, one constructs the solution $f_-$ of (1.19) described in Theorem 3.1. Omitting the details, we note only that

(3.76) $$f_-(u, \varphi) = A^-(\varphi) \psi_+^1(u, \varphi) + B^-(\varphi) \psi_-^1(u, \varphi),$$

for $\operatorname{Im} \zeta_1 < \operatorname{Im} \varphi < \operatorname{Im} \zeta_2$, $u \in \mathbb{R}$.

The functions $A^-$ and $B^-$ are analytic in the strip defined by (3.63) and admit the asymptotic representations

(3.77) $$A^-(\varphi) = e^{\theta_B(\varphi)} O(\varepsilon^\beta), \quad B^-(\varphi) = e^{\theta_B(\varphi)} (1 + O(\varepsilon^\beta)),$$

for $\varphi \in \mathcal{C}$ when $\varepsilon \to 0$, with the error estimate being uniform. Here, $\beta$ is any fixed number, $0 < \beta < 1$, and $\theta_B$ is the function defined by (3.36). $\square$

3.11.4. Finally, we check that the Wronskian of $f_\pm(u, \varphi)$ is independent of $\varphi$, and, thus, that they do form a consistent basis of solutions for the equation family (1.19). By (3.75) and (3.76),

$$w\{f_+, f_-\} = \left(A^+(\varphi)B^-(\varphi) - A^-(\varphi)B^+(\varphi)\right) w_1(\varphi).$$

Here $w_1(\varphi)$ is the Wronskian of the solutions $\psi_\pm^1$. Now, using the asymptotics of $A^+$ and $B^+$ described in Proposition 3.3, the asymptotics of $A^-$ and $B^-$ given by (3.77), and the asymptotics of $w_1$ (see Lemma 3.2), we get

$$(3.78) \qquad w\{f_+, f_-\} = w_0(\varphi) \exp\left(\ln k_E' + \int_{\varphi_0}^{\varphi} \left(\omega_0^+ + \omega_0^-\right) d\varphi\right)(1 + o(1)).$$

By (3.6)

$$\omega_0^+ + \omega_0^- = -\frac{\partial}{\partial \varphi} \ln \int_0^1 \psi_+^0(x, \varphi)\, \psi_-^0(x, \varphi))\, dx.$$

Using Lemma 3.2, we transform this expression to the form

$$\omega_0^+ + \omega_0^- = -\frac{\partial}{\partial \varphi} \ln\left(k_E'\, w_0(\varphi)\right).$$

This allows to get finally

$$w\{f_+, f_-\} = k_E'\, w_0|_{\varphi = \varphi_0} (1 + g(\varphi)), \quad g(\varphi) = o(1).$$

This formula shows only that the leading terms of the asymptotics of the Wronskian is independent of $\varphi$.

Let us discuss the term $g(\varphi)$ in the last formula. Since $f_\pm$ satisfy (1.20), it is $\varepsilon$-periodic in $\varphi$. Moreover, its estimate is uniform in $\varphi \in \mathcal{C}$.

Let us change the notations, denoting by $f_+$ the "old" function $f_+$ divided by $1 + g$. The "new" functions $f_\pm$, clearly, have all the properties announced in Theorem 3.1. This completes the proof of the theorem. □

## 4. Main constructions of the complex WKB method for equation (1.24)

Here we consider the family of equations (1.19) with $W(\varphi) = \cos \varphi$.

### 4.1. Complex momentum.

4.1.1. We begin by describing the Riemann surface $\Gamma$ of the complex momentum. We recall that the conformal properties of $\kappa$ and the structure of its Riemann surface depend on the spectral parameter $E$ and on the lengths of the spectral zones and the spectral gaps of the operator $H_0 = -\dfrac{d^2}{dx^2} + V$. For the sake of definiteness, we shall work under the condition (1.26)

The branching points of the complex momentum $\kappa$ are described in formula (3.1). As $\cos \varphi$ is periodic, each of the $E_l$ gives rise to an infinite lattice of branching points of the complex momentum. There are two branching points corresponding to each of $E_l$ and situated in the strip $0 \leq \mathrm{Re}\,\varphi < 2\pi$. Let us describe these branching points.

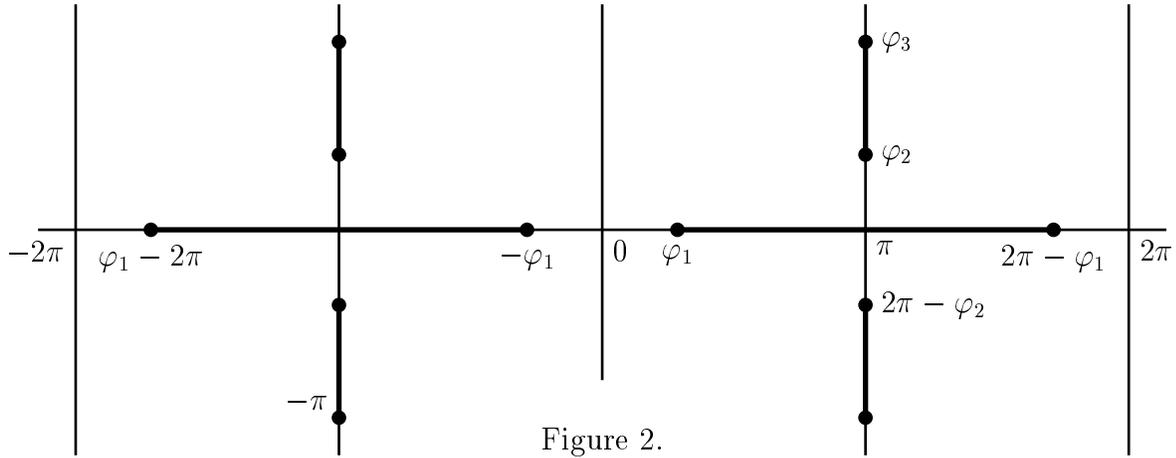

Figure 2.

The branching points corresponding to $E_1$ are real. They are symmetric with respect to the point $\pi$. We denote them by $\varphi_1$ and $2\pi - \varphi_1$ where $0 < \varphi_1 < \pi$.

The branching points corresponding to $E_l$, $l > 1$, are situated on the line $\operatorname{Re}\varphi = \pi$. They are also symmetric with respect to the point $\pi$. The points with positive imaginary part will be denoted by $\varphi_l$. Let $\eta_l = \operatorname{Im}\varphi_l > 0$.

The whole picture of the branching points is $2\pi$-periodic. It is partly shown on Fig. 2.

The conformal properties of the complex momentum near the branching points are determined by the behavior of the Bloch quasimomentum $k$ near the points $E_l$ (see subsection 2.1). In particular, the branching points $\varphi_1$ and $2\pi - \varphi_1$ related to $E_1$ are of square root type.

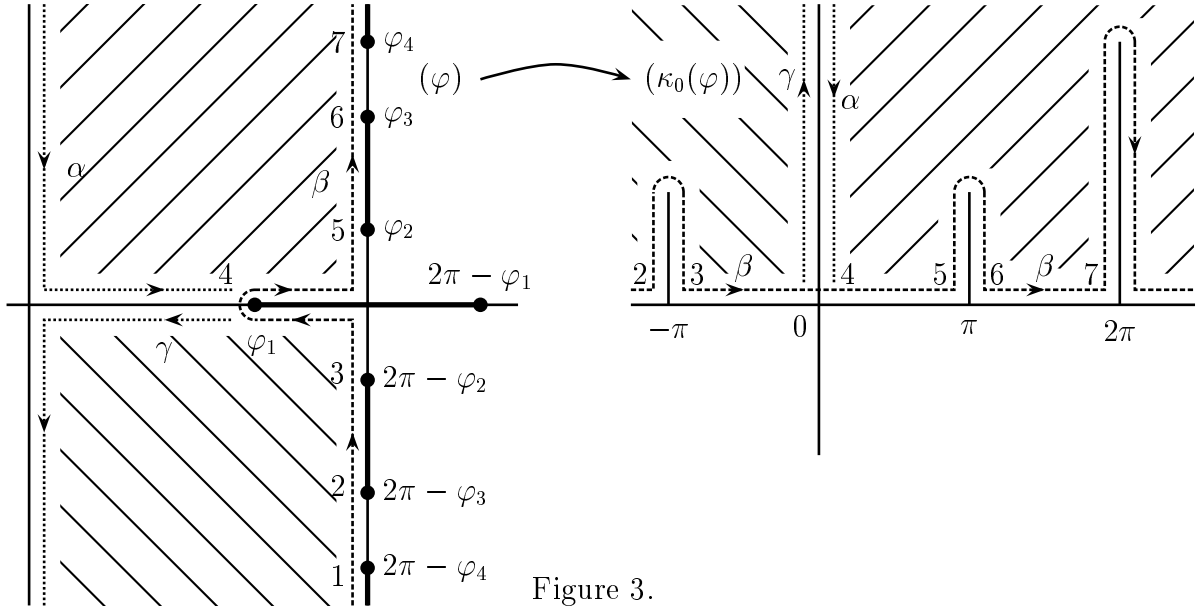

Figure 3.

4.1.2. One can fix a single valued branch of the complex momentum in a vicinity of $\varphi = 0$ by the conditions

(4.1) $$\operatorname{Im}\kappa_0(0) > 0, \quad \operatorname{Re}\kappa_0(0) = 0.$$

All the other branches defined there are then related to $\kappa_0$ by the formulae:

(4.2) $$\kappa_m^\pm(\varphi) = \pm\kappa_0(\varphi) + 2\pi m, \quad m \in \mathbf{Z}.$$

The branch $\kappa_0$ can be analytically continued along any curve not passing through any branching point. In Fig. 3, we draw several curves on the complex plane $\varphi \in \mathbb{C}$, and then, we show the images of these curves by $\kappa_0(\varphi)$. The results shown in these figures easily follow from the discussions of subsection 3.2.

Note that on the shadowed domain in Fig. 3, we have

$$\kappa_0(\varphi) = k_0(E - \cos\varphi),$$

where $k_0$ is the branch of the Bloch quasi-momentum defined in section 2.2.1. To get the following four pictures, one uses known elementary properties of the conformal mappings. We shall use Fig. 4 for drawing the cuts on the complex plane needed to define the branch $\kappa_0$ globally.

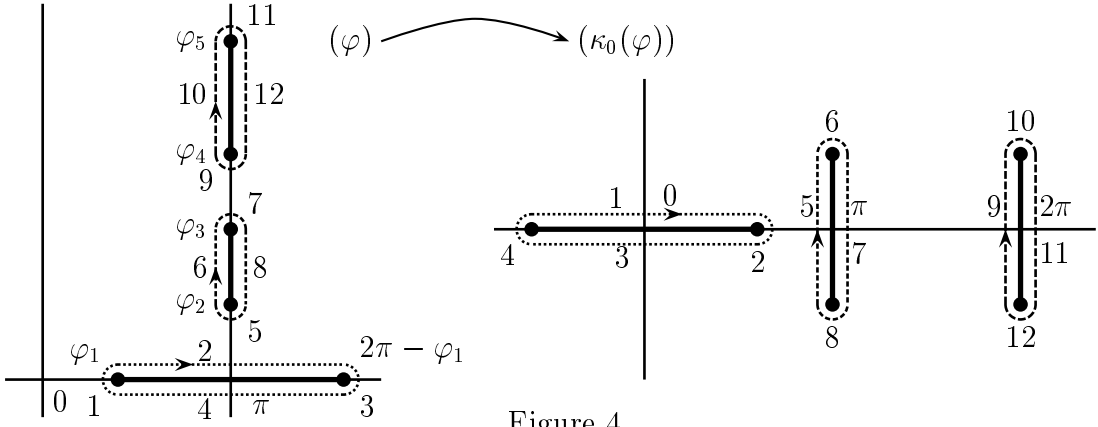

Figure 4.

Cut the complex plane $\varphi \in \mathbb{C}$ as in Fig. 2, e.g., first draw the horizontal cut between the points $\varphi_1$ and $2\pi - \varphi_1$, the vertical cuts between the pairs of the points $\varphi_{2l}$ and $\varphi_{2l+1}$, $l \in \mathbb{N}$, and then between all the pairs which can be obtained from the above ones by the reflection with respect to the real line and by the $2\pi$-translations. Denote the complex plane thus cut by $\mathbb{C}_{\kappa_0}$.

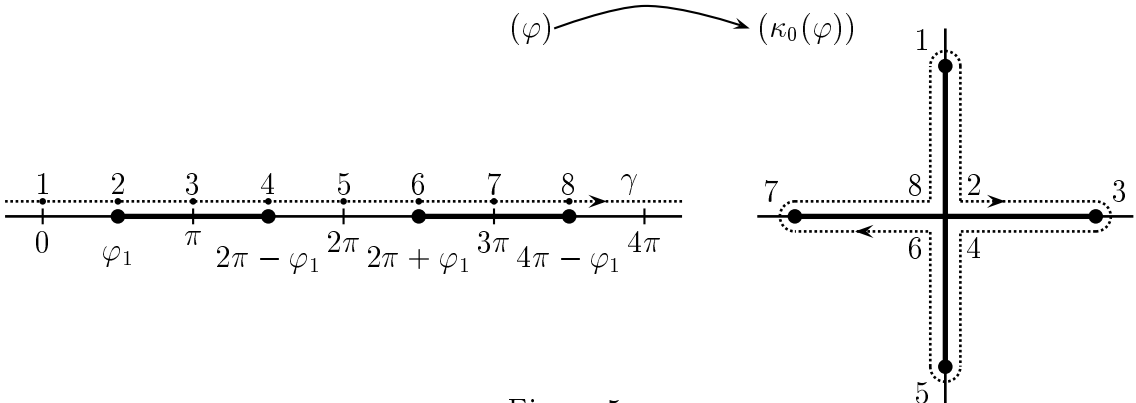

Figure 5.

The branch $\kappa_0$ can be continued up a single valued function defined on $\mathbb{C}_{\kappa_0}$. This functions has the following properties:

(4.3) $$\kappa_0(-\varphi) = \kappa_0(\varphi),$$

(4.4) $$\kappa_0(\overline{\varphi}) = -\overline{\kappa_0(\varphi)},$$

(4.5) $$\kappa_0(\varphi + 2\pi) = -\kappa_0(\varphi)$$

for $\varphi \in \mathbb{C}_{\kappa_0}$.

Relation (4.3) follows from the definition of $\kappa_0$, formula (4.4) is a consequence of the observation that $\kappa_0(\varphi) \in i\mathbb{R}_+$ if $0 < \varphi < \varphi_1$. The last relation comes from the analysis of Fig. 5.

We finish this subsection with the

**Lemma 4.1.** *If*

(4.6) $$\operatorname{Re}\varphi \neq \pi + 2\pi m, \quad m \in \mathbb{Z},$$

*then, we have*

(4.7) $$\kappa_0(\varphi) = \frac{i}{\sqrt{2}} e^{\mp i\varphi/2} + O\left(e^{-|\varphi|/2}\right), \quad \operatorname{Im}\varphi \to \pm\infty.$$

*This representation is uniform in $\operatorname{Re}\varphi$ outside any fixed $\delta$-vicinity of the lines $\operatorname{Re}\varphi = \pi + 2\pi m$.*

*Proof.* The asymptotics (2.12) and the relation between $k_0(E)$ and $\kappa_0(\varphi)$ for $\varphi$ being in the shadowed domain in Fig. 3 imply the asymptotics for $\operatorname{Im}\varphi \to +\infty$ in the case where $0 \leq \operatorname{Re}\varphi < \pi$. Relations (4.3)-(4.5) show that this asymptotics remains valid for all $\operatorname{Re}\varphi$ satisfying the condition (4.6). $\square$

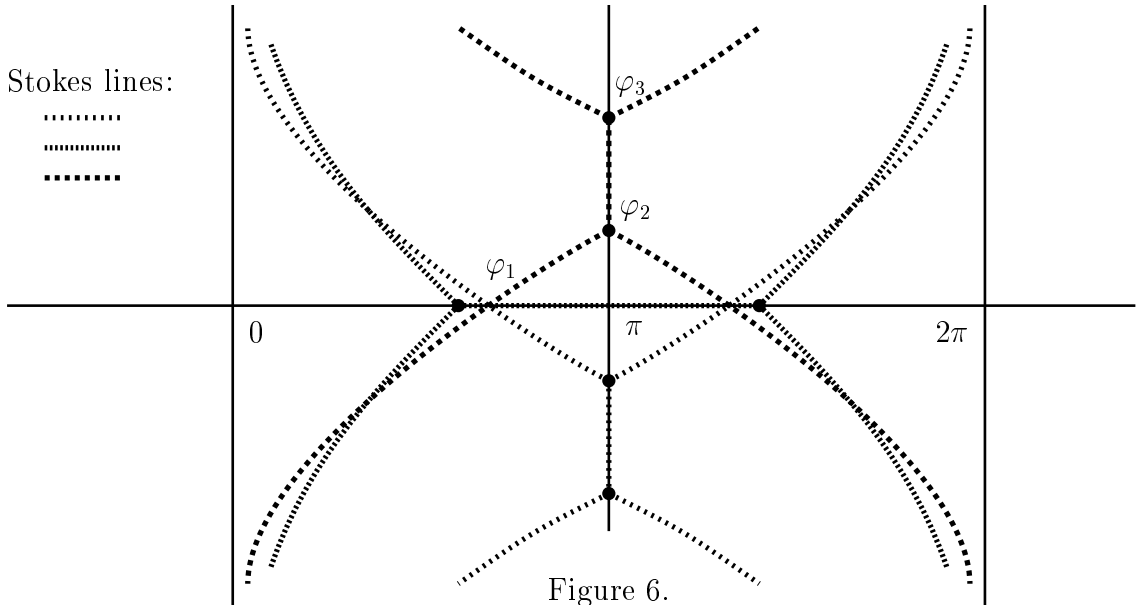

Figure 6.

**4.2. The Stokes lines.** The definition of the Stokes lines is fairly standard. The integral $\int \kappa \, d\varphi$ has the same branching points as the complex momentum. Let $\varphi_0$ be one of them. Consider the curves beginning at $\varphi_0$ and described by the

equation

$$\text{(4.8)} \qquad \operatorname{Im} \int_{\varphi_0}^{\varphi} \left( \kappa\left(\xi\right) - \kappa\left(\varphi_0\right) \right) d\xi = 0.$$

These curves are the *Stokes lines* beginning at $\varphi_0$. It follows from (4.2) that the Stokes line definition is independent of the choice of the branch of $\kappa$ in (4.8).

Assume, as before, that $E$ satisfies (1.26). In this case, at any branching point begin exactly three Stokes lines. The angle between any two of them at this point are equal to $\frac{2\pi}{3}$. Some of the Stokes lines are finite: they connect pairs of finite branching points; some of the Stokes lines are infinite: they go from finite branching points to the infinity. The asymptotes of the infinite lines can be described by means of Lemma 4.1.

The character of the Stokes lines in the domain $0 < \operatorname{Re} z < 2\pi$ is shown in Fig. 6.

Consider the branching point $\varphi_1$. Two of the Stokes lines beginning at $\varphi_1$ go to $\pm i\infty$ and have the asymptote $\operatorname{Re}\varphi = 0$. The last Stokes line starting at $\varphi_1$ is finite; it goes along the real axis to $2\pi - \varphi_1$. In the same way, for each $l > 1$, two of the Stokes lines beginning at the point $\varphi_l$ are infinite. And the last one is finite and runs along the line $\operatorname{Re}\varphi = \pi$ from one branching point to another.

The whole picture of the Stokes lines in the domain $0 \leq \operatorname{Re}\varphi \leq 2\pi$ is symmetric with respect to the real line as well as with respect to the line $\operatorname{Re}\varphi = \pi$.

If $\varphi$ is a branching point for $\kappa$, then so is the point $\varphi + 2\pi$. The Stokes lines passing through $\varphi + 2\pi$ are just the $2\pi$-translates of the Stokes lines beginning $\varphi$.

**4.3. Canonical domains.** Here we describe some canonical domains in the case when $W(\varphi) = \cos(\varphi)$.

4.3.1. In Fig. 7, we show two canonical domains. The first one (part A) corresponds to the case where $\zeta_1 = +i\infty$ and $\zeta_2 = -i\infty$. Part B corresponds to $\zeta_1$ and $\zeta_2$ finite.

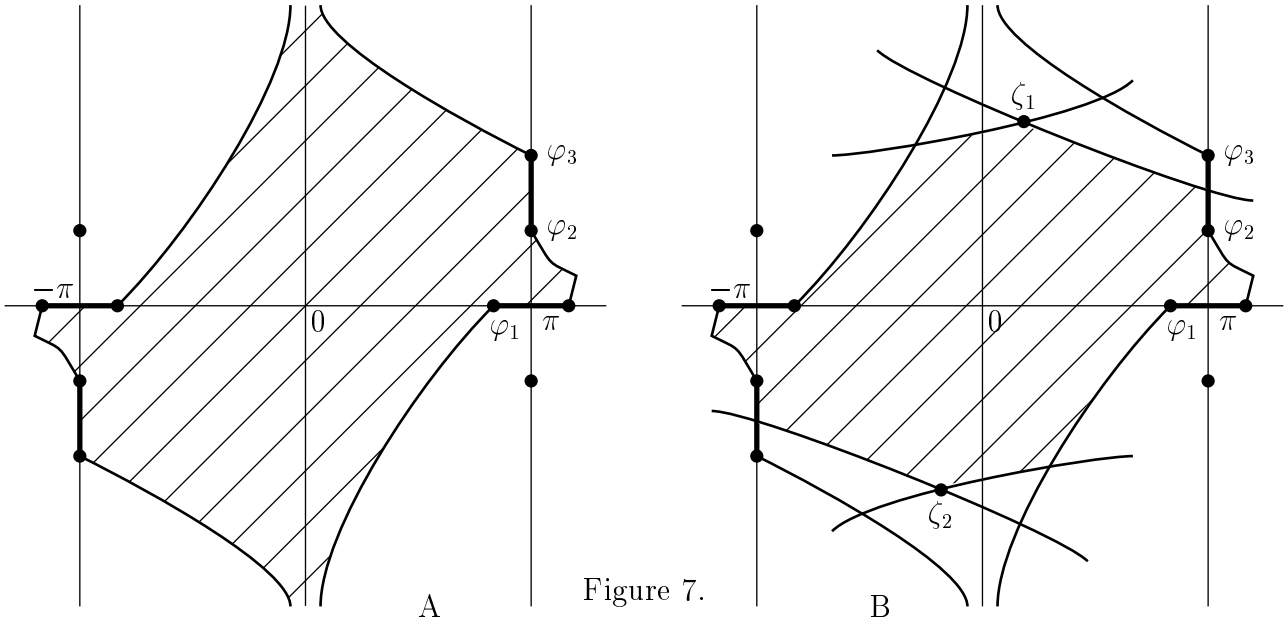

Figure 7.

This domain $K_0$ is canonical with respect to the branch $\kappa_0$. The boundary of $K_0$ consists of segments of Stokes lines. This domain exists if $E - 1 < E_1$.

To check that $K_0$ is canonical, one has to investigate the families of the lines Re $\int^\varphi \kappa_0 \, d\varphi = Const$ and Re $\int^\varphi (\kappa_0 - \pi) \, d\varphi = Const$ inside this domain.

In Fig. 8, we show these families of lines. In the left hand side part of Fig. 8, we show some of the lines Re $\int^\varphi \kappa_0 \, d\varphi = Const$, and, in the right hand side, we show some of the lines Re $\int^\varphi (\kappa_0 - \pi) \, d\varphi = Const$. The arrows indicate the directions of increase of the functions Im $\int^\varphi \kappa_0 \, d\varphi$ and Im $\int^\varphi (\kappa_0 - \pi) \, d\varphi$.

In Fig. 9, we show one more canonical domain $K_1$. This domain exists if $E_1 < E + 1 < E_2$. It is canonical with respect to the branch fixed by the condition $0 < \kappa_1(\pi) < \pi$. We note that $\kappa_1$ has the following symmetries: for $\varphi \in K_1$,

(4.9) $$\kappa_1(2\pi - \varphi) = \kappa_1(\varphi),$$

(4.10) $$\kappa_1(\overline{\varphi}) = \overline{\kappa_1(\varphi)}.$$

The function $\kappa_1$ can be continued on the complex plane cut along intervals $[2k\pi - \varphi_1, 2k\pi + \varphi_1]$ along the real line and as in Fig. 2 outside of the real line. In this case we have

(4.11) $$\kappa_1(\varphi + 2\pi) = -\kappa_1(\varphi).$$

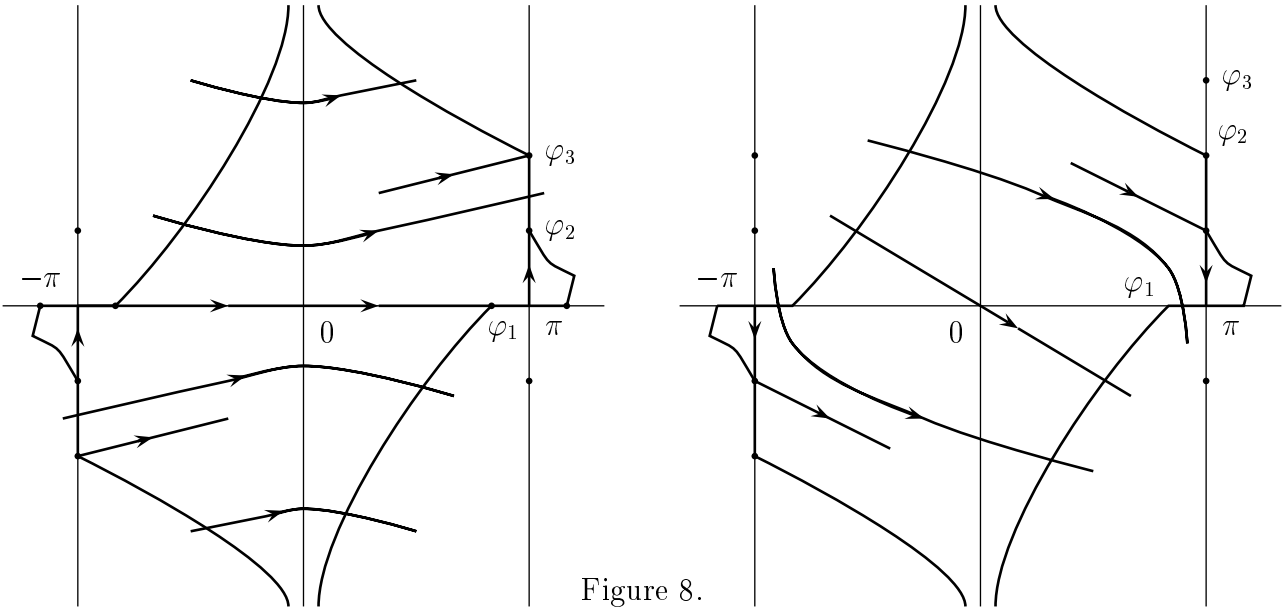

Figure 8.

Relations (4.3) – (4.5) and (4.9) – (4.11) show that any domain obtained of $K_0$ or $K_1$ by means of the $2\pi$-translations and the reflection with respect to the real line are canonical.

4.4. **Admissible sub-domains.** We first note that the admissible sub-domains corresponding to $K_0$ and $K_1$ are bounded. Moreover we do not need the notion of strictly canonical domains anymore. Indeed, one can show that, if $K$ is a canonical domain and $A$ is one of its admissible sub-domains, then there exists $B$ a strictly canonical domain such that $A \subset B \subset K$. Thus, the asymptotics given by Theorem 3.1 are valid on any fixed admissible sub-domain of a given canonical domain.

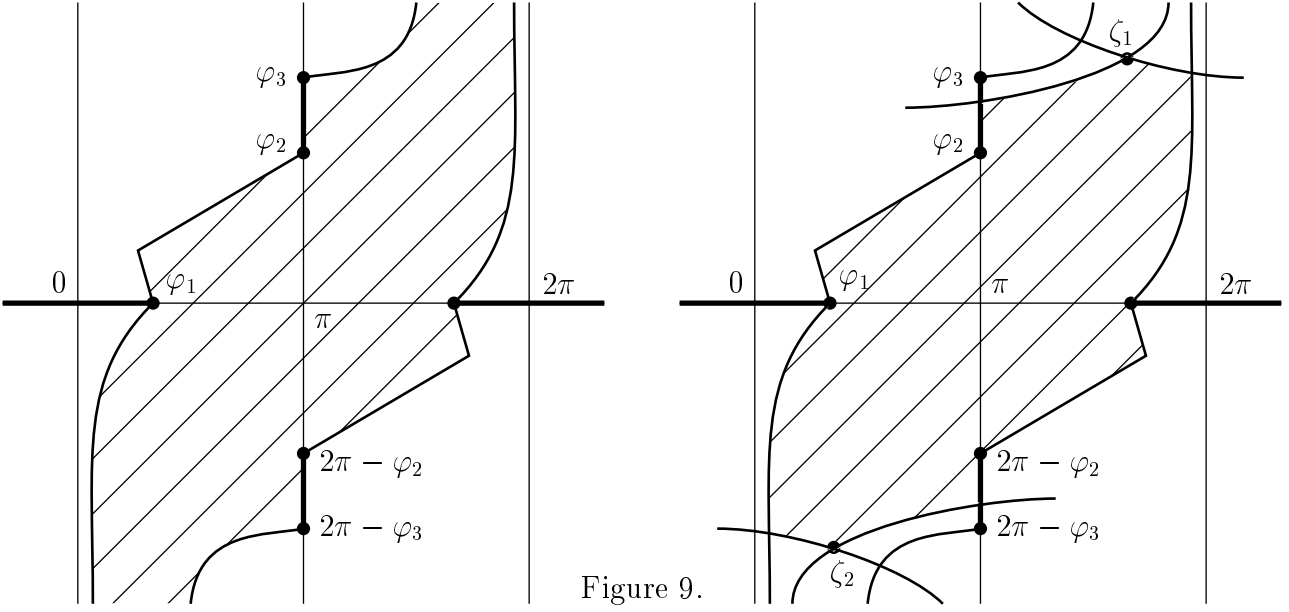

Figure 9.

## 5. The computation of the monodromy matrix for equation (1.24)

This section is devoted to the computation of the monodromy matrix associated to the consistent basis of solutions of (1.24) constructed in Theorem 3.1 for the canonical domain $K_0$ defined in section 4. Before we begin, let us underline that we will commit a slight abuse of language and write the monodromy matrix $M(\varphi)$ instead of $M(\phi) = M(\varphi/\varepsilon)$.

We now state the main theorem and then use the rest of the section to prove it. To describe the monodromy matrix, we define the following action integrals

$$t = exp\left(2\frac{i}{\varepsilon}\int_0^{\varphi_1} \kappa\right), \quad t_1 = exp\left(2\frac{i}{\varepsilon}\int_{\varphi_2}^{\pi}(\kappa - \pi)\right),$$
(5.1)
$$\phi_1 = \frac{2}{\varepsilon}\int_{\varphi_1}^{\pi} \kappa \quad \text{and} \quad \phi_2 = \frac{i\varepsilon}{2\pi}\int_{\varphi_1}^{\varphi_2}(\omega_- - \omega_+).$$

Here the functions $\kappa, \omega_\pm$ denote the branches of these functions first defined in $K_0$ and then continued along the path $\gamma$ (shown on Fig. 10) to the domain $K_0 \cup K_1$. We prove the

**Theorem 5.1.** *Define coefficients of the monodromy matrix $M(\varphi)$ by*

$$M(\varphi) = \begin{pmatrix} m_{11}(\varphi) & m_{12}(\varphi) \\ m_{21}(\varphi) & m_{22}(\varphi) \end{pmatrix}.$$

*Then, for any fixed $\delta > 0$, for $\varepsilon > 0$ small enough, uniformly for $|\operatorname{Im}\varphi| \leq \operatorname{Im}\varphi_3 - \delta$, we have*

$$m_{11}(\varphi) = t \cdot e^{i\phi_1}(1 + o(1)),$$

$$m_{12}(\varphi) = ie^{\int_0^{\varphi_1}\omega_- - \int_0^{\varphi_1}\omega_+} \cdot$$

$$\left(e^{i\phi_1}(1 + o(1)) - t_1 \cdot e^{\frac{2i\pi}{\varepsilon}(\varphi + \pi + \phi_2)}(1 + o(1))\right)$$

$$m_{21}(\varphi) = -ie^{-\int_0^{\varphi_1}\omega_- + \int_0^{\varphi_1}\omega_+}.$$

$$\left(e^{i\phi_1}(1+o(1)) - t_1 \cdot e^{-\frac{2i\pi}{\varepsilon}(\varphi+\pi+\phi_2)}(1+o(1))\right)$$

$$m_{22}(\varphi) = \frac{1}{t}e^{-i\phi_1}(1+o(1)) + \frac{1}{t}e^{i\phi_1}(1+o(1))$$
$$- \frac{t_1}{t}\left(e^{i\frac{2\pi}{\varepsilon}(\varphi+\pi+\phi_2)}(1+o(1)) + e^{-i\frac{2\pi}{\varepsilon}(\varphi+\pi+\phi_2)}(1+o(1))\right).$$

Theorem 1.2 is an immediate corollary of Theorem 5.1 if one sets

$$T = \begin{pmatrix} exp\left(\int_0^{\varphi_1}\omega_+\right) & 0 \\ 0 & exp\left(-\int_0^{\varphi_1}\omega_-\right) \end{pmatrix} \text{ and } C = \frac{2\pi}{\varepsilon}(\pi+\phi_2),$$

and if we prove that $\phi_2$ is real. This is done at the end of the proof of Theorem 5.1.

For the canonical domain $K_0$ (defined in the previous section), denote by $(f_{\pm}^0)$ the solutions of equation (1.24) constructed in Theorem 3.1. The monodromy matrix associated to this basis is then defined by

(5.2) $$F^0(u, \varphi + 2\pi) = M(\varphi)F^0(u, \varphi)$$

where $F^0 = (f_+^0, f_-^0)^T$.

5.1. **The basic ideas guiding the proof.** In a subset of $K_0$ where we know both the asymptotic behavior of $F^0(u, \varphi + 2\pi)$ and $F^0(u, \varphi)$, it will be easy to compute the asymptotics of $M$. Unfortunately a quick look at section 4.3 shows us that $K_0 \cap (K_0 + 2\pi) = \emptyset$. Therefore we need to introduce an auxiliary canonical domain.

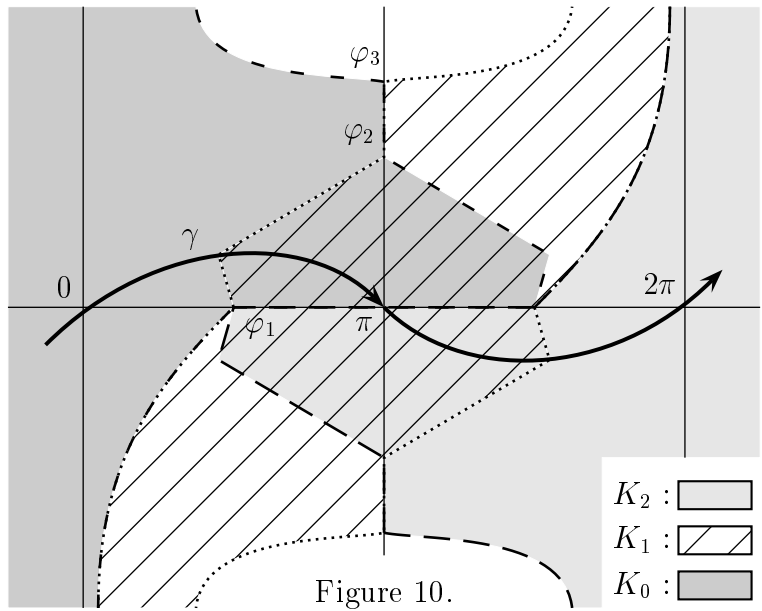

Figure 10.

5.1.1. *Auxiliary canonical domains and transition matrices.* The first auxiliary domain we need we have already defined: it is $K_1$ (see Fig. 9 and subsection 4.3).

Let us denote by $(f_\pm^1)$ the consistent basis constructed by Theorem 3.1 in $K_1$. As $(f_\pm^0)$ and $(f_\pm^1)$ are both basis of solutions of (1.19), we define

(5.3) $$F^0(u, \varphi) = T_1(\varphi) F^1(u, \varphi).$$

We will be able to compute the asymptotic behavior of $T_1$ in $K_0 \cap K_1$. We will see that we will actually be able to go beyond this simple remark using analyticity, periodicity and a continuation lemma.

The second auxiliary canonical domain needed will be denoted by $K_2$. It is defined to be $K_2 = K_0 + 2\pi$. This domain is canonical for the branch of $\kappa$ obtained by continuing $\kappa_0$ along the curve $\gamma$ shown on Fig. 10. Moreover it is immediate that if we define $f_\pm^2(\varphi) = f_\pm^0(\varphi - 2\pi)$, then for $K_2$, this defines a consistent basis. As $(f_\pm^2)$ and $(f_\pm^1)$ are both basis of solutions of (1.19), we define

(5.4) $$F^1(u, \varphi) = T_2(\varphi) F^2(u, \varphi).$$

As above, we will be able to compute the asymptotic behavior of the matrix $T_2$ in $K_1 \cap K_2$, and beyond that region.

Putting (5.3) and (5.4) together, we get that

$$F^0(u, \varphi) = T_1(\varphi) T_2(\varphi) F^2(u, \varphi) = T_1(\varphi) T_2(\varphi) F^0(u, \varphi - 2\pi).$$

Hence the monodromy matrix reads

(5.5) $$M(\varphi) = T_1(\varphi + 2\pi) T_2(\varphi + 2\pi).$$

**5.2. The asymptotics of $T_1$ and $T_2$.** These asymptotics are given by the

**Proposition 5.1.** *Let us define the coefficients of $T_1$ by*

$$T_1(\varphi) = \begin{pmatrix} a_1(\varphi) & b_1(\varphi) \\ c_1(\varphi) & d_1(\varphi) \end{pmatrix}.$$

*Then, for any fixed $0 < Y < \operatorname{Im} \varphi_3$, for $\varepsilon > 0$ small enough, uniformly for $|\operatorname{Im} \varphi| \le Y$, we have*

(5.6) $a_1(\varphi) = e^{\frac{i}{\varepsilon} \int_0^\pi \kappa_0 d\varphi + \int_0^\pi \omega_+^0 d\varphi}(1 + o(1))$,

(5.7) $b_1(\varphi) = e^{\frac{i}{\varepsilon} \int_0^{\varphi_1} \kappa_0 d\varphi} \cdot o(1)$,

(5.8) $$\begin{aligned} c_1(\varphi) = -i e^{-\frac{i}{\varepsilon} \int_0^\pi \kappa_0 d\varphi} &\left( e^{2\frac{i}{\varepsilon} \int_{\varphi_1}^\pi \kappa_0 d\varphi} e^{\int_0^{\varphi_1} \omega_-^0 d\varphi - \int_\pi^{\varphi_1} \omega_+^0 d\varphi}(1 + o(1)) \right. \\ &\left. - e^{-2\frac{i}{\varepsilon} \int_\pi^{\varphi_2}(\kappa_0 - \pi) d\varphi} e^{\int_0^{\varphi_2} \omega_-^0 d\varphi - \int_\pi^{\varphi_2} \omega_+^0 d\varphi} e^{-2\pi \frac{i}{\varepsilon}(\varphi - \pi)}(1 + o(1)) \right), \end{aligned}$$

(5.9) $d_1(\varphi) = e^{-\frac{i}{\varepsilon} \int_0^\pi \kappa_0 d\varphi + \int_0^\pi \omega_-^0 d\varphi}(1 + o(1))$.

*Here all the integrals are taken along contours living in $K_0 \cup K_1 \cup K_2$.*

and by the

**Proposition 5.2.** *Let us define the coefficients of $T_2$ by*

$$T_2(\varphi) = \begin{pmatrix} a_2(\varphi) & b_2(\varphi) \\ c_2(\varphi) & d_2(\varphi) \end{pmatrix}.$$

Then, for any $0 < Y < \operatorname{Im} \varphi_3$, for $\varepsilon > 0$ small enough, uniformly for $|\operatorname{Im} \varphi| \leq Y$, we have

(5.10) $\quad a_2(\varphi) = e^{\frac{i}{\varepsilon} \int_\pi^{2\pi} \kappa_1 d\varphi + \int_\pi^{2\pi} \omega_+^1 d\varphi}(1 + o(1)),$

(5.11) $\quad b_2(\varphi) = ie^{-\frac{i}{\varepsilon} \int_\pi^{2\pi} \kappa_1 d\varphi} \left( e^{2\frac{i}{\varepsilon} \int_\pi^{2\pi-\varphi_1} \kappa_1 d\varphi} e^{\int_\pi^{2\pi-\varphi_1} \omega_+^1 d\varphi - \int_{2\pi}^{2\pi-\varphi_1} \omega_-^1 d\varphi}(1 + o(1)) \right.$
$\quad \left. - e^{2\frac{i}{\varepsilon} \int_\pi^{2\pi-\varphi_2}(\kappa_1 - \pi) d\varphi} e^{\int_{2\pi-\varphi_2}^{2\pi} \omega_-^1 d\varphi - \int_{2\pi-\varphi_2}^{\pi} \omega_+^1 d\varphi} e^{2\pi \frac{i}{\varepsilon}(\varphi - \pi)}(1 + o(1)) \right),$

(5.12) $\quad c_2(\varphi) = e^{-\frac{i}{\varepsilon} \int_{2\pi}^{2\pi-\varphi_1} \kappa_1 d\varphi} \cdot o(1),$

(5.13) $\quad d_2(\varphi) = e^{-\frac{i}{\varepsilon} \int_\pi^{2\pi} \kappa_1 d\varphi + \int_\pi^{2\pi} \omega_-^1 d\varphi}(1 + o(1)).$

Here all the integrals are taken along contours living in $K_0 \cup K_1 \cup K_2$.

The next subsections will be devoted to the proof of Propositions 5.1 and 5.2. The proofs being very similar we will not dwell on the one of Proposition 5.2, but just give the necessary modifications of the proof of Proposition 5.1

### 5.3. The proof of Proposition 5.1.
By definition, we have

(5.14) $\quad a_1(\varphi) = \frac{1}{w_1} w(f_+^0, f_-^1), \quad b_1(\varphi) = -\frac{1}{w_1} w(f_+^0, f_+^1),$

(5.15) $\quad c_1(\varphi) = \frac{1}{w_1} w(f_-^0, f_-^1), \quad d_1(\varphi) = -\frac{1}{w_1} w(f_-^0, f_+^1),$

(5.16) $\quad w_1 = w(f_+^1, f_-^1) = f_-^1 \frac{d}{du} f_+^1 - f_+^1 \frac{d}{du} f_-^1.$

We recall here that $(f_\pm^{0,1})$ are the solution of (1.19) constructed in Theorem 3.1 for the canonical domains $K_{0,1}$. In an admissible sub-domain of $K_{0,1}$, by definition, $f_\pm^{0,1}$ have the following asymptotics when $\varepsilon \to 0$

(5.17) $\quad f_\pm^{0,1}(u, \varphi) = e^{\pm \frac{i}{\varepsilon} \int_{\varphi_{0,1}}^\varphi \kappa_{0,1} d\varphi}(\Psi_\pm^{0,1}(u, \varphi) + o(1))$

where $\varphi_0 = 0$, $\varphi_1 = \pi$,

$$\Psi_\pm^{0,1}(x, \varphi) = \sqrt{k_E'} e^{\int_{\varphi_{0,1}}^\varphi \omega_\pm^{0,1} d\varphi} \psi_\pm^{0,1}(x, \varphi)$$

and $\psi_\pm^1$ (resp. $\omega_\pm^1$) are the analytic continuations of $\psi_\pm^0$ (resp. $\omega_\pm^0$) along $\gamma$ from $K_0$ to $K_1$ (see Fig. 10 and 12).

This implies that

(5.18) $\quad w_1 = w(f_+^1, f_-^1) + o(1) = k'(E + 1) \cdot w_0(\pi) + o(1)$

where $w_0$ is the analytic continuation of the Wronskian $w(f_+^0, f_-^0)$ from $K_0$ to $K_1$ along $\gamma$ (see Fig. 10).

#### 5.3.1. *Periodicity.*
One remarkable fact about $(f_\pm^{0,1})$ is that

$$f_\pm^{0,1}(u, \varphi + \varepsilon) = f_\pm^{0,1}(u + 1, \varphi).$$

Using the fact that the Wronskians taken in (5.14) are independent of the point $u$, we get that these Wronskians are $\varepsilon$-periodic. So that the coefficients of $T_1$ are $\varepsilon$-periodic. Hence to get a total control of a coefficients of $T_1$, say $a_1$, in a horizontal strip, we only need to control it in some vertical sub-strip on width $\varepsilon$.

5.3.2. *Analyticity and Fourier series.* Moreover $a_1(\varphi)$ (and the other coefficients as well) is obviously analytic in $\varphi$ by Theorem 3.1. This allows us to expand $a_1$ into an exponentially converging Fourier series. To control this series, we need to control both positive and negative Fourier coefficients, hence we need to have a control on the Wronskians in (5.14) et al. beyond the obvious zone $K_0 \cap K_1$ (see Fig. 10). Actually we will get a uniform control in a strip of the form

$$-\Delta \leq \operatorname{Im}\varphi \leq \operatorname{Im}\varphi_3 - \delta, \quad \Delta > \delta > 0.$$

Therefore we divide this strip into three different smaller strips called (I), (II) and (III) (see Fig. 11). Each of these strip will require a different type of computation.

We will now start with describing the behavior of the main objects of our construction in each of the different strips.

In the strip (II), we have a non-empty intersection for $K_0$ and $K_1$. In the intersection, we have

(5.19) $$\kappa_1 = \kappa_0, \quad \Psi_\pm^1 = \Psi_\pm^0 \quad \text{and} \quad \omega_\pm^1 = \omega_\pm^0.$$

Consider the strip (I). Denote by $\gamma$ the common boundary of $K_0$ and $K_1$ in this strip (see Fig. 13). It is one of the Stokes lines beginning at $\varphi = \varphi_1$. On $\gamma$, we get

(5.20) $$\kappa_1(\varphi + 0) = -\kappa_0(\varphi - 0)$$
(5.21) $$\psi_\pm^1(u, \varphi + 0) = \psi_\mp^0(u, \varphi - 0)$$
(5.22) $$\omega_\pm^1(\varphi + 0) = \omega_\mp^0(\varphi - 0)$$

where $\varphi + 0$ (resp. $\varphi - 0$) denotes the limit taken from the right (resp. left). Formula (5.20) follows from (2.10), (5.21) from (2.6) and (5.22) from the definition of $\omega_\pm$.

We also notice that the conformal properties of $\kappa_0$ (see subsection 4.1) and (5.20) show that, for $\varphi \in$ (I) and $\operatorname{Im}\varphi < 0$, we have

$$-\operatorname{Im}\kappa_0(\varphi) = \operatorname{Im}\kappa_1(\varphi) < 0.$$

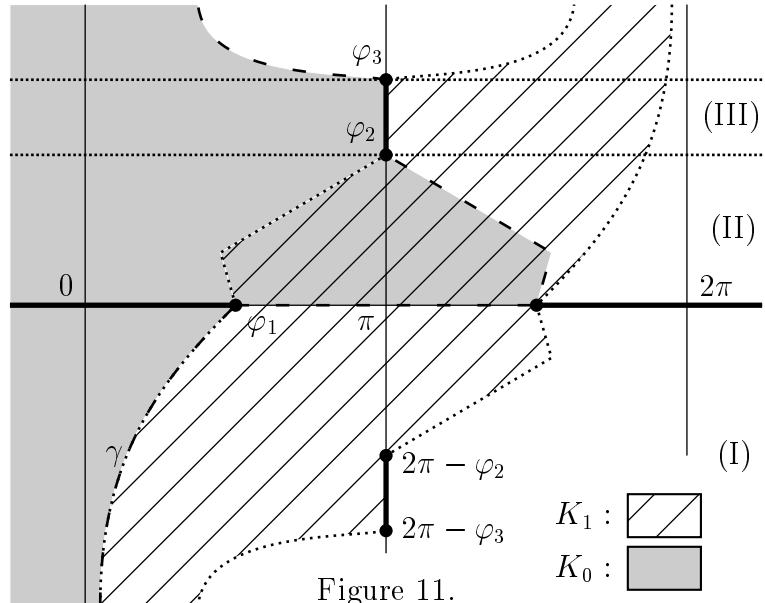

Figure 11.

In the strip (III), the common boundary of $K_0$ and $K_1$ consists in the interval $[\varphi_2, \varphi_3]$. This interval is the Stokes line joining $\varphi_2$ to $\varphi_3$. Along $[\varphi_2, \varphi_3]$, we have

(5.23) $$\kappa_1(\varphi + 0) = 2\pi - \kappa_0(\varphi - 0)$$

(5.24) $$\psi_\pm^1(u, \varphi + 0) = \psi_\mp^0(u, \varphi - 0)$$

(5.25) $$\omega_\pm^1(\varphi + 0) = \omega_\mp^0(\varphi - 0).$$

As in region (I), for $\operatorname{Im} \varphi_2 < \operatorname{Im} \varphi < \operatorname{Im} \varphi_3$, we have

(5.26) $$-\operatorname{Im} \kappa_0(\varphi) = \operatorname{Im} \kappa_1(\varphi) < 0.$$

5.3.3. *A continuation lemma.* To make the strategy outlined above work, we need to solve a technical problem that is that the asymptotics given by Theorem 3.1 are only valid in admissible sub-domains of the canonical domains. So, if we want to continue the asymptotics across the boundary of a canonical domain to join, say, an admissible sub-domain of $K_0$ to an admissible sub-domain of $K_1$, we will have to bridge a gap of width $2\delta$. Therefore we will use the continuation lemma presented in the next subsection. We prove the

**Lemma 5.1.** *Let $\varphi_-, \varphi_+$ be two fixed points such that*
- $\operatorname{Im} \varphi_- = \operatorname{Im} \varphi_+$,
- *there is no branching point of $\varphi \mapsto \kappa(\varphi)$ on the interval $[\varphi_-, \varphi_+]$.*

*Fix a continuous branch of $\kappa$ on $[\varphi_-, \varphi_+]$. Pick $\varphi_0 \in (\varphi_-, \varphi_+)$. Let $f(u, \varphi)$, $f_\pm(u, \varphi)$ be normal (i.e $\varepsilon$-periodic in $\varphi$) solutions of (1.19) for $\varphi \in [\varphi_-, \varphi_+]$ and $u \in [-U, U]$ such that:*

1. $f(u, \varphi) = e^{\frac{i}{\varepsilon} \int_{\varphi_0}^{\varphi} \kappa d\varphi}(\Psi_+(u, \varphi) + o(1))$ *for $\varphi \in [\varphi_-, \varphi_0]$ when $\varepsilon \to 0$, and the asymptotic is differentiable in $u$.*
2. $f_\pm(u, \varphi) = e^{\pm\frac{i}{\varepsilon} \int_{\varphi_0}^{\varphi} \kappa d\varphi}(\Psi_\pm(u, \varphi) + o(1))$ *for $\varphi \in [\varphi_-, \varphi_+]$ when $\varepsilon \to 0$, and the asymptotic is differentiable in $u$. Here $(\Psi_\pm)$ are chosen as in Theorem 3.1.*

*Then,*
- *if $\operatorname{Im}(\kappa(\varphi)) > 0$ for all $\varphi \in [\varphi_-, \varphi_+]$, there exists $C > 0$ such that, for $\varepsilon > 0$ small enough and $\varphi \in [\varphi_-, \varphi_+]$, we have*

$$\left| \frac{df}{du}(u, \varphi) \right| + |f(u, \varphi)| \leq C e^{\frac{1}{\varepsilon} \int_{\varphi_0}^{\varphi} |\operatorname{Im}\kappa| d\varphi}.$$

- *if $\operatorname{Im}(\kappa(\varphi)) < 0$ for all $\varphi \in [\varphi_-, \varphi_+]$, then, for $\varphi \in [\varphi_-, \varphi_+]$, we have*

$$f(u, \varphi) = e^{\frac{i}{\varepsilon} \int_{\varphi_0}^{\varphi} \kappa d\varphi}(\Psi_+(u, \varphi) + o(1)),$$

*and the asymptotic is differentiable in $u$.*

**Remark 5.1.** This lemma is a precise version of a well known WKB heuristic that says that, if a function $f$ admits on $[\varphi_-, \varphi_0]$ an asymptotic representation that is defined and exponentially increasing on $[\varphi_-, \varphi_+]$, then the $f$ admits the same representation on $[\varphi_0, \varphi_+]$.

As throughout this paper, the estimates and asymptotics used for the solutions of equation (1.19) are differentiable in $u$, we will only write the relevant estimates and asymptotics for the solution. We will not bother repeating them for $\dfrac{df}{du}$.

As a last remark, let us say that in the rest of this section, to simplify the notations, we will denote by $C$ every positive constant that is independent of $\varphi$ and $\varepsilon$.

**Proof.** By assumption 2 of Lemma 5.1, $(f_\pm)$ are linearly independent solutions of (1.19). So that we can write

$$f = a(\varphi) f_+ + b(\varphi) f_- \tag{5.27}$$

where

$$a(\varphi) = \frac{w(f, f_-)}{w(f_+, f_-)} \quad \text{and} \quad b(\varphi) = \frac{w(f_+, f)}{w(f_+, f_-)}.$$

As in the proof of Theorem 3.1, we have

$$|w(f_+, f_-)| = |k'_E w_{0|\varphi = \varphi_0}(1 + o(1))| = C \neq 0.$$

Hence, on $[\varphi_-, \varphi_0]$, by assumptions 1. and 2., we have

$$a(\varphi) = 1 + o(1). \tag{5.28}$$

$a$ is $\varepsilon$-periodic in $\varphi$ so that (5.28) holds on $[\varphi_-, \varphi_+]$.

Let us now estimate $b$. We start with the case when $\text{Im}(\kappa(\varphi)) > 0$ for all $\varphi \in [\varphi_-, \varphi_+]$. Then, on $[\varphi_-, \varphi_0]$, we have

$$|w(f_+, f)| \leq C \left| e^{\frac{2i}{\varepsilon} \int_{\varphi_0}^{\varphi} \kappa d\varphi} \right| \leq C e^{\frac{2}{\varepsilon} \int_{\varphi}^{\varphi_0} \text{Im}(\kappa) d\varphi}.$$

So, for $\varphi \in [\varphi_0 - \varepsilon, \varphi_0]$, we get that $|b(\varphi)| \leq C$; hence, by $\varepsilon$-periodicity, we have $|b(\varphi)| \leq C$ for $\varphi \in [\varphi_-, \varphi_+]$. On the other hand, on $[\varphi_0, \varphi_+]$, we have

$$|f_+| \leq C e^{-\frac{1}{\varepsilon} \int_{\varphi_0}^{\varphi} \text{Im}(\kappa) d\varphi}, |f_-| \leq C e^{\frac{1}{\varepsilon} \int_{\varphi_0}^{\varphi} \text{Im}(\kappa) d\varphi}.$$

So that, by (5.27), we get the announced result on $[\varphi_0, \varphi_+]$.

If $\text{Im}(\kappa(\varphi)) < 0$ for all $\varphi \in [\varphi_-, \varphi_+]$ then, for $\varphi \in [\varphi_-, \varphi_0]$, we estimate

$$|b(\varphi)| \leq C \left| e^{\frac{2i}{\varepsilon} \int_{\varphi_0}^{\varphi} \kappa d\varphi} \right|$$

Hence using this estimate for $\varphi \in [\varphi_- - \varepsilon, \varphi_-]$ and the $\varepsilon$-periodicity, for $\varphi \in [\varphi_-, \varphi_+]$, we get

$$|b(\varphi)| \leq C e^{-\frac{2}{\varepsilon} \int_{\varphi_-}^{\varphi_0} |\text{Im}(\kappa)| d\varphi}.$$

Hence, by (5.27), on $[\varphi_0, \varphi_+]$, we have

$$\begin{aligned} f &= a(\varphi) f_+ + b(\varphi) f_- \\ &= e^{\frac{i}{\varepsilon} \int_{\varphi_0}^{\varphi} \kappa d\varphi} \left( \Psi_+(u, \varphi) + O\left( e^{-\frac{2i}{\varepsilon} \int_{\varphi_0}^{\varphi} \kappa d\varphi} e^{-\frac{2}{\varepsilon} \int_{\varphi_-}^{\varphi_0} |\text{Im}(\kappa)| d\varphi} \right) + o(1) \right) \\ &= e^{\frac{i}{\varepsilon} \int_{\varphi_0}^{\varphi} \kappa d\varphi} \left( \Psi_+(u, \varphi) + o(1) \right). \end{aligned}$$

This ends the proof of Lemma 5.1. $\square$

5.3.4. *One more auxiliary canonical domain.* To use Lemma 5.1, we introduce one more auxiliary canonical domain. We denote it by $\tilde{K}_0$. Let $\tilde{\kappa}_0$ be the continuation of $\kappa_0$ through the cut $(\varphi_2, \varphi_3)$. Then $\tilde{K}_0$ is defined to be the canonical domain for $\tilde{\kappa}_0$ comprising the vertical interval $(\varphi_2, \varphi_3)$. The part of $\tilde{K}_0$ that will be needed to compute $T_1$ is represented on Fig. 12. It will enable us to continue the asymptotics of a solution corresponding to the canonical domains $K_0$ to the domain $K_1$ along the horizontal lines $\operatorname{Im} \varphi = \eta$ when $\operatorname{Im} \varphi_2 < \eta < \operatorname{Im} \varphi_3$.

We are now ready to start with the computations of the coefficients of $T_1$.

5.3.5. *The asymptotics of $d_1(\varphi)$.* We compute the asymptotics of $w(f_-^0, f_+^1)$ when $\varepsilon \to 0$.

Region (I). In view of (5.26) by Lemma 5.1, we can continue continue the asymptotics given by Theorem 3.1 for $f_-^0$ from $[\varphi_-, \varphi_0]$ to the whole interval $[\varphi_-, \varphi_+]$ (see Fig. 13). The consistent basis needed to apply Lemma 5.1 is obtained by noticing that, as $\kappa_0(\overline{\varphi}) = -\overline{\kappa_0(\varphi)}$, the domain $\overline{K_0}$ symmetric of $K_0$ with respect to the real axis also is a canonical domain. This domain comprises a neighborhood of $\gamma$, the common boundary of $K_0$ and $K_1$ in region (I) (see Fig. 13).

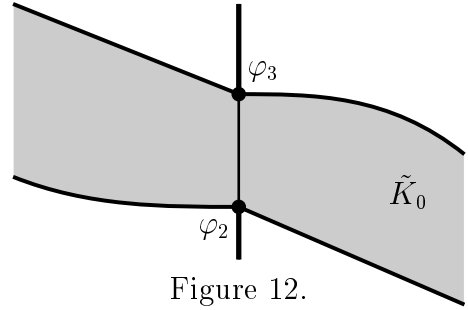

Figure 12.

By what has been said above, on $[\varphi_b, \varphi_+]$, we have

$$f_-^0 = e^{-\frac{i}{\varepsilon}\int_0^{\varphi_1} \kappa_0 d\varphi - \frac{i}{\varepsilon}\int_{\varphi_1}^{\varphi_b} \kappa_0 d\varphi} e^{-\frac{i}{\varepsilon}\int_{\varphi_b}^{\varphi} \kappa_0 d\varphi}(\Psi_-^0 + o(1)).$$

By (5.17) and as $\kappa_1 = -\kappa_0$ in the strip (I), we know that

$$f_+^1 = e^{\frac{i}{\varepsilon}\int_\pi^{\varphi_1} \kappa_0 d\varphi - \frac{i}{\varepsilon}\int_{\varphi_1}^{\varphi_b} \kappa_0 d\varphi} e^{-\frac{i}{\varepsilon}\int_{\varphi_b}^{\varphi} \kappa_0 d\varphi}(\Psi_+^1 + o(1)).$$

The curve $\gamma$ joining $\varphi_1$ to $\varphi_b$ is a Stokes line for $\kappa_0$ (and $\kappa_1$), hence

$$|e^{-\frac{i}{\varepsilon}\int_{\varphi_1}^{\varphi_b} \kappa_0 d\varphi}| = 1.$$

Moreover we notice that for $\varphi \in [\varphi_1, \pi]$, $\kappa_0(\varphi) \in \mathbb{R}$.
Now as $\varphi$ has to be in an admissible sub-domain of $K_0$, we know that $\operatorname{Re}(\varphi - \varphi_b) \geq \delta$; hence, we get that

$$(5.29) \qquad |w(f_-^0, f_+^1)| \leq C|e^{-\frac{i}{\varepsilon}\int_0^\pi \kappa_0 d\varphi}|e^{C\delta/\varepsilon} = C|e^{-\frac{i}{\varepsilon}\int_0^{\varphi_1} \kappa_0 d\varphi}|e^{C\delta/\varepsilon}.$$

This was done on the segment $[\varphi_b, \varphi_+]$. By $\varepsilon$-periodicity, it will hold everywhere on the line $\{\operatorname{Im}(\varphi) = \operatorname{Im}(\varphi_b)\}$. The estimate obtained above is valid as long as the point $\varphi_-$ stays in the admissible sub-domain of $K_0$ and if the interval $[\varphi_-, \varphi_+]$ stays off some vicinity of the branching points. Hence, we get that, for $-\Delta \leq \operatorname{Im} \varphi \leq -\delta$ ($\Delta, \delta > 0$), we have

$$(5.30) \qquad |d_1(\varphi)| \leq C|e^{-\frac{i}{\varepsilon}\int_0^{\varphi_1} \kappa_0 d\varphi}|e^{C\delta/\varepsilon}.$$

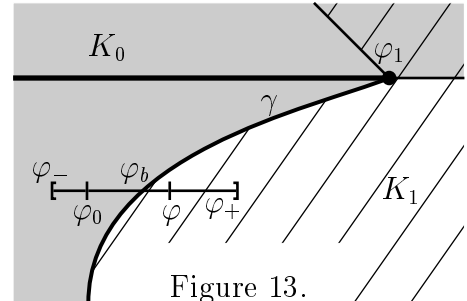

Figure 13.

Region (II). In this case, by (5.19), for $\varphi \in K_0 \cap K_1$, we have
$$f_-^0 = e^{-\frac{i}{\varepsilon}\int_0^\pi \kappa_0 d\varphi + \int_0^\pi \omega_-^0 d\varphi} e^{-\frac{i}{\varepsilon}\int_\pi^\varphi \kappa_0 d\varphi}(\Psi_-^1(u,\varphi) + o(1)).$$

Comparing this with the formula
$$f_+^1 = e^{\frac{i}{\varepsilon}\int_\pi^\varphi \kappa_0 d\varphi}(\Psi_+^1(u,\varphi) + o(1)),$$

one easily obtains
$$w(f_-^0, f_+^1) = -e^{-\frac{i}{\varepsilon}\int_0^\pi \kappa_0 d\varphi} e^{\int_0^\pi \omega_-^0 d\varphi}(w(f_+^1, f_-^1) + o(1)).$$

Hence, as long as we stay away from the turning points $\varphi_1$ and $\varphi_2$, in region (II), we get that

(5.31) $\quad d_1(\varphi) = e^{-\frac{i}{\varepsilon}\int_0^\pi \kappa_0 d\varphi} e^{\int_0^\pi \omega_-^0 d\varphi}(1 + o(1)), \ \delta \leq \operatorname{Im}\varphi \leq \operatorname{Im}\varphi_2 - \delta \ (\delta > 0)$

Here the integrals are taken along a curve in $K_0 \cap K_1$.

Region (III). In this region, we will satisfy ourselves with an estimate on $w(f_-^0, f_+^1)$. As $\operatorname{Im}\kappa_0 > 0$ in this sector, by Lemma 5.1, we can continue the asymptotics given by Theorem 3.1 for $f_-^0$ from $[\varphi_-,\varphi_0]$ to the whole interval $[\varphi_-,\varphi_+]$ (see figure 14). The consistent basis needed to apply Lemma 5.1 is obtained by using Theorem 3.1 for the canonical domain $\tilde{K}_0$. This domain comprises a neigh- 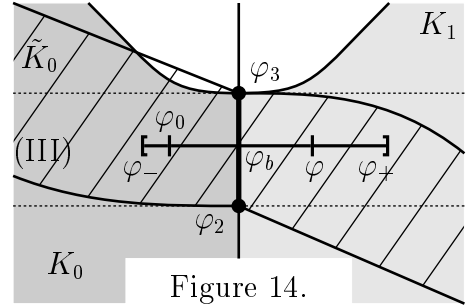

Figure 14.

borhood of the common boundary of $K_0$ and $K_1$ that is the segment $(\varphi_2, \varphi_3)$ (see Fig. 11, 12 and 14). On the segment $[\varphi_b, \varphi_+]$, we have
$$f_-^0 = e^{-\frac{i}{\varepsilon}\int_0^{\varphi_2}\kappa_0 d\varphi - \frac{i}{\varepsilon}\int_{\varphi_2}^{\varphi_b}\kappa_0 d\varphi} e^{-\frac{i}{\varepsilon}\int_{\varphi_b}^\varphi \kappa_0 d\varphi}(\Psi_-^0 + o(1))$$

Hence, for $\varphi \in [\varphi_-, \varphi_+]$, we get
$$|f_-^0| \leq C\left|e^{-\frac{i}{\varepsilon}\int_0^{\varphi_2}\kappa_0 d\varphi - \frac{i}{\varepsilon}\int_{\varphi_2}^\varphi \kappa_0 d\varphi}\right|$$

We also have
$$|f_+^1| \leq C\left|e^{\frac{i}{\varepsilon}\int_\pi^{\varphi_1}\kappa_0 d\varphi - \frac{i}{\varepsilon}\int_{\varphi_1}^\varphi \kappa_0 d\varphi} e^{\frac{i}{\varepsilon}(\varphi - \varphi_2)2\pi}\right|.$$

To get this estimate, we have used (5.23). The last two estimates imply
$$|d_1(\varphi)| \leq C\left|e^{-\frac{i}{\varepsilon}\int_0^{\varphi_2}\kappa_0 d\varphi}\right| \cdot \left|e^{-\frac{i}{\varepsilon}\int_{\varphi_2}^\pi \kappa_0 d\varphi}\right| \cdot \left|e^{-2\frac{i}{\varepsilon}\int_{\varphi_2}^\varphi (\kappa_0 - \pi)d\varphi}\right|.$$

Now we note that
- for $\varphi_1 \leq \varphi \leq \pi$, $\kappa_0 \in \mathbb{R}$.
- for $\varphi \in [\varphi_2, \varphi_3]$, $\kappa_0 = \pi + it$ where $t \geq 0$.

Assuming that $|\varphi - \varphi_b|$ is of order $\delta$, we obtain

(5.32) $\quad |d_1(\varphi)| \leq C|e^{-\frac{i}{\varepsilon}\int_0^{\varphi_1}\kappa_0 d\varphi}|e^{C\delta/\varepsilon}.$

Since $d_1$ is $\varepsilon$-periodic, this estimate is valid in the strip $\{\operatorname{Im}\varphi_2 + \delta \leq \operatorname{Im}\varphi \leq \operatorname{Im}\varphi_3 - \delta\}$ ($\delta > 0$).

Conclusion: as $d_1(\varphi)$ is $\varepsilon$-periodic in $\varphi$ and analytic in a strip $\{-\Delta \leq \operatorname{Im}\varphi \leq \operatorname{Im}\varphi_3 - \delta\}$, we can expand $d_1$ in a Fourier series with exponentially decreasing coefficients

$$d_1(\varphi) = \sum_{n \in \mathbb{Z}} \delta_n e^{2i\pi n \frac{\varphi}{\varepsilon}}$$

where, for any $\varphi_0 \in \{-\Delta \leq \operatorname{Im}\varphi \leq \operatorname{Im}\varphi_3 - \delta\}$,

$$\delta_n = \frac{1}{\varepsilon} \int_{\varphi_0}^{\varphi_0+\varepsilon} d_1(\varphi) e^{-2i\pi n \frac{\varphi}{\varepsilon}} d\varphi.$$

Using estimates (5.30) and (5.32) in the regions (I) and (III), we get that, for $n \neq 0$, we have

$$|\delta_n| \leq C |e^{-\frac{i}{\varepsilon}\int_0^\pi \kappa_0 d\varphi}| e^{-2\pi |n| Y_0/\varepsilon} e^{C\delta/\varepsilon}.$$

where $Y_0 = \min(\Delta, \operatorname{Im}\varphi_3 - \delta, Y)$ and $\delta$ is chosen such that $C\delta < 2\pi(Y - Y_0)$.
Using estimate (5.31) in region (II), we obtain

$$\delta_0 = e^{-\frac{i}{\varepsilon}\int_0^\pi \kappa_0 d\varphi} e^{\int_0^\pi \omega_-^0 d\varphi}(1 + o(1)).$$

Then, uniformly for $\varphi$ in the strip $\{|\operatorname{Im} z| \leq Y\}$, we have

$$d_1(\varphi) = e^{-\frac{i}{\varepsilon}\int_0^\pi \kappa_0 d\varphi + \int_0^\pi \omega_-^0 d\varphi}(1 + o(1))$$

5.3.6. *The asymptotics of $b_1(\varphi)$.* We will only estimate $w(f_+^0, f_+^1)$ when $\varepsilon \to 0$.

Region (I). Let $\gamma$ be the Stokes line beginning at $\varphi = \varphi_1$ and bordering $K_0$ and $K_1$. As $\operatorname{Im}\kappa_0 > 0$, Lemma 5.1 will only give us an estimation on $f_+^0$ when we cross $\gamma$ along $[\varphi_-, \varphi_+]$. This estimation is

$$|f_+^0| \leq C |e^{\frac{i}{\varepsilon}\int_0^{\varphi_1} \kappa_0 d\varphi + \frac{i}{\varepsilon}\int_{\varphi_1}^{\varphi_b} \kappa_0 d\varphi + \frac{i}{\varepsilon}\int_{\varphi_b}^{\varphi_0} \kappa_0 d\varphi}| e^{\frac{1}{\varepsilon}\int_{\varphi_0}^{\varphi} |\operatorname{Im}\kappa_0| d\varphi}.$$

As on $\gamma$ and on $[\varphi_-, \varphi_+]$, $\kappa_0 = -\kappa_1$, we can estimate $f_+^1$ by

$$|f_+^1| \leq C |e^{\frac{i}{\varepsilon}\int_\pi^{\varphi_1} \kappa_0 d\varphi - \frac{i}{\varepsilon}\int_{\varphi_1}^{\varphi_b} \kappa_0 d\varphi - \frac{i}{\varepsilon}\int_{\varphi_b}^{\varphi_0} \kappa_0 d\varphi}| |e^{-\frac{i}{\varepsilon}\int_{\varphi_0}^{\varphi} \kappa_0 d\varphi}|.$$

These two estimates imply

$$|w(f_+^0, f_+^1)| \leq C |e^{\frac{i}{\varepsilon}\int_0^{\varphi_1} \kappa_0 d\varphi}| |e^{\frac{2}{\varepsilon}\int_{\varphi_0}^{\varphi} |\operatorname{Im}\kappa_0| d\varphi}|.$$

for $\varphi$ on $[\varphi_b, \varphi_+]$ and inside the admissible domain for $K_1$. Now the only restriction we have on $\varphi_0$ is that it has to be in the admissible domain for $K_0$; hence we know that we can choose $|\varphi_0 - \varphi| \sim \delta$ where $\delta$ is the distance that defines the admissible domains (see (3.5)). It may be chosen as small as needed. So, for $-\Delta \leq \operatorname{Im}\varphi \leq -\delta$ ($\Delta, \delta > 0$), we have

(5.33) $$|b_1(\varphi)| \leq Ce^{C\frac{\delta}{\varepsilon}} |e^{\frac{i}{\varepsilon}\int_0^\pi \kappa_0 d\varphi}|.$$

Region (II). As we only estimate $B_1(\varphi)$, we don't need to compute it in the region.

Region (III). As $\operatorname{Im}(\kappa_0) > 0$, Lemma 5.1 will only give us an estimation on $f_+^0$ in $K_1$. The point $\varphi_b$ being as on Fig. 14 the intersection of $[\varphi_2, \varphi_3]$ and $[\varphi_-, \varphi_+]$, we get

$$|f_+^0| \leq C |e^{\frac{i}{\varepsilon}\int_0^\pi \kappa_0 d\varphi + \frac{i}{\varepsilon}\int_\pi^{\varphi_2} \kappa_0 d\varphi + \frac{i}{\varepsilon}\int_{\varphi_2}^{\varphi_b} \kappa_0 d\varphi}| e^{\frac{i}{\varepsilon}\int_{\varphi_b}^{\varphi_0} \kappa_0 d\varphi} e^{\frac{1}{\varepsilon}\int_{\varphi_0}^{\varphi} |\operatorname{Im}\kappa_0| d\varphi}$$

for $\varphi$ being in an $\delta$-vicinity of $\varphi_0$. Similarly, we have
$$|f_+^1| \leq C|e^{\frac{i}{\varepsilon}\int_\pi^{\varphi_2}\kappa_1 d\varphi + \frac{i}{\varepsilon}\int_{\varphi_2}^{\varphi_b}\kappa_1 d\varphi + \frac{i}{\varepsilon}\int_{\varphi_b}^{\varphi_0}\kappa_1 d\varphi}|e^{\frac{i}{\varepsilon}\int_{\varphi_0}^\varphi \kappa_1 d\varphi}$$
under the same conditions on $\varphi$.

We know that along $[\varphi_2, \varphi_3]$ and on $[\varphi_0, \varphi_+]$, we have $\kappa_0 = 2\pi - \kappa_1$; hence
$$|w(f_+^0, f_+^1)| \leq C|e^{\frac{i}{\varepsilon}\int_0^\pi \kappa_0 d\varphi}||e^{2\frac{i}{\varepsilon}\int_\pi^{\varphi_2}\kappa_0 d\varphi + 2\pi\frac{i}{\varepsilon}\int_{\varphi_2}^{\varphi_b}d\varphi + 2\pi\frac{i}{\varepsilon}\int_{\varphi_b}^{\varphi_0}d\varphi}||e^{\frac{2}{\varepsilon}\int_{\varphi_0}^\varphi |\operatorname{Im}\kappa_0|d\varphi}|$$
as $\varphi - \varphi_0 \geq \delta$. We now use the facts that $\kappa_0$ is real positive along $[\pi, \varphi_2]$ and that $\operatorname{Im}(\varphi_b - \varphi_2) > 0$. Hence, in (III), we have
$$|b_1(\varphi)| \leq Ce^{C\frac{\delta}{\varepsilon}}|e^{\frac{i}{\varepsilon}\int_0^\pi \kappa_0 d\varphi}|e^{-\eta/\varepsilon}$$
where $\eta$ is a fixed positive constant (independent of both $\delta$ and $\varepsilon$) as $\operatorname{Im}(\int_\pi^{\varphi_2}\kappa_0 d\varphi + \int_{\varphi_2}^{\varphi_b}d\varphi) > 0$ Hence, for $\delta$ small enough,
$$b_1(\varphi) = e^{\frac{i}{\varepsilon}\int_0^{\varphi_1}\kappa_0 d\varphi}e^{-c/\varepsilon}.$$

Conclusion: using the same method as before, in the strip $\{|\operatorname{Im}\varphi| \leq Y_0\}$, we get that
$$b_1(\varphi) = e^{\frac{i}{\varepsilon}\int_0^{\varphi_1}\kappa_0 d\varphi} \cdot o(1).$$

5.3.7. *The asymptotics of $a_1(\varphi)$.* We compute the asymptotics of $w(f_+^0, f_-^1)$ when $\varepsilon \to 0$. The computations are done in the same way as for $b_1$ and $d_1$.

Region (I). We get that,
$$|a_1(\varphi)| \leq Ce^{C\frac{\delta}{\varepsilon}}|e^{\frac{i}{\varepsilon}\int_0^\pi \kappa_0 d\varphi}|,$$
for $-\Delta \leq \operatorname{Im}\varphi \leq -\delta$ ($\Delta, \delta > 0$).

Region (II). We get that
$$a_1(\varphi) = e^{\frac{i}{\varepsilon}\int_0^\pi \kappa_0 d\varphi}e^{\int_0^\pi \omega_+^0 d\varphi}(1 + o(1)),$$
for $0 < \delta \leq \operatorname{Im}\varphi \leq \operatorname{Im}\varphi_2 - \delta$.

Region (III). We get that,
$$|a_1(\varphi)| \leq Ce^{2C\frac{\delta}{\varepsilon}}|e^{\frac{i}{\varepsilon}\int_0^\pi \kappa_0 d\varphi}|.$$
for $\operatorname{Im}\varphi_2 + \delta \leq \operatorname{Im}\varphi \leq \operatorname{Im}\varphi_3 - \delta$ ($\delta > 0$).

Conclusion: estimating the Fourier coefficients of $a_1$, as above, we get that, uniformly for $\varphi$ in the strip $\{|\operatorname{Im}z| \leq Y\}$,
$$a_1(\varphi) = e^{\frac{i}{\varepsilon}\int_0^\pi \kappa_0 d\varphi + \int_0^\pi \omega_+^0 d\varphi}(1 + o(1)).$$

5.3.8. *The asymptotics of $c_1(\varphi)$.* We compute the asymptotics of $w(f_-^0, f_-^1)$ when $\varepsilon \to 0$.

Region (I). As $\operatorname{Im}(\kappa_0) > 0$, we can use Lemma 5.1 to continue the asymptotic expansion for $f_-^0$ given by Theorem 3.1 in $K_0$ across $\partial K_0$ along a horizontal line (see Fig. 13). This will allow us to get a precise asymptotic of $c_1$. We use formulae (5.20) – (5.22). The asymptotic expansion for $f_-^0$ on $[\varphi_b, \varphi_+]$ is
$$f_-^0(u, \varphi) = e^{-\frac{i}{\varepsilon}\int_0^\varphi \kappa_0 d\varphi} \cdot \left(\sqrt{k_E'}e^{\int_0^\varphi \omega_-^0 d\varphi}(\Psi_-^0 + o(1))\right).$$

Here we integrate along paths going from 0 to $\varphi$ staying below the point $\varphi_1$. We use (5.20) and compute

$$\text{(5.34)} \quad \int_0^\varphi \kappa_0 d\varphi = \int_0^{\varphi_1} \kappa_0 d\varphi - \int_{\varphi_1}^\varphi \kappa_1 d\varphi$$
$$= \int_0^{\varphi_1} \kappa_0 d\varphi - \int_{\varphi_1}^\pi \kappa_1 d\varphi - \int_\pi^\varphi \kappa_1 d\varphi.$$

In the first two integrals, we integrate along the real axis; in the last integral, the contour is located in $K_1$. Similarly

$$\int_0^\varphi \omega_-^0 d\varphi = \int_0^{\varphi_1} \omega_-^0 d\varphi + \int_{\varphi_1}^\varphi \omega_+^1 d\varphi,$$

where the first integral is taken along the real line; the integration contour in the second integral lies in $K_1$. Finally, we get

$$\text{(5.35)} \quad \int_0^\varphi \omega_-^0 d\varphi = \int_0^{\varphi_1} \omega_-^0 d\varphi - \int_{\varphi_1}^\pi \omega_+^1 d\varphi + \int_\pi^\varphi \omega_+^1 d\varphi,$$

We also note that, on the different sides of $\gamma$ in region (I) (see Fig. 13), we have

$$\text{(5.36)} \quad \sqrt{k_E'}|_{\varphi+0} = i\sqrt{k_E'}|_{\varphi-0}$$

Combining (5.34) – (5.36), we get

$$f_-^0(u,\varphi) = -i e^{-\frac{i}{\varepsilon}\int_0^{\varphi_1} \kappa_0 d\varphi + \frac{i}{\varepsilon}\int_{\varphi_1}^\pi \kappa_0 d\varphi + \int_0^{\varphi_1} \omega_-^0 - \int_\pi^{\varphi_1} \omega_+^0} e^{\frac{i}{\varepsilon}\int_\pi^\varphi \kappa_1 d\varphi}(\Psi_+^1(u,\varphi) + o(1)).$$

On the other hand, we know that

$$f_-^1(u,\varphi) = e^{-\frac{i}{\varepsilon}\int_\pi^\varphi \kappa_1 d\varphi}(\Psi_-^1(u,\varphi) + o(1)).$$

Hence we easily obtain

$$w(f_-^0, f_-^1) = -i \cdot e^{-\frac{i}{\varepsilon}\int_0^\pi \kappa_0 d\varphi} e^{2\frac{i}{\varepsilon}\int_{\varphi_1}^\pi \kappa_0 d\varphi} e^{\int_0^{\varphi_1} \omega_-^0 d\varphi - \int_\pi^{\varphi_1} \omega_+^0 d\varphi}(w_1 + o(1)),$$

that is

$$\text{(5.37)} \quad c_1(\varphi) = -i \cdot e^{-\frac{i}{\varepsilon}\int_0^\pi \kappa_0 d\varphi} e^{2\frac{i}{\varepsilon}\int_{\varphi_1}^\pi \kappa_0 d\varphi} e^{\int_0^{\varphi_1} \omega_-^0 d\varphi - \int_\pi^{\varphi_1} \omega_+^0 d\varphi}(1 + o(1))$$

where the integrals are taken along the real line.
This holds in the strip $\{-\Delta \leq \operatorname{Im}\varphi \leq -\delta\}$ ($\Delta > \delta > 0$).

Region (II). As for $b_1$, we will not need any information on $c_1$ in this region.

Region (III). As $\operatorname{Im}(\kappa_0) > 0$, we will again be able to use Lemma 5.1 to compute a precise asymptotic for $c_1$. So that, for $\varphi \in [\varphi_-, \varphi_+]$, we know that

$$\text{(5.38)} \quad f_-^0 = e^{-\frac{i}{\varepsilon}\int_0^\varphi \kappa_0 d\varphi}\left(\sqrt{k_E'} e^{\int_0^\varphi \omega_-^0 d\varphi}\psi_-^0 + o(1)\right)$$

where $[\varphi_-, \varphi_+]$ is defined in Fig. 14. In (5.38), the integration contour first goes from 0 to $\varphi_-$ in $K_0$, then to $\varphi$ along $[\varphi_-, \varphi_+]$. Using (5.24) – (5.25), we rewrite the integrals in (5.38) to get

$$\int_0^\varphi \kappa_0 d\varphi = \int_0^\pi \kappa_0 d\varphi + 2\int_\pi^{\varphi_2} \kappa_0 d\varphi + 2\pi(\varphi - \varphi_2) - \int_\pi^\varphi \kappa_1 d\varphi,$$

where, in the first integral, we integrate along a path in $K_0 \cap K_1$, and in the next two, along a path in $K_1$.

The same way we get

$$\int_0^\varphi \omega_-^0 d\varphi = \int_0^\pi \omega_-^0 d\varphi + \int_\pi^{\varphi_2} (\omega_-^0 - \omega_+^0) d\varphi + \int_\pi^\varphi \omega_+^0 \varphi.$$

We also note that, for $\varphi \in [\varphi_2, \varphi_3]$, we have

(5.39) $$\sqrt{k'_E}_{|\varphi+0} = -i\sqrt{k'_E}_{|\varphi-0}.$$

Combining these formulae, we obtain

$$f_-^0 = i e^{-\frac{i}{\varepsilon}\int_0^\pi \kappa_0 d\varphi - 2\frac{i}{\varepsilon}\int_\pi^{\varphi_2}(\kappa_0-\pi)d\varphi} e^{-2\pi\frac{i}{\varepsilon}(\varphi-\pi)} e^{\int_0^{\varphi_2} \omega_-^0 d\varphi - \int_\pi^{\varphi_2} \omega_+^0 d\varphi} e^{\frac{i}{\varepsilon}\int_\pi^\varphi \kappa_1 d\varphi}(\Psi_+^1 + o(1)).$$

Using this in conjunction with

$$f_-^1 = e^{-\frac{i}{\varepsilon}\int_\pi^\varphi \kappa_1 d\varphi}(\Psi_-^1 + o(1)),$$

we get

(5.40) $$c_1(\varphi) = i e^{-\frac{i}{\varepsilon}\int_0^\pi \kappa_0 d\varphi - 2\frac{i}{\varepsilon}\int_\pi^{\varphi_2}(\kappa_0-\pi)d\varphi} e^{-2\pi\frac{i}{\varepsilon}(\varphi-\pi)} e^{\int_0^{\varphi_2} \omega_-^0 d\varphi - \int_\pi^{\varphi_2} \omega_+^0 d\varphi}(1+o(1)).$$

This asymptotic is valid for $\operatorname{Im}\varphi_2 + \delta \leq \varphi \leq \operatorname{Im}\varphi_3 - \delta$ ($\delta > 0$).

Conclusion: now to get global information on $c_1$, we use the same method as above; let

$$c_1(\varphi) = \sum_{n\in\mathbb{Z}} \gamma_n e^{2i\pi n \frac{\varphi-\pi}{\varepsilon}}.$$

Then, by the asymptotic (5.37), for $n > 0$, we have

$$|\gamma_n| \leq C |e^{-\frac{i}{\varepsilon}\int_0^\pi \kappa_0 d\varphi}| e^{-\frac{2\pi}{\varepsilon}|n|Y_0}.$$

The same formula gives

$$\gamma_0 = i e^{-\frac{i}{\varepsilon}\int_0^\pi \kappa_0 d\varphi + 2\frac{i}{\varepsilon}\int_{\varphi_1}^\pi \kappa_0 d\varphi} e^{\int_0^{\varphi_1} \omega_-^0 d\varphi - \int_\pi^{\varphi_1} \omega_+^0 d\varphi}(1+o(1)).$$

We use the asymptotic given by region (III) to estimate $(\gamma_n)$ for $n < 0$. For some $\varphi_0$ in (III), we compute

$$\gamma_{-1} = \frac{1}{\varepsilon}\int_{\varphi_0}^{\varphi_0+\varepsilon} c_1(\varphi) e^{2i\pi \frac{\pi-\varphi}{\varepsilon}} d\varphi$$

$$= i e^{-\frac{i}{\varepsilon}\int_0^\pi \kappa_0 d\varphi - 2\frac{i}{\varepsilon}\int_\pi^{\varphi_2}(\kappa_0-\pi)d\varphi} e^{\int_0^{\varphi_2} \omega_-^0 d\varphi - \int_\pi^{\varphi_2} \omega_+^0 d\varphi}(1+o(1))$$

Now we recall that the segment $[\varphi_2, \varphi_3]$ is a Stokes line for $\kappa_0$ and $\kappa_0(\varphi_2) = \pi$; hence $\operatorname{Im}\int_\pi^{\varphi_2}(\kappa_0 - \pi)d\varphi = 0$. For $n < -1$, we use (5.40) to estimate

$$|\gamma_n| = \left|\frac{1}{\varepsilon}\int_{\varphi_0}^{\varphi_0+\varepsilon} c_1(\varphi) e^{2i\pi \frac{\pi-\varphi}{\varepsilon}} d\varphi\right|$$

$$\leq C \left|e^{-\frac{i}{\varepsilon}\int_0^\pi \kappa_0 d\varphi} e^{-2\frac{i}{\varepsilon}\int_\pi^{\varphi_2}(\kappa_0-\pi)d\varphi} e^{\int_0^{\varphi_2} \omega_-^0 d\varphi - \int_\pi^{\varphi_2} \omega_+^0 d\varphi}\right| \frac{1}{\varepsilon}\left|\int_{\varphi_0}^{\varphi_0+\varepsilon} e^{2i(n-1)\pi\frac{\varphi}{\varepsilon}}\right|$$

$$\leq C |e^{-\frac{i}{\varepsilon}\int_0^\pi \kappa_0 d\varphi}| e^{-\frac{2\pi}{\varepsilon}|n+1|Y_0}.$$

Putting all this together, we get that, for $\varphi$ in the strip $\{|\mathrm{Im} z| \leq Y_0\}$, we have

$$c_1(\varphi) = -ie^{-\frac{i}{\varepsilon}\int_0^\pi \kappa_0 d\varphi}\left(e^{2\frac{i}{\varepsilon}\int_{\varphi_1}^\pi \kappa_0 d\varphi}e^{\int_0^{\varphi_1} \omega_-^0 d\varphi - \int_\pi^{\varphi_1} \omega_+^0 d\varphi}(1+o(1))\right.$$
$$\left.-e^{-2\frac{i}{\varepsilon}\int_\pi^{\varphi_2}(\kappa_0-\pi)d\varphi}e^{\int_0^{\varphi_2}\omega_-^0 d\varphi - \int_\pi^{\varphi_2}\omega_+^0 d\varphi}e^{-2\pi\frac{i}{\varepsilon}(\varphi-\pi)}(1+o(1))\right)$$

where $o(1)$ is uniform in the strip $\{|\mathrm{Im} z| \leq Y_0\}$. This ends the proof of Proposition 5.1.

5.4. **The asymptotics of $T_2$**. As we have already written, we will not do the detailed computations for $T_2$. We will only describe the geometry that is slightly different. Let us write

$$T_2(\varphi) = \begin{pmatrix} a_2(\varphi) & b_2(\varphi) \\ c_2(\varphi) & d_2(\varphi) \end{pmatrix},$$

then

(5.41) $$a_2(\varphi) = \frac{1}{w_2}w(f_+^1, f_-^2), \quad b_2(\varphi) = -\frac{1}{w_2}w(f_+^1, f_+^2),$$

(5.42) $$c_2(\varphi) = \frac{1}{w_2}w(f_-^1, f_-^2), \quad d_2(\varphi) = -\frac{1}{w_2}w(f_-^1, f_+^2),$$

(5.43) $$w_2 = w(f_+^2, f_-^2) = f_-^2\frac{d}{du}f_+^2 - f_+^2\frac{d}{du}f_-^2.$$

$T_2$ shares many of the properties of $T_1$ (though not in the same domain); in particular:

- its coefficients are $\varepsilon$-periodic,
- they are analytic in some strip.

Now due to the geometry of $K_1$ and $K_2$ (see Fig. 15), we will only be able to control these coefficients in some strip of the form

$$-\mathrm{Im}\,\varphi_3 + \delta \leq \mathrm{Im}\,\varphi \leq \Delta, \quad \Delta > \delta > 0.$$

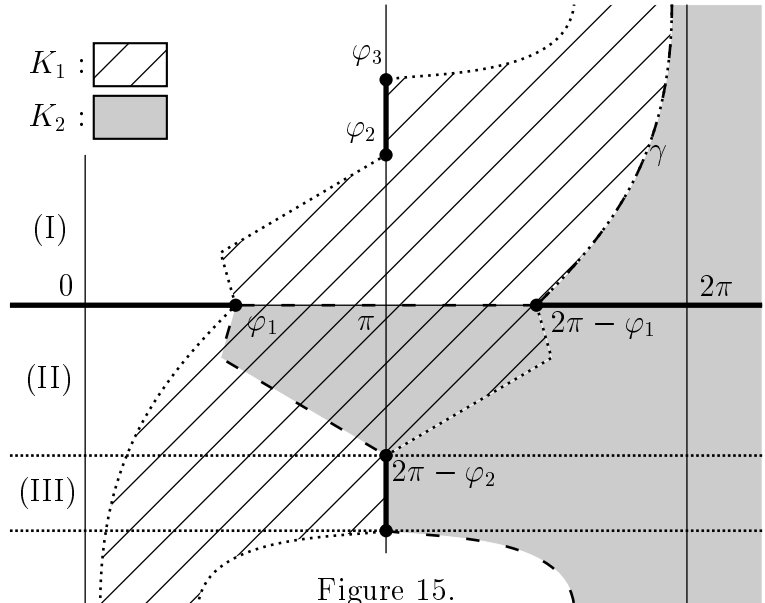

Figure 15.

Again we divide the strip $\{-\operatorname{Im}\varphi_3 + \delta \leq \operatorname{Im}\varphi \leq \Delta\}$ into three smaller strips called respectively (I), (II) and (III) (see Fig. 15).

Let us briefly describe the main objects of our construction in each of these strips.

In the strip (II), we have a non-empty intersection for $K_1$ and $K_2$. In the intersection, we have

$$\kappa_1 = \kappa_2, \quad \psi_\pm^1 = \psi_\pm^2 \quad \text{and} \quad \omega_\pm^1 = \omega_\pm^2.$$

Consider strip (I). Denote by $\gamma$ the common boundary of $K_2$ and $K_1$ in this strip (see Fig. 16). It is one of the Stokes lines beginning at $\varphi = 2\pi - \varphi_1$. On $\gamma$, we get

$$\kappa_1(\varphi - 0) = -\kappa_2(\varphi + 0)$$
$$\psi_\pm^1(u, \varphi - 0) = \psi_\mp^2(u, \varphi + 0)$$
$$\omega_\pm^1(\varphi - 0) = \omega_\mp^2(\varphi + 0).$$

These formulae immediately follow from equations (5.20) – (5.22) as, for $\varphi \in K_0$, we have

(5.44) $\qquad \kappa_2(2\pi - \varphi) = \kappa_0(\varphi).$

We also notice that the conformal properties of $K_0$ (see subsection 4.1) and (5.20) show that, for $\varphi \in$(I) and $\operatorname{Im}\varphi < 0$, we have

$$-\operatorname{Im}\kappa_2(\varphi) = \operatorname{Im}\kappa_1(\varphi) < 0.$$

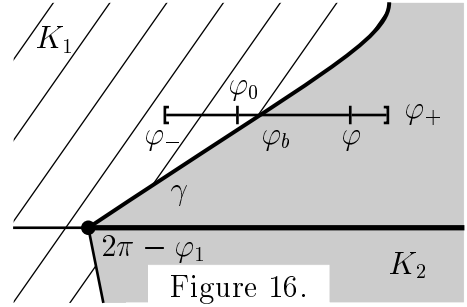

Figure 16.

In the strip (III), the common boundary of $K_2$ and $K_1$ consists in the interval $[2\pi - \varphi_2, 2\pi - \varphi_3]$. This interval is the Stokes line joining $2\pi - \varphi_2$ to $2\pi - \varphi_3$. Along $[2\pi - \varphi_2, 2\pi - \varphi_3]$, we have

$$\kappa_1(\varphi - 0) = 2\pi - \kappa_2(\varphi + 0)$$
$$\psi_\pm^1(u, \varphi - 0) = \psi_\mp^2(u, \varphi + 0)$$
$$\omega_\pm^1(\varphi - 0) = \omega_\mp^2(\varphi + 0).$$

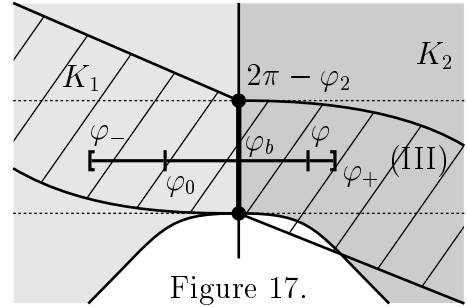

Figure 17.

As in region (I), for $-\operatorname{Im}\varphi_3 < \operatorname{Im}\varphi < -\operatorname{Im}\varphi_2$, we have

$$-\operatorname{Im}\kappa_2(\varphi) = \operatorname{Im}\kappa_1(\varphi) < 0.$$

To use Lemma 5.1 to compute $T_2$, we will use another auxiliary canonical domain, namely the symmetric of $\tilde{K}_0$ with respect to the point $\pi$. It comprises a neighborhood of $(2\pi - \varphi_2, 2\pi - \varphi_3)$. It is canonical with respect to $\tilde{\kappa}_2$ the branch of the quasi-momentum obtained from $\kappa_2$ by analytic continuation through the interval $[2\pi - \varphi_2, 2\pi - \varphi_3]$ by (5.44). How to use Lemma 5.1 is partly described in Fig. 16 and 17. We also point out that, in this case, (5.39) becomes

$$\sqrt{k'_E}_{|\varphi+0} = i\sqrt{k'_E}_{|\varphi-0}.$$

We will not detail the computations of the coefficients of the matrix $T_2$ any further.

**Remark 5.2.** Let us just notice that one could have used the symmetry between the domain $K_0 \cup K_1$ and $K_1 \cup K_2$, and between the associated consistent basis to compute $T_2$ from $T_1$.

5.5. **The computation of** $M$. Let $K := K_0 \cup K_1 \cup K_2$. This domain is represented on Fig. 18. Its natural cuts are the branches of Stokes lines shown on Fig. 18 by bold lines. On this domain, we define a unique branch of the complex quasi-momentum $\kappa$ by

$$\kappa := \kappa_0 = \kappa_1.$$

This branch is analytic on $K$. In the same way, on $K$, we define a unique analytic function $\omega$ by

$$\omega_\pm := \omega_\pm^0 = \omega_\pm^1.$$

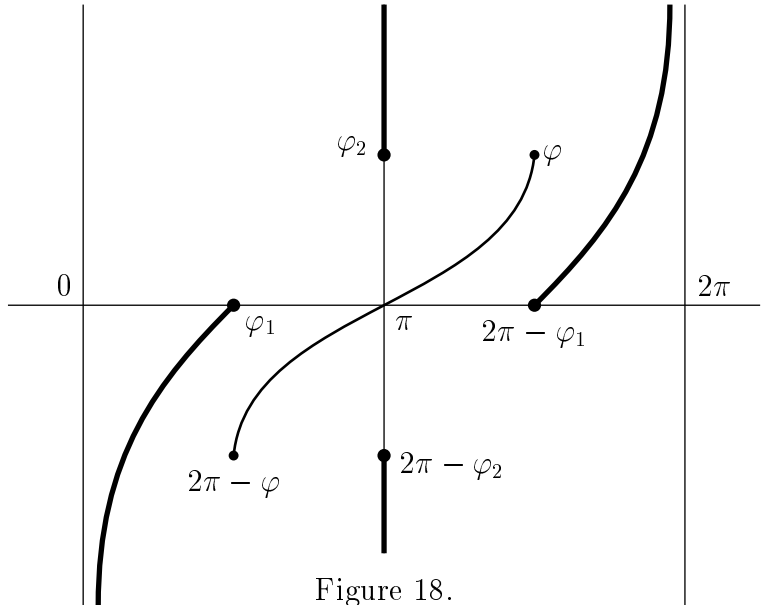

Figure 18.

These branches have the following symmetry properties:

(5.45) $\qquad \kappa(2\pi - \varphi) = \kappa(\varphi) \quad \text{and} \quad \omega_\pm(2\pi - \varphi) = -\omega_\pm(\varphi).$

The first equality is just (4.9); the second is obtained using (3.6) and the fact that, as $\cos(2\pi - \varphi) = \cos(\varphi)$, we have $p_\pm(x, 2\pi - \varphi) = p_\pm(x, \varphi)$.

Rewriting the asymptotics of the matrices $T_1$ and $T_2$ in terms of $\kappa$ and $\omega_\pm$ and using the symmetry property (5.45), one computes the product of these matrices. It is an elementary but long and tedious calculation that we will not reproduce here. This ends the proof of Theorem 5.1. $\qquad \square$

We end this section with proving that $\phi_2$ is real. Therefore we compute

$$\overline{\phi_2} = -\frac{i\varepsilon}{2\pi} \overline{\int_{\varphi_1}^{\varphi_2} (\omega_+ - \omega_-)} = \frac{i\varepsilon}{2\pi} \int_{\varphi_1}^{2\pi - \varphi_2} (\omega_+ - \omega_-)$$

$$= \frac{i\varepsilon}{2\pi} \int_{\varphi_1}^{\varphi_2} (\omega_+ - \omega_-),$$

as $\overline{\omega_+(\varphi)} = \omega_-(\overline{\varphi})$, $\overline{\varphi}_2 = 2\pi - \varphi_2$ and $\int_{\varphi_2}^{2\pi-\varphi_2} (\omega_+ - \omega_-) = 0$.

## 6. Spectral application

### 6.1. Reduction to a difference equation.
Here, we discuss the passage from the family of equations (1.1) to the difference equation (1.9) and justify the relation (1.11).

Consider the vector $\Psi$ composed of the consistent basis solutions of (1.1) as in (1.7). The definition of the monodromy matrix implies that

(6.1)
$$\begin{aligned}\Psi\left(x+2\pi n/\varepsilon, \phi + 2\pi n/\varepsilon\right) &= \\ &= M\left(\phi + 2\pi(n-1)/\varepsilon\right) M\left(\phi + 2\pi(n-2)/\varepsilon\right) \cdots \\ &\qquad \cdots M\left(\phi + 2\pi/\varepsilon\right) M\left(\phi\right) \Psi\left(x, \phi\right).\end{aligned}$$

Assume that $\chi$ is a fundamental solution of equation (1.9) defined by the monodromy matrix $M$ corresponding to $\Psi$. By definition, $\chi(\phi + h)\chi^{-1}(\phi) = M(\phi)$. Since $M(\phi + 2\pi l/\varepsilon) = M(\phi + lh)$ (for $l \in \mathbb{N}$), together with (6.1), this implies the relation
$$\Psi\left(x + 2\pi n/\varepsilon, \phi + 2\pi n/\varepsilon\right) = \chi\left(\phi + nh\right)\chi^{-1}(\phi)\Psi\left(x, \phi\right).$$
Since $\Psi$ is 1-periodic in $\phi$, it immediately leads to (1.11).

### 6.2. Bloch solutions of difference equations.

6.2.1. To prove Theorem 1.3 we have to construct some solutions of equation (1.9). So, we recall some basic ideas (see [8]) related to difference equations of the form

(6.2)
$$\psi\left(\phi + h\right) = M\left(\phi\right)\psi\left(\phi\right), \quad \phi \in \mathbb{R},$$

where $h$ is a positive number and $M(\phi) \in \mathrm{SL}\left(2, \mathbb{C}\right)$ is a given matrix function 1-periodic in $\phi$.

If $\psi$ is a fundamental solution (i.e. $\det\psi = 1$), then any other solution can be uniquely represented in the form
$$\chi\left(\phi\right) = \psi\left(\phi\right) P\left(\phi\right), \quad \phi \in \mathbb{R}.$$
where $P$ is an $h$-periodic matrix function, $P\left(\phi + h\right) = P\left(\phi\right)$. In particular, one can write

(6.3)
$$\psi(\phi + 1) = \psi\left(\phi\right) M_1^T\left(\phi\right).$$

This relation defines the *monodromy matrix* $M_1$ corresponding to the fundamental solution $\psi$. Clearly,

(6.4)
$$M_1(\phi + h) = M_1(\phi), \qquad \det M_1(\phi) \equiv 1, \ \forall \phi \in \mathbb{R}.$$

If the monodromy matrix corresponding to a fundamental solution is diagonal,
$$M_1(\phi) = \begin{pmatrix} \beta\left(\phi\right) & 0 \\ 0 & \beta^{-1}(\phi) \end{pmatrix},$$
we call $\psi$, a *Bloch solution*, and the coefficient $\beta$ is its *Floquet multiplier*.

We pick the determination of the function $\ln \beta(\phi)$ so that all its possible discontinuity jumps be smaller than $2\pi$. If the limit

(6.5) $$\theta(\phi) = \lim_{n\to\infty} \frac{1}{n} \sum_{l=0}^{n} \ln \beta(\phi + l)$$

exists and its convergence is uniform in $\phi$, we call this limit the *complex Lyapunov exponent* of the Bloch solution $\psi$.

We call a Bloch solution monotonous if $\operatorname{Re}\theta(\phi) \geq Const > 0$ for all $\phi$.

We finish this section by formulating a simple condition for the existence of the monotonous Bloch solutions [8].

Set

(6.6) $$\rho(\phi) = \frac{M_{12}(\phi)}{M_{12}(\phi - h)}, \quad v(\phi) = M_{11}(\phi) + \rho(\phi) M_{22}(\phi - h).$$

Let

$$\rho_- = \inf_{\phi \in \mathbb{R}} |\rho(\phi)|, \quad \rho_+ = \sup_{\phi \in \mathbb{R}} |\rho(\phi)|, \quad v_- = \inf_{\phi \in \mathbb{R}} |v(\phi)|.$$

For $f(\phi) \neq 0$ being a periodic continuous function on $\mathbb{R}$, we denote by $\operatorname{ind} f$ the integer equal to the increment of $\frac{\arg f}{2\pi}$ over one period. Clearly, $\operatorname{ind} f$ is independent of $\phi$. One has

**Proposition 6.1.** *Let $M$ be continuous in $\phi$, and let $M_{12} \neq 0$. If*

$$\rho_- > 0, \quad \rho_+ < \left(\frac{v_-}{2}\right)^2, \quad \text{and} \quad \operatorname{ind} v = \operatorname{ind} \rho = 0,$$

*then*

1. *equation (6.2) has a fundamental Bloch solution continuous in $\phi$;*
2. *the Lyapunov exponent satisfies the estimates*

$$\ln(v_-/2 + \sqrt{(v_-/2)^2 - \rho_+}) \leq h \operatorname{Re}(\theta(\phi)) \leq \ln(v_+ + v_-/2 - \sqrt{(v_-/2)^2 - \rho_+}).$$

6.3. **Bloch vectors for the family of equations (1.1).** In terms of Bloch solutions of (1.9), one can construct a very natural set of solutions of the family of equation (1.1). Let $\Psi$ be a vector composed of consistent basis solutions of (1.1), and let $M$ be the corresponding monodromy matrix. Suppose that $\chi$ is a Bloch solution of equation (1.9). Set

$$F(x, \phi) = \chi^{-1}(\phi) \Psi(x, \phi).$$

The components $F_1$ and $F_2$ of the vector $F$ satisfy (1.1). First, we note that $w(F_1, F_2) = \det \chi^{-1}(\phi) w(\Psi_1, \Psi_2)$, and thus,

(6.7) $$w(F_1, F_2) \equiv 1, \quad \forall \phi \in \mathbb{R}.$$

Written in terms of $F$, formula (1.11) takes the form

(6.8) $$F(x + 2\pi n/\varepsilon, \phi) = F(x, \phi - nh), \quad \forall x, \phi \in \mathbb{R}.$$

Finally, one can easily see that

(6.9) $$F(x, \phi + 1) = \begin{pmatrix} \beta^{-1}(\phi) & 0 \\ 0 & \beta(\phi) \end{pmatrix} F(x, \phi), \quad \forall \phi \in \mathbb{R},$$

where $\beta$ is the Floquet multiplier of the Bloch solution $\chi$.

We shall call any vector of solutions of (1.1) possessing the properties (6.7)–(6.9) a *Bloch vector* of this family of equation, and the coefficient $1/\beta$ from (6.9) its *Floquet multiplier*.

6.4. **Bloch vectors and the spectra of (1.1).** As it was done above for Bloch solutions of equation (6.2), we define the complex Lyapunov exponent of a Bloch vector of (1.1) in terms of its Floquet multiplier. If the Lyapunov exponent verifies the estimate of the form $\operatorname{Re}\theta(x) \geq Const > 0$ for all $x$, we call the Bloch vector *monotonous*.

Clearly, if there is a monotonous Bloch solution of (1.9), then there exists a monotonous Bloch vector of (1.1). We shall use the following elementary observation.

**Lemma 6.1.** *Assume that, for a given $E$, there exists a monotonous Bloch vector for (1.1), and that its components are locally bounded in $x$ and $\phi$. Then, for any $\phi \in \mathbb{R}$, $E$ does not belong to the spectrum of (1.1) in $L^2(\mathbb{R})$.*

*Proof.* For any given $\phi$, $F_{1,2}$ are solutions of (1.1). By (6.7), these solutions are linearly independent. Estimate the first of them. In view of (6.8) – (6.9),

$$|F_1(x + 2n\pi/\varepsilon, \phi)| = \prod_{l=0}^{N-1} |\beta(\phi + l)||F_1(x, \phi - nh + N)|.$$

Let $0 \leq x \leq 2\pi/\varepsilon$. Choose $N$ so that $0 \leq \phi - nh + N \leq 1$. Then

$$|F_1(x + 2n\pi/\varepsilon, \phi)| \leq C \exp\left(-\sum_{l=1}^{N} \ln|\beta(\phi + l)|\right).$$

Since the limit (6.5) exists and is uniform in $x$, and the Bloch vector is monotonous, then
$$|F_1(x + 2n\pi/\varepsilon, \phi)| \leq C\, e^{-n\delta},$$
for sufficiently large $n$ with some positive $\delta$. In the same way, one checks that
$$|F_2(x + 2n\pi/\varepsilon, \phi)| \leq C e^{n\delta}$$
for sufficiently large $-n$. The linear independence of the solutions $F_{1,2}$ and their estimates imply that the energy $E$ is in the resolvent set of the operator defined by the left hand side of (1.1) (see [10]).

6.5. **Exponential localization of the spectrum.** Now, we can immediately prove Theorem 1.3. Consider equation (6.2) with the matrix $M$ described in the Theorem 1.2. In this case,

$$\rho = 1 + o(1), \quad v(z) = \tfrac{2}{t}\left(\cos(\phi_1 + o(1)) - t_1(1 + o(1))\cos z\right)(1 + o(1))$$

with $z$ chosen as in Theorem 1.2. All the error terms are locally analytic in $E$. The assumptions of Proposition 6.1 are satisfied if

$$|\cos(\phi_1 + o(1))| \geq (t + t_1)(1 + o(1)).$$

This condition guarantees the existence of monotonous Bloch solutions of (1.9). Now, combining Proposition 6.1 and Lemma 6.1, we see that, for $E$ being in the interval described by (1.26), the spectrum of (1.1) can be situated only on the subintervals where

(6.10) $\qquad\qquad\qquad |\cos(\phi_1 + o(1))| \leq (t + t_1)(1 + o(1)).$

Remind that $\phi_1$, $t$, $t_1$ are analytic in $E$ in a vicinity of any point $E$ satisfying (1.26). We show that

(6.11)
$$\phi_1'(E) \neq 0,$$

for any $E$ satisfying (1.26). This function is defined by (5.1). It can be represented in the form
$$\phi_1 = \frac{1}{\varepsilon} \int_{\varphi_1}^{2\pi - \varphi_1} \kappa_1 \, d\varphi.$$

Here $\kappa_1$ is the branch of the complex momentum corresponding to the canonical domain $K_1$, see section 5. Recall that it is related to a branch of the Bloch quasi-momentum of (1.14) by the formula
$$\kappa_1(\varphi) = k\left(E - \cos\varphi\right).$$

So that, for $E$ satisfying (1.26),

(6.12)
$$\phi_1'(E) = \frac{1}{\varepsilon} \int_{\varphi_1}^{2\pi - \varphi_1} k'\left(E - \cos\varphi\right) d\varphi.$$

Omitting the elementary details, we only note that, since $k$ as a function of $\varphi$ has square root branching points at $\varphi = \varphi_1$ and $\varphi = 2\pi - \varphi_1$, this integral (6.12) converges. Moreover, clearly, up to its sign, it does not depend on the choice of the branch of $k$. For any $\varphi \in (\varphi_1, 2\pi - \varphi_1)$,
$$E_1 \leq E - \cos\varphi.$$

Moreover, under the condition (1.26),
$$E - \cos\varphi < E_2 - \delta, \quad \delta > 0, \quad \varphi \in \mathbb{R}$$

As $k'(E - \cos\varphi)$ is real and these two inequalities imply that $k'(E - \cos\varphi) \neq 0$ on the interval $(\varphi_1, 2\pi - \varphi_1)$, we get the desired result.

The proof of Theorem 1.3 follows from (6.10) – (6.11). We now establish the properties of the subintervals described by (6.10). Since $t$ and $t_1$ are exponentially small and $\phi_1' \neq 0$, the length of these intervals are exponentially small as well. They are situated in $o\left(\varepsilon\right)$-vicinities of the points $E^{(l)}$ defined by
$$\cos\phi_1(E^{(l)}) = 0.$$

Since $\phi_1' \neq 0$, then the distances between these points are of order $\varepsilon$, and near each of these points, there is only one of the subintervals containing the spectrum. Recalling the definitions of $t$ and $t_1$ (cf (5.1)), one sees that, up to the factor $1 + o\left(1\right)$, these functions are constant on the exponentially small intervals. This implies that the length of the subinterval situated near to $E^{(l)}$ is equal to $(t(E^{(l)}) + t_1(E^{(l)}))\,(1+o\,(1))$. This completes the proof. $\square$

(Alexander Fedotov) DEPARTEMENT OF MATHEMATICAL PHYSICS, ST PETERSBURG STATE UNIVERSITY, 1, ULIANOVSKAJA, 198904 ST PETERSBURG-PETRODVORETZ, RUSSIA
*E-mail address*: `fedotov@mph.niif.spb.su`

(Frédéric Klopp) DÉPARTEMENT DE MATHÉMATIQUE, INSTITUT GALILÉE, L.A.G.A., URA 742 C.N.R.S, UNIVERSITÉ DE PARIS-NORD, AVENUE J.-B. CLÉMENT, F-93430 VILLETANEUSE, FRANCE
*E-mail address*: `klopp@math.univ-paris13.fr`